\documentclass[draft,preprint,12pt,numbers,sort&compress]{elsarticle}
\usepackage{mathrsfs}
\usepackage{amssymb}
\usepackage{amsthm}
\usepackage{mathrsfs}
\usepackage[centertags]{amsmath}
\usepackage{amsfonts}

\usepackage{color}

\input amssym.def
\input amssym.tex

\date{}

\newtheorem{Theorem}{Theorem}[section]

\newtheorem{Lemma}{Lemma}[section]

\newcommand\R{\mbox{\bf R}}

\newcommand\SR{\mbox{\scriptsize\bf R}}

\newcommand{\definition}{{\lower .5ex
  \hbox{$\>\>\stackrel{\triangle}{=}\>\>$} }}

\begin{document}

\baselineskip=22pt
\thispagestyle{empty}

\mbox{}
\bigskip

\begin{center}
{\Large \bf The Cauchy problem for  two dimensional  generalized}\\[1ex]
{\Large \bf  Kadomtsev-Petviashvili-I equation  in anisotropic   Sobolev spaces}\\[1ex]

{Wei Yan\footnote{Email:yanwei19821115@sina.cn}$^{a,d}$,\quad Yongsheng Li
\footnote{Email:yshli@scut.edu.cn}$^b$,\quad Jianhua  Huang
\footnote{ Email: jhhuang32@nudt.edu.cn}$^c$,\quad Jinqiao  Duan\footnote{ Email: duan@iit.edu}$^{d}$}\\[1ex]

{$^a$School of Mathematics and Information Science and  Henan Engineering Laboratory for
 Big Data Statistical Analysis and Optimal Control, Henan  Normal University,}\\
{Xinxiang, Henan 453007, P. R. China}\\[2ex]

{$^b$School of Mathematics, South  China  University of  Technology,}\\
{Guangzhou, Guangdong 510640,  China}\\[1ex]

{$^c$College of Science, National University of Defense  Technology,}\\
{ Changsha, Hunan 410073,  China}\\[1ex]

{$^d$Department of Applied Mathematics, Illinois Institute of Technology,}\\[1ex]
{Chicago, IL 60616, USA}\\[1ex]

\end{center}
\noindent{\bf Abstract.} The goal of  this paper is three-fold. Firstly,
 we prove that the Cauchy problem
  for generalized KP-I equation
\begin{eqnarray*}
 u_{t}+|D_{x}|^{\alpha}\partial_{x}u+\partial_{x}^{-1}\partial_{y}^{2}u+\frac{1}{2}\partial_{x}(u^{2})=0,\alpha\geq4
 \end{eqnarray*}
is locally well-posed  in the anisotropic Sobolev spaces $H^{s_{1},\>s_{2}}(\R^{2})$
  with $s_{1}>-\frac{\alpha-1}{4}$ and $s_{2}\geq 0$. Secondly,  we prove that
the problem is globally well-posed in $H^{s_{1},\>0}(\R^{2})$ with $s_{1}>-\frac{(\alpha-1)(3\alpha-4)}{4(5\alpha+3)}$ if $4\leq \alpha \leq5$.
Finally, we prove that
the problem is globally well-posed in $H^{s_{1},\>0}(\R^{2})$ with $s_{1}>-\frac{\alpha(3\alpha-4)}{4(5\alpha+4)}$ if $\alpha>5.$
Our result improves the result of
Saut and  Tzvetkov (J. Math. Pures Appl. 79(2000), 307-338.) and Li and Xiao (J. Math. Pures Appl. 90(2008), 338-352.).

\bigskip
\bigskip

{\large\bf 1. Introduction}
\bigskip

\setcounter{Theorem}{0} \setcounter{Lemma}{0}

\setcounter{section}{1}

In this paper, we consider the Cauchy problem for the fifth-order KP-I equation
\begin{eqnarray}
&&  u_{t}+|D_{x}|^{\alpha}\partial_{x}u+\partial_{x}^{-1}\partial_{y}^{2}u+\frac{1}{2}\partial_{x}(u^{2})=0,\label{1.01}\\
&&u(x,y,0)=u_{0}(x,y)\label{1.02}
\end{eqnarray}
in anisotropic Sobolev space $H^{s_{1},s_{2}}(\R^{2})$ defined in page 6.
(\ref{1.01}) occurs in the modeling of certain long dispersive waves
 \cite{AS,KB1991,KB-1991}.
(\ref{1.01}) is the higher-order version of the KP equation
\begin{eqnarray}
 u_{t}+\beta\partial_{x}^{3}u+\partial_{x}^{-1}\partial_{y}^{2}u
 +\frac{1}{2}\partial_{x}(u^{2})=0,\alpha\neq 0.\label{1.03}
\end{eqnarray}
When $\beta<0$, (\ref{1.03}) is the KP-I equation. When $\beta>0$,
(\ref{1.03}) is the KP-II equation.
The KP-I and KP-II equations arise in physical contexts as models for the propagation
of dispersive long waves with weak transverse
effects \cite{KP}, which are two-dimensional extensions
 of the Korteweg-de-Vries
equation \cite{KP}.

  Many people have investigated the Cauchy problem for KP equation,
 for instance, see \cite{Bourgain-GAFA-KP,BS,CKS-GAFA,CIKS,GPW,Hadac2008,Hadac2009,
 HN,IMCPDE,IMEJDE,IMT,IKT,ILP,ILM-CPDE,IM2011,Kenig,KL,
 LJDE,MST-Duke,MST2002,MST2004,MST2011,ST-CMP,TDCDS,TADE,TT,TCPDE,TDIE,TIMRN,Z}
 and the references therein. Bourgain \cite{Bourgain-GAFA-KP} established
   the global well-posedness of  the Cauchy problem for the
 KP-II equation  in $L^{2}(\R^{2})$ and $L^{2}(\mathbf{T}^{2}).$ Takaokao and  Tzvetkov \cite{TT} and
 Isaza and  Mej\'{\i}a \cite{IMCPDE} established the local well-posedness
  of KP-II equation in $H^{s_{1},s_{2}}(\R^{2})$
 with $s_{1}>-\frac13$ and $s_{2}\geq0.$ Takaoka \cite{TDCDS} established the local well-posedness
 of KP-II equation in $H^{s_{1},0}(\R^{2})$
 with $s_{1}>-\frac12$ under the assumption that $D_{x}^{-\frac{1}{2}+\epsilon}u_{0}\in L^{2}$
 with the suitable chosen $\epsilon$, where $D_{x}^{-\frac{1}{2}+\epsilon}$ is Fourier multiplier
 operator with multiplier $|\xi|^{-\frac{1}{2}+\epsilon}.$
 Hadac  et al. \cite{Hadac2009}  established the small data global well-posedness
 and scattering result of  KP-II equation in the homogeneous anisotropic Sobolev space
$\dot{H}^{-\frac{1}{2},\>0}(\R^{2})$ which can be seen in \cite{Hadac2009} and arbitrary large initial data local well-posedness
 in both homogeneous Sobolev space $\dot{H}^{-\frac{1}{2},\>0}(\R^{2})$ and inhomogeneous anisotropic  Sobolev space
$H^{-\frac{1}{2},\>0}(\R^{2})$. It is proved that the Cauchy problem for KP-I equation is globally well-posed
in the  second energy spaces on both $\R^{2}$ and $\mathbf{T^{2}}$ \cite{Kenig,MST2002,MST2004}.
For KP-I equation, Molinet et al. \cite{MST-Duke} proved that the Picard
  iterative methods fails in standard Sobolev space and in anisotropic
Sobolev space  since the flow map fails to be real-analytic at the origin in these spaces.
By  introducing  some  resolution  spaces and bootstrap  inequality as well as the energy estimate,
Ionescu et al. \cite{IKT} established the  global  well-posedness of  KP-I in the natural energy space
$$E^{1}=\left\{u_{0} \in L^{2}(\R^{2}),\partial_{x}u_{0} \in L^{2}(\R^{2}),
\partial_{x}^{-1}\partial_{y}u_{0}\in L^{2}(\R^{2})\right\}.$$
Molinet et al. \cite{MST2007} proved that the Cauchy problem for the KP-I equation
 is locally well-posed
in $H^{s,\>0}(\R^{2})$ with $s>\frac{3}{2}.$
Guo et al. \cite{GPW} proved that the Cauchy problem for the KP-I equation is
 locally well-posed
in $H^{1,\>0}(\R^{2}).$ Zhang \cite{Z} proved that periodic  KP-I initial value problem
 is locally well-posed
in the Besov type space
$
B_{2,1}^{\frac{1}{2}}(\mathbf{T}^{2}).
$
It is worth noticing that the resonant function of KP-I equation does not possess the
 good property as its of KP-II equation.

When $\alpha=4,$ (\ref{1.01}) reduces to the fifth order KP-I equation
\begin{eqnarray}
u_{t}+\partial_{x}^{5}u+\partial_{x}^{-1}\partial_{y}^{2}u+\frac{1}{2}\partial_{x}(u^{2})=0.\label{1.03}
\end{eqnarray}
Some people have studied the Cauchy problem for (\ref{1.03}), for instance, see \cite{ST2000, CLM,LX,GHF}.
Saut and Tzvetkov \cite{ST2000} proved that the Cauchy problem for (\ref{1.01})
 is globally well-posed
for initial data $u_{0} \in L^{2}(\R^{2})$ with finite energy.
Chen et al. \cite{CLM} proved that the problem for (\ref{1.01}) is locally well-posed
 in $E^{s}$ with $0<s\leq 1,$ where
\begin{eqnarray*}
E^{s}=\left\{u_{0}\in E^{s}:\|u_{0}\|_{E^{s}}
=\left\|\left(1+|\xi|^{2}+\left|\frac{\mu}{\xi}\right|\right)^{s}
\mathscr{F}_{xy}u_{0}(\xi,\mu)\right\|_{L^{2}}<\infty\right\}.
\end{eqnarray*}
By using the Fourier restriction norm method and sufficiently exploiting
 the geometric structure of the resonant set of (\ref{1.01}) to dispose the high-high frequency
   interaction,  Li and Xiao \cite{LX} proved that
 the Cauchy problem for (\ref{1.01}) is  globally well-posed in
  $L^{2}(\R^{2}).$
Recently, Guo et al. \cite{GHF} proved that the Cauchy problem
 for (\ref{1.03}) is locally well-posed
in $H^{s_{1},s_{2}}(\R^{2})$ with $s_{1}\geq-\frac34,s_{2}\geq0$.
Saut and Tzvetkov \cite{ST1999} proved that the fifth-order KP-II equation
\begin{eqnarray}
u_{t}-\partial_{x}^{5}u+\partial_{x}^{-1}\partial_{y}^{2}u
+\frac{1}{2}\partial_{x}(u^{2})=0\label{1.04}
\end{eqnarray}
is locally well-posed
in $H^{s_{1},s_{2}}(\R^{2})$ with $s_{1}>-\frac14,s_{2}\geq0$. Isaza et al.
\cite{ILM-CPAA} proved that the Cauchy problem for
(\ref{1.04}) is locally well-posed
in $H^{s_{1},s_{2}}(\R^{2})$ with $s_{1}>-\frac54,s_{2}\geq0$ and globally well-posed
in $H^{s_{1},0}(\R^{2})$ with $s_{1}>-\frac47.$ Recently, Li and Shi \cite{LS} proved that
the Cauchy problem for
(\ref{1.04}) is locally well-posed
in $H^{s_{1},s_{2}}(\R^{2})$ with $s_{1}\geq-\frac54,s_{2}\geq0$.

Recently, Linares et al. \cite{LPS} proved various ill-posedness
 and wellposedness
results on the Cauchy problem
\begin{eqnarray}
u_{t}+uu_{x}-D_{x}^{\alpha}u_{x}+\gamma\partial_{x}^{-1}u_{yy}=0,
\gamma\in \R,0<\alpha \leq 1.\label{1.05}
\end{eqnarray}
To the best of our knowledge, the Cauchy problem for (\ref{1.01})
 in low regularity space
 is yet  to  be  answered with $\alpha\geq4$.
 The main reason is that
the resonant function of KP-I type equation does not enjoy
 the same good property as its of  KP-II type equation.

In this paper, inspired by \cite{CKS-GAFA,ST2000,LX,ILM-CPAA},
by using the Fourier restriction norm method introduced in
\cite{Beals,Bourgain93,KM,RR}
and developed in \cite{KPV1993,KPV1996}, the Cauchy-Schwartz
inequality
 and Strichartz estimates
 as well as  suitable splitting of domains,  we prove
 that the Cauchy problem for (\ref{1.01})
is locally well-posed  in the anisotropic Sobolev spaces $H^{s_{1},\>s_{2}}(\R^{2})$
  with $s_{1}>-\frac{\alpha-1}{4}$ and $s_{2}\geq 0$; combining the
   local well-posness result
  of this paper with the I-method introduced in
  \cite{CKSTT2001,CKSTT2003},
   we also prove that
the problem is globally well-posed in $H^{s_{1},\>0}(\R^{2})$ with
 $s_{1}>-\frac{(\alpha-1)(3\alpha-4)}{4(5\alpha+3)}$
 if $4\leq \alpha \leq5$
and prove that
the problem is globally well-posed in $H^{s_{1},\>0}(\R^{2})$ with $s_{1}>-\frac{\alpha(3\alpha-4)}{4(5\alpha+4)}$ if $\alpha>5.$
 Thus, our result improves the result of \cite{ST2000,LX}.

We introduce some notations before giving the main results.
 Throughout this paper, we assume that
$C$ is a positive constant which may depend upon $\alpha$
 and  vary from line to line. $a\sim b$ means that there exist constants $C_{j}>0(j=1,2)$ such that $C_{1}|b|\leq |a|\leq C_{2}|b|$.
 $a\gg b$ means that there exist a positive constant $C^{\prime}$ such that  $|a|> C^{\prime}|b|.$ $0<\epsilon\ll1$ means that $0<\epsilon<\frac{1}{100\alpha}$.
 We define
  \begin{eqnarray*}
  &&\langle\cdot\rangle:=1+|\cdot|,\\
  &&\phi(\xi,\mu):=\xi|\xi|^{\alpha}+\frac{\mu^{2}}{\xi},\\
  &&\sigma:=\tau+\phi(\xi,\mu),\sigma_{j}=\tau_{j}+\phi(\xi_{j},\mu_{j})(j=1,2),\\
  &&\mathscr{F}u(\xi,\mu,\tau):=\frac{1}{(2\pi)^{\frac{3}{2}}}\int_{\SR^{3}}e^{-ix\xi-iy\mu-it\tau}u(x,y,t)dxdydt,\\
  &&\mathscr{F}_{xy}f(\xi,\mu):=\frac{1}{2\pi}\int_{\SR^{2}}e^{-ix\xi-iy\mu}f(x,y)dxdy,\\
  &&\mathscr{F}^{-1}u(\xi,\mu,\tau):=\frac{1}{(2\pi)^{\frac{3}{2}}}\int_{\SR^{3}}e^{ix\xi+iy\mu+it\tau}u(x,y,t)dxdydt,\\
  &&D_{x}^{a}u(x,y,t):=\frac{1}{(2\pi)^{\frac{3}{2}}}\int_{\SR^{2}}|\xi|^{a}\mathscr{F}u(\xi,\mu,\tau)e^{ix\xi+iy\mu+it\tau}d\xi d\mu d\tau,\\
 && W(t)f:=\frac{1}{2\pi}\int_{\SR^{2}}e^{ix\xi+iy\mu+it\phi(\xi,\mu)}\mathscr{F}_{xy}f(\xi,\mu)d\xi d\mu.
  \end{eqnarray*}
   Let $\eta$ be a bump function with compact support in $[-2,2]\subset \R$
  and $\eta=1$ on $(-1,1)\subset \R$.
 For each integer $j\geq1$, we define $\eta_{j}(\xi)=\eta(2^{-j}\xi)-\eta(2^{1-j}\xi),$
  $\eta_{0}(\xi)=\eta(\xi),$
 $\eta_{j}(\xi,\mu,\tau)=\eta_{j}(\sigma),$ thus, $\sum\limits_{j\geq0}\eta_{j}(\sigma)=1.$
 $\psi(t)$ is a smooth function
 supported in $[0,2]$ and equals
 $1$ in $[0,1]$.
 Let $I\subset \R$, $\chi_{I}(x)=1$ if $x\in I$; $\chi_{I}(x)=0$
  if $x$ doesnot belong to $I$.  We  denote by ${\rm mes}(E)$ the Lebesgue measure of a set $E.$
We define $|\xi_{\rm min}|:={\rm min}\left\{|\xi|,|\xi_{1}|,|\xi_{2}|\right\}$ and
$|\xi_{\rm max}|:={\rm max}\left\{|\xi|,|\xi_{1}|,|\xi_{2}|\right\}.$
 We define
 \begin{eqnarray*}
 \|f\|_{L_{t}^{r}L_{xy}^{p}}:=\left(\int_{\SR}\left(\int_{\SR^{2}}|f|^{p}dxdy\right)^{\frac{r}{p}}dt\right)^{\frac{1}{r}}.
 \end{eqnarray*}
 The anisotropic Sobolev space $H^{s_{1},s_{2}}$ is defined as follows:
 \begin{eqnarray*}
 H^{s_{1},s_{2}}(\R^{2}):=\left\{u_{0}\in \mathscr{S}^{'}(\R^{2}):\quad \|u_{0}\|_{H^{s_{1},s_{2}}(\SR^{2})}=\left\|\langle\xi\rangle^{s_{1}}
 \langle\mu\rangle^{s_{2}}\mathscr{F}_{xy}u_{0}(\xi,\mu)\right\|_{L_{\xi\mu}^{2}}\right\}.
 \end{eqnarray*}
 Space
$
  X_{b}^{s_{1},s_{2}}
$ is defined by
$$
X_{b}^{s_{1},s_{2}}:= \left\{u\in  \mathscr{S}^{'}(\R^{3})
 :\, \|u\|_{X_{b}^{s_{1},s_{2}}}
 =  \left\|\langle\xi\rangle^{s_{1}} \langle\mu\rangle^{s_{2}}
 \left\langle\sigma\right\rangle^{b}\mathscr{F}u(\xi,\mu,\tau)
 \right\|_{L_{\tau\xi\mu}^{2}(\SR^{3})}<\infty\right\}.
$$
The space $ X_{b}^{s_{1},s_{2}}([0,T])$ denotes the restriction
 of $X_{b}^{s_{1},s_{2}}$ onto the finite time interval $[0,T]$ and
is equipped with the norm
 \begin{equation*}
    \|u\|_{X_{b}^{s_{1},s_{2}}([0,T])} =\inf \left\{\|g\|_{X_{b}^{s_{1},s_{2}}}
    :g\in X_{b}^{s_{1},s_{2}}, u(t)=g(t)
 \>\> {\rm for} \>  t\in [0,T]\right\}.
 \end{equation*}
For $s<0$ and $N\in N^{+},N\geq100$, inspired by \cite{CKSTT2001, CKSTT2003}, we define an operator $I_{N}$  by
 $\mathscr{F}I_{N}u(\xi,\mu,\tau)=M(\xi)\mathscr{F}u(\xi,\mu,\tau)$,
where $M(\xi)=1$ if $|\xi|<N$; $M(\xi)=\left(\frac{|\xi|}{N}\right)^{s}$ if $|\xi|\geq N.$

The main results of this paper are as follows.

\begin{Theorem}\label{Thm1}(Local well-posedness)

\noindent Let $|\xi|^{-1}\mathscr{F}_{xy}u_{0}(\xi,\mu)\in \mathscr{S}^{'}(\R^{2})$. Then,  the Cauchy problem  for (\ref{1.01}) are locally well-posed in  the anisotropic Sobolev spaces
 $H^{s_{1},\>s_{2}}(\R^{2})$ with $s_{1}>-\frac{\alpha-1}{4},$$\>s_{2}\geq0.$
\end{Theorem}
\noindent{\bf Remark 1.}Note that the resonant function of generalized KP-II equation
is
\begin{eqnarray}
&&R_{II}(\xi_{1},\xi_{2},\mu_{1},\mu_{2}):
=\left[|\xi|^{\alpha}\xi-\xi_{1}|\xi_{1}|^{\alpha}-\xi_{2}|\xi_{2}|^{\alpha}
+\frac{\xi_{1}\xi_{2}}{\xi}\left(\frac{\mu_{1}}{\xi_{1}}-\frac{\mu_{2}}{\xi_{2}}\right)^{2}\right].\label{1.07}
\end{eqnarray}
Combining Lemma 2.8 of \cite{YLZ} with (\ref{1.07}),  we have that
\begin{eqnarray}
\left|R_{II}(\xi_{1},\xi_{2},\mu_{1},\mu_{2})\right|\geq \left||\xi|^{\alpha}\xi-\xi_{1}|\xi_{1}|^{\alpha}-\xi_{2}|\xi_{2}|^{\alpha}\right|.\label{1.08}
\end{eqnarray}
However, the resonant function of the generalized KP-I equation is
\begin{eqnarray}
&&\sigma-\sigma_{1}-\sigma_{2}=R_{I}(\xi_{1},\xi_{2},\mu_{1},\mu_{2}):=\phi(\xi,\mu)-\phi(\xi_{1},\mu_{1})-\phi(\xi_{2},\mu_{2})\nonumber\\&&
=\left[|\xi|^{\alpha}\xi-\xi_{1}|\xi_{1}|^{\alpha}-\xi_{2}|\xi_{2}|^{\alpha}
-\frac{\xi_{1}\xi_{2}}{\xi}\left(\frac{\mu_{1}}{\xi_{1}}-\frac{\mu_{2}}{\xi_{2}}\right)^{2}\right],\label{1.09}
\end{eqnarray}
thus, we consider
\begin{eqnarray}
&&\left|R_{I}(\xi_{1},\xi_{2},\mu_{1},\mu_{2})\right|\geq \frac{\left||\xi|^{\alpha}\xi-\xi_{1}|\xi_{1}|^{\alpha}-\xi_{2}|\xi_{2}|^{\alpha}\right|}{\alpha},\label{1.010}\\
&&\left|R_{I}(\xi_{1},\xi_{2},\mu_{1},\mu_{2})\right|< \frac{\left||\xi|^{\alpha}\xi-\xi_{1}|\xi_{1}|^{\alpha}-\xi_{2}|\xi_{2}|^{\alpha}\right|}{\alpha},\label{1.011}
\end{eqnarray}
respectively. When (\ref{1.010}) is valid, we follow the procedure of Lemma 3.1 in \cite {YLZ}
to obtain Lemma 3.1 of this paper. When (\ref{1.011}) is valid, inspired by \cite{CKS-GAFA,LX}, we sufficiently exploit the geometric structure
of (\ref{1.011}) to Lemma 3.1 of this paper. When $\alpha=4$, we improve the local well-posednes result of \cite{LX}.

\begin{Theorem}\label{Thm2}(Global well-posedness when $4\leq \alpha\leq 5$)

\noindent
Let $4\leq\alpha\leq5$ and  $|\xi|^{-1}\mathscr{F}_{xy}u_{0}(\xi,\mu)\in \mathscr{S}^{'}(\R^{2})$. Then, the Cauchy problem  for (\ref{1.01}) are globally well-posed in $H^{s_{1},\>0}(\R^{2})$
with $s_{1}>-\frac{(\alpha-1)(3\alpha-4)}{4(5\alpha+3)}.$
\end{Theorem}

\noindent{\bf Remark 2.} In proving Theorem 1.2, we only present the proof
of case $-\frac{(\alpha-1)(3\alpha-4)}{4(5\alpha+3)}<s_{1}<0$ since
the global well-posedness of  the Cauchy problem  for (\ref{1.01})
 in $H^{s_{1},s_{2}}(\R^{2})$
with $s_{1}\geq0,s_{2}\geq0$ can be easily proved with the aid of $L^{2}$
conservation law of (\ref{1.01}).
When $\alpha=4$, we improve the global well-posedness result of \cite{LX}.
The establishment of Lemma 3.2 plays the crucial role in proving Theorem 1.2.
\noindent When (\ref{1.010}) is valid,  we follow the method of Lemma 3.2
 of \cite{YLZ} to deal with  case (\ref{1.010}).
\noindent When (\ref{1.011}) is valid,  we use the technique used in
 Lemma 3.1 of this paper to deal with case (\ref{1.011}).

\begin{Theorem}\label{Thm3}(Global well-posedness when $\alpha >5$)

\noindent
Let $\alpha>5$ and  $|\xi|^{-1}\mathscr{F}_{xy}u_{0}(\xi,\mu)\in \mathscr{S}^{'}(\R^{2})$.
 Then, the Cauchy problem  for (\ref{1.01}) are globally well-posed in $H^{s_{1},\>0}(\R^{2})$
with $s_{1}>-\frac{\alpha(3\alpha-4)}{4(5\alpha+4)}.$
\end{Theorem}
In proving Theorem 1.3, we only present the proof
of case $-\frac{\alpha(3\alpha-4)}{4(5\alpha+4)}<s_{1}<0$ since
the global well-posedness of  the Cauchy problem  for (\ref{1.01})
 in $H^{s_{1},s_{2}}(\R^{2})$
with $s_{1}\geq0,s_{2}\geq0$ can be easily proved with the aid of $L^{2}$
conservation law of (\ref{1.01}).

The rest of the paper is arranged as follows. In Section 2,  we give some
preliminaries. In Section 3, we establish two crucial bilinear estimates.
 In Section 4, we prove
the Theorem 1.1. In Section 5, we prove
the Theorem 1.2. In Section 6, we prove
the Theorem 1.3.

\bigskip
\bigskip

 \noindent{\large\bf 2. Preliminaries }

\setcounter{equation}{0}

\setcounter{Theorem}{0}

\setcounter{Lemma}{0}

\setcounter{section}{2}

In this section, motivated by \cite{Bourgain-GAFA-KP, MST2011}, we give
   Lemmas 2.1-2.6 which play a significant role in establishing
  Lemmas 3.1, 3.2. Lemma 2.2 in combination with Lemma 3.1  yields Theorem 1.1. Lemma 2.7 in combination with Lemmas 3.1, 3.2 yields
  Lemma 5.1.

\begin{Lemma}\label{Lemma2.1}
 Let $b>|a|\geq 0$. Then, we have that
 \begin{eqnarray}
 &&\int_{-b}^{b}\frac{dx}{\langle x+a\rangle^{\frac{1}{2}}}\leq
  Cb^{\frac{1}{2}},\label{2.01}\\
 &&\int_{\SR}\frac{dt}{\langle t\rangle^{\gamma}\langle t-a\rangle^{\gamma}}
 \leq C\langle a\rangle^{-\gamma},\gamma>1,\label{2.02}\\
 &&\int_{\SR}\frac{dt}{\langle t\rangle^{\gamma}|t-a|^{\frac{1}{2}}}
 \leq C\langle a\rangle^{-\frac{1}{2}},\gamma\geq1,\label{2.03}\\
 &&\int_{-K}^{K}\frac{dx}{|x|^{\frac{1}{2}}|a-x|^{\frac{1}{2}}}
 \leq C\frac{K^{\frac{1}{2}}}{|a|^{\frac{1}{2}}}.\label{2.04}
 \end{eqnarray}
 \end{Lemma}
 {\bf Proof.} The  conclusion of (\ref{2.01}) is given in (2.4) of Lemma 2.1 in \cite{ILM-CPAA}.
(\ref{2.02})-(\ref{2.03})  can be seen in Proposition 2.2 of \cite{ST2000}. (\ref{2.04})
can be seen in line 24 of page 6562 in \cite{Hadac2008}.

This completes the proof of Lemma 2.1.
\begin{Lemma}\label{Lemma2.2}
Let $T\in (0,1)$ and $s_{1},s_{2} \in \R$ and $-\frac{1}{2}<b^{\prime}\leq0\leq b\leq b^{\prime}+1$.
Then, for $h\in X_{b^{\prime}}^{s_{1},s_{2}},$  we have that
\begin{eqnarray}
&&\left\|\psi(t)S(t)\phi\right\|_{X_{b}^{s_{1},s_{2}}}\leq C\|\phi\|_{H^{s_{1},\>s_{2}}},\label{2.05}\\
&&\left\|\psi\left(\frac{t}{T}\right)\int_{0}^{t}S(t-\tau)h(\tau)d\tau\right\|_{X_{b}^{s_{1},\>s_{2}}}\leq C
T^{1+b^{\prime}-b}\|h\|_{X_{b^{\prime}}^{s_{1},\>s_{2}}}.\label{2.06}
\end{eqnarray}
\end{Lemma}

For the proof of Lemma 2.2, we refer the readers to \cite{G,Bourgain93,KPV1993}
 and  Lemmas 1.7, 1.9 of  \cite{Grunrock}.

\begin{Lemma}\label{Lemma2.3}
Let $b>\frac{1}{2}$. Then, we have
\begin{eqnarray}
\left\|D_{x}^{\frac{\alpha-2}{8}}u\right\|_{L_{t}^{4}L_{xy}^{4}(\SR^{3})}
\leq C\|u\|_{X_{b}^{0,0}}.\label{2.07}
\end{eqnarray}
\end{Lemma}
\noindent{\bf Proof.} Combining Lemma 3.1 of \cite{Hadac2008} with Lemma 3.3 of \cite{G}, we have that Lemma 2.3 is valid.

This completes the proof of Lemma 2.3.

Motivated by \cite{TCPDE,TDCDS,TDIE,TADE,TT,IMCPDE}
 and Theorem 3.3 of \cite{Hadac2008}, we present the proof of  Lemma 2.4.
\begin{Lemma}\label{Lemma2.4}
Let
\begin{eqnarray*}
&&|\sigma-\sigma_{1}-\sigma_{2}|=\left|\xi|\xi|^{\alpha}-\xi_{1}|\xi_{1}|^{\alpha}-\xi_{2}|\xi_{2}|^{\alpha}-\frac{\xi_{1}\xi_{2}}{\xi}
\left|\frac{\mu_{1}}{\xi_{1}}-\frac{\mu_{2}}{\xi_{2}}\right|^{2}\right|\nonumber\\&&
\geq \frac{\left||\xi|^{\alpha}\xi-\xi_{1}|\xi_{1}|^{\alpha}-\xi_{2}|\xi_{2}|^{\alpha}\right|}{\alpha}
\end{eqnarray*}
and
\begin{eqnarray*}
&&\hspace{-1cm}\mathscr{F}P_{\frac{1}{4}}(u_{1},u_{2})(\xi,\mu,\tau)=\int_{\SR^{3}}\chi_{|\xi_{1}|\leq \frac{|\xi_{2}|}{4}}
(\xi_{1},\mu_{1},\tau_{1},\xi,\mu,\tau)\prod\limits_{j=1}^{2}\mathscr{F}u_{j}(\xi_{j},\mu_{j},\tau_{j})
d\xi_{1}d\mu_{1}d\tau_{1}.
\end{eqnarray*}
For $b>\frac{1}{2}$, we have
\begin{eqnarray}
\left\|P_{\frac{1}{4}}(u_{1},u_{2})\right\|_{L_{T}^{2}L_{xy}^{2}}
\leq C\left\||D_{x}|^{\frac{1}{2}}u_{1}\right\|_{X_{b}^{0,0}}\left\||D_{x}|^{-\frac{\alpha}{4}}u_{2}\right\|_{X_{b}^{0,0}}.\label{2.08}
\end{eqnarray}
\end{Lemma}
\noindent {\bf Proof.} Let
\begin{eqnarray*}
f_{1}(\xi_{1},\mu_{1},\tau_{1})=|\xi_{1}|^{\frac{1}{2}}\langle \sigma_{1}\rangle^{b}\mathscr{F}u_{1}(\xi_{1},\mu_{1},\tau_{1}),
f_{2}(\xi_{2},\mu_{2},\tau_{2})=|\xi_{2}|^{-\frac{\alpha}{4}}\langle \sigma_{2}\rangle^{b}\mathscr{F}u_{2}(\xi_{2},\mu_{2},\tau_{2}).
\end{eqnarray*}
 To obtain
(\ref{2.08}), it suffices to prove that
\begin{eqnarray}
\left\|\int_{\SR^{3}}\chi_{|\xi_{1}|\leq \frac{|\xi_{2}|}{4}}
\frac{|\xi_{1}|^{-\frac{1}{2}}|\xi_{2}|^{\frac{\alpha}{4}}f_{1}(\xi_{1},\mu_{1},\tau_{1})
f_{2}(\xi_{2},\mu_{2},\tau_{2})}{\prod\limits_{j=1}^{2}\langle\sigma_{j}\rangle^{b}}
d\xi_{1}d\mu_{1}d\tau_{1}\right\|_{L_{\tau\xi\mu}^{2}}\leq C\prod_{j=1}^{2}\|f_{j}\|_{L_{\tau\xi\mu}^{2}}.\label{2.09}
\end{eqnarray}
To obtain (\ref{2.09}), by duality, it suffices to prove that
\begin{eqnarray}
&&\left|\int_{\SR^{6}}\chi_{|\xi_{1}|\leq \frac{|\xi_{2}|}{4}}
\frac{|\xi_{1}|^{-\frac{1}{2}}|\xi_{2}|^{\frac{\alpha}{4}}f(\xi,\mu,\tau)\prod\limits_{j=1}^{2}f_{j}(\xi_{j},\mu_{j},\tau_{j})
}{\prod\limits_{j=1}^{2}\langle\sigma_{j}\rangle^{b}}
d\xi_{1}d\mu_{1}d\tau_{1}d\xi d\mu d\tau\right|\nonumber\\&&\leq C\|f\|_{L_{\tau\xi\mu}^{2}}\prod_{j=1}^{2}\|f_{j}\|_{L_{\tau\xi\mu}^{2}}.\label{2.010}
\end{eqnarray}
We define
\begin{eqnarray}
I(\xi,\mu,\tau):=\int_{\SR^{3}}\chi_{|\xi_{1}|\leq \frac{|\xi_{2}|}{4}}
\frac{|\xi_{1}|^{-1}|\xi_{2}|^{\frac{\alpha}{2}}}{\prod\limits_{j=1}^{2}\langle\sigma_{j}\rangle^{b}}
d\xi_{1}d\mu_{1}d\tau_{1}.\label{2.011}
\end{eqnarray}
For fixed $(\xi,\mu,\tau),$ we make the change of variables ${\rm L}:(\xi_{1},\mu_{1},\tau_{1})
\longrightarrow (\Delta,\sigma_{1},\sigma_{2})$,
where
\begin{eqnarray*}
\Delta:=\xi|\xi|^{\alpha}-\xi_{1}|\xi_{1}|^{\alpha}-\xi_{2}|\xi_{2}|^{\alpha},
\sigma_{1}:=\tau_{1}+\phi(\xi_{1},\mu_{1}),\sigma_{2}:=\tau_{2}+\phi(\xi_{2},\mu_{2}).
\end{eqnarray*}
By using a direct computation, since $\sigma=\tau+\phi(\xi,\mu),$ we have that
\begin{eqnarray}
\sigma_{1}+\sigma_{2}-\sigma=-\Delta+\frac{(\xi_{1}\mu_{2}-\mu_{1}\xi_{2})^{2}}{\xi\xi_{1}\xi_{2}}.\label{2.012}
\end{eqnarray}
Thus, we have that the Jacobian determinant
equals
\begin{eqnarray}
&&\frac{\partial(\Delta,\sigma_{1},\sigma_{2})}{\partial(\xi_{1},\mu_{1},\tau_{1})}
=-2(\alpha+1)\left(\xi_{1}^{\alpha}-\xi_{2}^{\alpha}\right)
\left(\frac{\mu_{1}}{\xi_{1}}-\frac{\mu_{2}}{\xi_{2}}\right)\nonumber\\
&&=-2(\alpha+1)\left(\xi_{1}^{\alpha}-\xi_{2}^{\alpha}\right)
(\sigma_{1}+\sigma_{2}-\sigma+\Delta)^{\frac{1}{2}}
\left(\frac{\xi}{\xi_{1}\xi_{2}}\right)^{\frac{1}{2}}.\label{2.013}
\end{eqnarray}
Notice that it is possible to divide the integration into a finite number of open subsets $W_{i}$ such that
${\rm L}$ is an injective $C^{1}$-function in $W_{i}$ with non-zero Jacobian determinant.
From (\ref{2.013}), since $\frac{|\xi_{2}|}{4}\geq |\xi_{1}|$ and $|\Delta|\sim |\xi_{1}||\xi_{2}|^{\alpha},$ we have that
\begin{eqnarray}
&&\left|\frac{\partial(\Delta,\sigma_{1},\sigma_{2})}{\partial(\xi_{1},\mu_{1},\tau_{1})}\right|=2(\alpha+1)\left|\left(\xi_{1}^{2}-\xi_{2}^{2}\right)
\left(\xi_{1}^{2}+\xi_{2}^{2}\right)\left(\frac{\mu_{1}}{\xi_{1}}-\frac{\mu_{2}}{\xi_{2}}\right)\right|\nonumber\\
&&=\left|\left(\xi_{1}^{\alpha}-\xi_{2}^{\alpha}\right)
(\sigma_{1}+\sigma_{2}-\sigma+\Delta)^{\frac{1}{2}}\left(\frac{\xi}{\xi_{1}\xi_{2}}\right)^{\frac{1}{2}}\right|\nonumber\\
&&\sim |\xi_{1}|^{-1}|\xi_{2}|^{\frac{\alpha}{2}}|\Delta|^{\frac{1}{2}}\left|\sigma_{1}+\sigma_{2}-\sigma+\Delta\right|^{\frac{1}{2}}.\label{2.014}
\end{eqnarray}
Since $\left|\sigma_{1}+\sigma_{2}-\sigma\right|\geq \frac{|\Delta|}{\alpha}$,
 by using the change of variables $(\xi_{1},\mu_{1},\tau_{1})\longrightarrow (\Delta,\sigma_{1},\sigma_{2})$ and (\ref{2.04}), we have that
\begin{eqnarray}
&&I(\xi,\mu,\tau):=\int_{\SR^{3}}\chi_{|\xi_{1}|\leq \frac{|\xi_{2}|}{4}}
\frac{|\xi_{2}|^{\frac{\alpha}{2}}|\xi_{1}|^{-1}}{\prod\limits_{j=1}^{2}\langle\sigma_{j}\rangle^{2b}}
d\xi_{1}d\mu_{1}d\tau_{1}\nonumber\\
&&\leq C\int_{\SR^{3}}\frac{\chi_{|\Delta|\leq 4|\sigma_{1}+\sigma_{2}-\sigma|d\Delta d\sigma_{1}d\sigma_{2}}}{|\Delta|^{\frac{1}{2}}
\left|\sigma_{1}+\sigma_{2}-\sigma+\Delta\right|^{\frac{1}{2}}\prod\limits_{j=1}^{2}\langle\sigma_{j}\rangle^{2b}}\nonumber\\
&&=\int_{\SR^{2}}\left(\int_{\SR}\frac{\chi_{|\Delta|\leq 4|\sigma_{1}+\sigma_{2}-\sigma|d\Delta }}{|\Delta|^{\frac{1}{2}}
\left|\sigma_{1}+\sigma_{2}-\sigma+\Delta\right|^{\frac{1}{2}}}\right)\frac{d\sigma_{1}d\sigma_{2}}
{\prod\limits_{j=1}^{2}\langle\sigma_{j}\rangle^{2b}}\leq C\int_{\SR^{2}}\frac{d\sigma_{1}d\sigma_{2}}
{\prod\limits_{j=1}^{2}\langle\sigma_{j}\rangle^{2b}}\leq C.\label{2.015}
\end{eqnarray}
Combining (\ref{2.010}) with (\ref{2.015}), by using the Cauchy-Schwartz inequality twice, we have that
\begin{eqnarray}
&&\hspace{-1cm}\int_{\SR^{6}}\chi_{|\xi_{1}|\leq \frac{|\xi_{2}|}{4}}
\frac{|\xi_{2}|^{\frac{\alpha}{4}}|\xi_{1}|^{-\frac{1}{2}}f_{1}(\xi_{1},\mu_{1},\tau_{1})
f_{2}(\xi_{2},\mu_{2},\tau_{2})f(\xi,\mu,\tau)}{\prod\limits_{j=1}^{2}\langle\sigma_{j}\rangle^{b}}
d\xi_{1}d\mu_{1}d\tau_{1}d\xi d\mu d\tau\nonumber\\
&&\hspace{-1cm}\leq C\left[\sup \limits_{\xi,\mu,\tau}I(\xi,\mu,\tau)\right]^{\frac{1}{2}}
\|f\|_{L_{\tau\xi\mu}^{2}}\left(\prod_{j=1}^{2}\|f_{j}\|_{L_{\tau\xi\mu}^{2}}\right)\leq C\|f\|_{L_{\tau\xi\mu}^{2}}
\left(\prod_{j=1}^{2}\|f_{j}\|_{L_{\tau\xi\mu}^{2}}\right).\label{2.016}
\end{eqnarray}

This completes the proof of Lemma 2.4.

Inspired by Proposition 3.5 of \cite{Hadac2008}, we present the proof of Lemma 2.5.

\begin{Lemma}\label{Lemma2.5}
Assume
\begin{eqnarray*}
&&|\sigma-\sigma_{1}-\sigma_{2}|=\left|\xi|\xi|^{\alpha}-\xi_{1}|\xi_{1}|^{\alpha}-\xi_{2}|\xi_{2}|^{\alpha}-\frac{\xi_{1}\xi_{2}}{\xi}
\left|\frac{\mu_{1}}{\xi_{1}}-\frac{\mu_{2}}{\xi_{2}}\right|^{2}\right|\nonumber\\&&
\geq \frac{\left||\xi|^{\alpha}\xi-\xi_{1}|\xi_{1}|^{\alpha}-\xi_{2}|\xi_{2}|^{\alpha}\right|}{\alpha}
\end{eqnarray*}
and
 $b>\frac{1}{2}$, then,
 we have that
\begin{eqnarray}
\left|\int_{\SR^{6}}\frac{|\xi_{1}|^{-\frac{1}{2}}
|\xi_{2}|^{\frac{\alpha}{4}}F(\xi,\mu,\tau)\prod\limits_{j=1}^{2}F_{j}(\xi_{j},\mu_{j},\tau_{j})}
{\prod\limits_{j=1}^{2}\langle\sigma_{j}\rangle^{b}}dV\right|
\leq C\|f\|_{L_{\tau\xi\mu}^{2}}\left(\prod\limits_{j=1}^{2}\|f_{j}\|_{L_{\tau\xi\mu}^{2}}\right)\label{2.017}
 \end{eqnarray}
and
\begin{eqnarray}
\left|\int_{\SR^{6}}\frac{|\xi_{1}|^{-\frac{1}{2}}
|\xi|^{\frac{\alpha}{4}}F(\xi,\mu,\tau)\prod\limits_{j=1}^{2}F_{j}(\xi_{j},\mu_{j},\tau_{j})}
{\langle\sigma_{1}\rangle^{b}\langle\sigma\rangle^{b}}dV\right|
\leq C\|f\|_{L_{\tau\xi\mu}^{2}}\left(\prod\limits_{j=1}^{2}\|f_{j}\|_{L_{\tau\xi\mu}^{2}}\right)\label{2.018}
 \end{eqnarray}
and
\begin{eqnarray}
\left|\int_{\SR^{6}}\frac{|\xi|^{-\frac{1}{2}}
|\xi_{2}|^{\frac{\alpha}{4}}F(\xi,\mu,\tau)\prod\limits_{j=1}^{2}F_{j}(\xi_{j},\mu_{j},\tau_{j})}
{\langle\sigma\rangle^{b}\langle\sigma_{2}\rangle^{b}}dV\right|
\leq C\|f\|_{L_{\tau\xi\mu}^{2}}\left(\prod\limits_{j=1}^{2}\|f_{j}\|_{L_{\tau\xi\mu}^{2}}\right)\label{2.019}
\end{eqnarray}
and
\begin{eqnarray}
\left|\int_{\SR^{6}}\frac{|\xi|^{-\frac{1}{2}}
|\xi_{1}|^{\frac{\alpha}{4}}F(\xi,\mu,\tau)\prod\limits_{j=1}^{2}F_{j}(\xi_{j},\mu_{j},\tau_{j})}
{\langle\sigma\rangle^{b}\langle\sigma_{1}\rangle^{b}}dV\right|
\leq C\|f\|_{L_{\tau\xi\mu}^{2}}\left(\prod\limits_{j=1}^{2}\|f_{j}\|_{L_{\tau\xi\mu}^{2}}\right)\label{2.020}
\end{eqnarray}
where $dV=d\xi_{1}d\mu_{1}d\sigma_{1}d\xi_{2}d\mu_2d\sigma_{2}$.
\end{Lemma}
\noindent {\bf Proof.} We firstly prove (\ref{2.017}). When $\frac{|\xi_{2}|}{4}\geq |\xi_{1}|$, from Lemma 2.4,
we have that (\ref{2.017}) is valid. When $\frac{|\xi_{2}|}{4}< |\xi_{1}|$, since $|\xi_{1}|^{-\frac{1}{2}}|\xi_{2}|\leq C
|\xi_{1}|^{\frac{1}{4}}|\xi_{2}|^{\frac{1}{4}}$, from Lemma 2.3, we know that (\ref{2.017}) is valid.
Let $\xi_{1}=\xi_{1}^{\prime},\mu_{1}=\mu_{1}^{\prime},\tau_{1}=\tau_{1}^{\prime}$
and $-\xi_{2}=\xi^{\prime},-\tau_{2}=\tau^{\prime},-(\mu-\mu_{1})=\mu^{\prime}$
and
 $-\xi=\xi^{\prime}-\xi_{1}^{\prime},-\mu=\mu^{\prime}-\mu_{1}^{\prime},-\tau=\tau^{\prime}-\tau_{1}^{\prime}$
and $\sigma_{2}^{\prime}=\tau_{2}^{\prime}-\phi^{\prime}(\xi_{2}^{\prime},\mu_{2}^{\prime}),\sigma_{1}=\sigma_{1}^{\prime}
=\tau_{1}^{\prime}-\phi(\xi_{1}^{\prime},\mu_{1}^{\prime})$. Thus, $-\sigma=\sigma_{2}^{\prime},\sigma_{1}=\sigma_{1}^{\prime}.$
Let
\begin{eqnarray*}
H(\xi_{1}^{\prime},\mu_{1}^{\prime},\tau_{1}^{\prime},\xi^{\prime},\mu^{\prime},
\tau^{\prime})=f_{1}(\xi_{1}^{\prime},\mu_{1}^{\prime},\tau_{1}^{\prime})
f_{2}(-\xi^{\prime},-\mu^{\prime},-\tau^{\prime})f(-\xi_{2}^{\prime},
-\mu_{2}^{\prime},-\tau_{2}^{\prime}).
\end{eqnarray*}
To obtain (\ref{2.018}), it suffices to prove that
\begin{eqnarray}
&&\left|\int_{\SR^{6}}\frac{|\xi_{1}^{\prime}|^{-\frac{1}{2}}
|\xi_{2}^{\prime}||H(\xi_{1}^{\prime},\mu_{1}^{\prime},
\tau_{1}^{\prime},\xi^{\prime},\mu^{\prime},
\tau^{\prime})}
{\langle\sigma_{1}^{\prime}\rangle^{b}\langle\sigma_{2}^{\prime}\rangle^{b}}
d\xi_{1}^{\prime}d\mu_{1}^{\prime}d\tau_{1}^{\prime}d\xi^{\prime}
d\mu ^{\prime} d\tau^{\prime}\right|\nonumber\\&&
\leq C\|f\|_{L_{\tau\xi\mu}^{2}}\left(\prod\limits_{j=1}^{2}
\|f_{j}\|_{L_{\tau\xi\mu}^{2}}\right)\label{2.021}.
\end{eqnarray}
Obviously, (\ref{2.021}) follows from (\ref{2.017}). By using a
  proof similar to (\ref{2.018}), we obtain that (\ref{2.019})-(\ref{2.020}) are valid.

This ends the proof of Lemma 2.5.

\begin{Lemma}\label{Lemma2.6}
Let $I,J$ be two intervals on the real line and $f:J\longrightarrow \R$ be a smooth function. Then
\begin{eqnarray}
{\rm mes}\left\{ x\in J, f(x)\in I\right\}\leq \frac{{\rm mes}I}{\inf\limits_{\xi \in J}|f^{\prime}(\xi)|}.\label{2.022}
\end{eqnarray}
\end{Lemma}

Lemma 2.6 can be seen in Lemma 3.8 of \cite{MP}.

\begin{Lemma}\label{Lemma2.7}
Let $0<b_{1}<b_{2}<\frac12$. Then, we have that
\begin{eqnarray}
&&\left\|\chi_{I}(\cdot)u\right\|_{X_{b_{1}}^{0,0}}\leq C\left\|u\right\|_{X_{b_{2}}^{0,0}},\label{2.023}\\
&&\left\|\chi_{I}(\cdot)u\right\|_{X_{-b_{2}}^{0,0}}\leq C\left\|u\right\|_{X_{-b_{1}}^{0,0}}.\label{2.024}
\end{eqnarray}
\end{Lemma}

For the proof of Lemma 2.7, we refer the readers to  Lemma 3.1. of \cite{IMEJDE}.

\begin{Lemma}\label{Lemma2.8}
Let $\phi_{\alpha}(\xi)=\xi|\xi|^{\alpha}$, $\xi=\xi_{1}+\xi_{2}$ and $\alpha\geq4$ and
\begin{eqnarray}
r_{\alpha}(\xi,\xi_{1})=\phi_{\alpha}(\xi)-\phi_{\alpha}(\xi_{1})-\phi_{\alpha}(\xi_{2}).\label{2.018}
\end{eqnarray}
Then $|r_{\alpha}(\xi,\xi_{1})|\sim |\xi_{\rm min}||\xi_{\rm max}|^{\alpha}$.
\end{Lemma}

For the proof of Lemma 2.8, we refer the readers to Lemma 3.4 of \cite{Hadac2008}.

\bigskip
\bigskip

\noindent{\large\bf 3. Bilinear estimates}

\setcounter{equation}{0}

 \setcounter{Theorem}{0}

\setcounter{Lemma}{0}

 \setcounter{section}{3}
 In this section, we give the proof of  Lemmas 3.1, 3.2.
 Lemma 3.1 is used to prove Theorems 1.1. Lemma 3.2 in combination with I-method yields Theorems 1.2.
  Lemma 3.3 is used to prove Lemma
 5.1.

 \begin{Lemma}\label{Lemma3.1}
Let $s_{1}\geq-\frac{\alpha-1}{4}+4\alpha\epsilon,s_{2}\geq0$, $b=\frac{1}{2}+\epsilon$, $b^{\prime}=-\frac{1}{2}+2\epsilon$ and $u_{j}\in X_{\frac{1}{2}+\epsilon}^{s_{1},s_{2}}(j=1,2)$.
Then, we have that
\begin{eqnarray}
&&\|\partial_{x}(u_{1}u_{2})\|_{X_{b^{\prime}}^{s_{1},s_{2}}}\leq C
\prod_{j=1}^{2}\|u_{j}\|_{X_{b}^{s_{1},s_{2}}}.\label{3.01}
\end{eqnarray}
\end{Lemma}
\noindent{\bf Proof.}  To prove (\ref{3.01}),  by duality, it suffices to  prove that
\begin{eqnarray}
&&\left|\int_{\SR^{3}}\bar{u}\partial_{x}(u_{1}u_{2})dxdydt\right|\leq
C\|u\|_{X_{-b^{\prime}}^{-s_{1},-s_{2}}}\left(\prod_{j=1}^{2}
\|u_{j}\|_{X_{b}^{s_{1},s_{2}}}\right).\label{3.02}
\end{eqnarray}
for $u\in X_{-b^{\prime}}^{-s_{1},-s_{2}}.$
We define
\begin{eqnarray}
&&F(\xi,\mu,\tau):=\langle\xi\rangle^{-s_{1}}\langle\mu\rangle^{-s_{2}}
\langle \sigma\rangle^{-b^{\prime}}\mathscr{F}u(\xi,\mu,\tau),\nonumber\\&&
F_{j}(\xi_{j},\mu_{j},\tau_{j}):=\langle\xi_{j}\rangle^{s_{1}}\langle\mu_{j}\rangle^{s_{2}}
\langle \sigma_{j}\rangle^{b}
\mathscr{F}u_{j}(\xi_{j},\mu,\tau_{j})(j=1,2),\nonumber\\
&&dV=d\xi_{1}d\mu_{1}d\tau_{1}d\xi d\mu d\tau\label{3.03}
\end{eqnarray}
and
\begin{eqnarray}
D:=\left\{(\xi_1,\mu_{1},\tau_1,\xi,\mu,\tau)\in {\rm R^6},
\xi=\sum_{j=1}^{2}\xi_j,\mu=\sum_{j=1}^{2}\mu_{j},\tau=\sum_{j=1}^{2}\tau_j\right\}.\label{3.04}
\end{eqnarray}
To obtain (\ref{3.02}), from (\ref{3.03}), it suffices to prove that
\begin{eqnarray}
\int_{D}\frac{|\xi|\langle\xi\rangle^{s_{1}}\langle\mu\rangle^{s_{2}}
F(\xi,\mu,\tau)\prod\limits_{j=1}^{2}F_{j}(\xi_{j},\mu_{j},\tau_{j})}{\langle\sigma_{j}\rangle^{-b^{\prime}}
\prod\limits_{j=1}^{2}\langle\xi_{j}\rangle^{s_{1}}\langle\mu_{j}\rangle^{s_{2}}\langle\sigma_{j}\rangle^{b}}dV
\leq C
\|F\|_{L_{\tau\xi\mu}^{2}}\prod_{j=1}^{2}\|F_{j}\|_{L_{\tau\xi\mu}^{2}}.\label{3.05}
\end{eqnarray}
Without loss of generality, by using the symmetry,  we assume that
$|\xi_{1}|\geq |\xi_{2}|$  and   $F(\xi,\mu,\tau)\geq 0,F_j(\xi_{j},\mu_{j},\tau_{j})\geq 0(j=1,2)$
and
\begin{eqnarray*}
D^{*}:=\left\{(\xi_1,\mu_{1},\tau_1,\xi,\mu,\tau)\in D,|\xi_{2}|\leq|\xi_{1}|\right\}.
\end{eqnarray*}
We define
\begin{eqnarray*}
&&\hspace{-0.8cm}\Omega_1=\left\{(\xi_1,\mu_{1},\tau_1,\xi,\mu,\tau)\in D^{*},
 |\xi_2|\leq |\xi_{1}|\leq 2A\right\},\\
&&\hspace{-0.8cm} \Omega_2=\{ (\xi_1,\mu_{1},\tau_1,\xi,\mu,\tau)\in D^{*},
|\xi_1|\geq 2A, |\xi_{1}|\gg|\xi_{2}|,|\xi_{2}|\leq 2A\},\\
&&\hspace{-0.8cm} \Omega_3=\{ (\xi_1,\mu_{1},\tau_1,\xi,\mu,\tau)\in D^{*},
|\xi_1|\geq 2A, |\xi_{1}|\gg|\xi_{2}|,|\xi_{2}|> 2A\},\\
&&\hspace{-0.8cm}\Omega_4=\{(\xi_1,\mu_{1},\tau_1,\xi,\mu,\tau)\in D^{*},
|\xi_{1}|\geq 2A,4|\xi|\leq |\xi_{2}|\sim|\xi_{1}|,|\xi|\leq 2A,\xi_{1}\xi_{2}<0\},\\
&&\hspace{-0.8cm}\Omega_5=\{(\xi_1,\mu_{1},\tau_1,\xi,\mu,\tau)\in D^{*},
|\xi_{1}|\geq 2A,4|\xi|\leq |\xi_{2}|\sim|\xi_{1}|,|\xi|> 2A,\xi_{1}\xi_{2}<0\},\\
&&\hspace{-0.8cm}\Omega_6=\{(\xi_1,\mu_{1},\tau_1,\xi,\mu,\tau)\in D^{*},
 |\xi_{1}|\geq 2A, |\xi_{1}|\sim |\xi_{2}|,\xi_{1}\xi_{2}<0,|\xi|\geq \frac{|\xi_{2}|}{4}\},\\
 &&\hspace{-0.8cm}\Omega_7=\left\{(\xi_1,\mu_{1},\tau_1,\xi,\mu,\tau)\in D^{*},
 |\xi_{1}|\geq 2A, |\xi_{2}|\sim |\xi_{1}|,\xi_{1}\xi_{2}>0\right\}.
\end{eqnarray*}
Obviously, $D^{*}\subset\bigcup\limits_{j=1}^{7}\Omega_{j}.$
We define
\begin{equation}
    K_{1}(\xi_{1},\mu_{1},\tau_{1},\xi,\mu,\tau):=\frac{|\xi|
    \langle\xi\rangle^{s_{1}}\langle\mu\rangle^{s_{2}}}{\langle\sigma_{j}\rangle^{-b^{\prime}}
\prod\limits_{j=1}^{2}\langle\xi_{j}\rangle^{s_{1}}
\langle\mu_{j}\rangle^{s_{2}}\langle\sigma_{j}\rangle^{b}}\label{3.06}
\end{equation}
and
\begin{eqnarray*}
{\rm Int_{j}}:=\int_{\Omega_{j}} K_{1}(\xi_{1},\mu_{1},\tau_{1},\xi,\mu,\tau)F(\xi,\mu,\tau)
\prod_{j=1}^{2}F_{j}(\xi_{j},\mu_{j},\tau_{j})
dV
\end{eqnarray*}
where $1\leq j\leq 7, j\in N$ and $dV$  is defined as in (\ref{3.03}).
Since $s_{2}\geq0$ and $\mu=\sum\limits_{j=1}^{2}\mu_{j}$, we have that
 $\langle \mu\rangle^{s_{2}} \leq \prod\limits_{j=1}^{2}\langle\mu_{j}\rangle^{s_{2}}$,
thus, we have that
\begin{eqnarray}
K_{1}(\xi_{1},\mu_{1},\tau_{1},\xi,\mu,\tau)\leq\frac{|\xi|\langle\xi\rangle^{s_{1}}}
{\langle\sigma\rangle^{-b^{\prime}}
\prod\limits_{j=1}^{2}\langle\xi_{j}\rangle^{s_{1}}\langle\sigma_{j}\rangle^{b}}.\label{3.07}
\end{eqnarray}
We only prove (\ref{3.01}) with $-\frac{\alpha-1}{4}+16\alpha\epsilon\leq s_{1}<0$ since case $s\geq0$ can be easily proved.

\noindent (1). Region $\Omega_{1}.$ In this region $|\xi|\leq |\xi_{1}|+|\xi_{2}|\leq 4A,$ this case can be proved similarly to case
$low+low\longrightarrow low$ of pages 344-345 of Theorem 3.1 in \cite{LX}.

 \noindent (2). Region $\Omega_{2}.$ In this region, we have that $|\xi|\sim |\xi_{1}|$.

\noindent By using the Cauchy-Schwartz inequality with respect to $\xi_{1},\mu_{1},\tau_{1}$,
 from (\ref{3.06}),  we have that
\begin{eqnarray}
&&{\rm Int_{2}}\leq C\int_{\SR^{3}}\frac{|\xi|}{\langle \sigma \rangle^{-b^{\prime}}}
\left(\int_{\SR^{3}}\frac{d\xi_{1}d\mu_{1}d\tau_{1}}{\prod\limits_{j=1}^{2}
\langle\sigma_{j}\rangle^{2b}}\right)^{\frac{1}{2}}\nonumber\\&&\qquad\qquad\qquad \times\left(\int_{\SR^{3}}
\prod\limits_{j=1}^{2}\left|F_{j}(\xi_{j},\mu_{j},\tau_{j})\right|^{2}
d\xi_{1}d\mu_{1}d\tau_{1}\right)^{\frac{1}{2}}F(\xi,\mu,\tau)d\xi d\mu d\tau.\label{3.08}
\end{eqnarray}
By using (\ref{2.02}), we have that
\begin{eqnarray}
&&\hspace{-1.4cm}\frac{|\xi|}{\langle \sigma \rangle^{-b^{\prime}}}\left(\int_{\SR^{3}}
\frac{d\xi_{1}d\mu_{1}d\tau_{1}}{\prod\limits_{j=1}^{2}
\langle\sigma_{j}\rangle^{2b}}\right)^{\frac{1}{2}}
\leq C\frac{|\xi|}{\langle \sigma \rangle^{-b^{\prime}}}\left(\int_{\SR^{2}}\frac{d\xi_{1}d\mu_{1}}{
\langle\tau+\phi(\xi_{1},\mu_{1})+\phi(\xi_{2},\mu_{2})\rangle^{2b}}\right)^{\frac{1}{2}}\label{3.09}.
\end{eqnarray}
Let $\nu=\tau+\phi(\xi_{1},\mu_{1})+\phi(\xi_{2},\mu_{2})$ and
 $\Delta=|\xi|\xi|^{\alpha}-\xi_{1}|\xi_{1}|^{\alpha}-\xi_{2}|\xi_{2}|^{\alpha},$ since $|\xi_{1}|\gg|\xi_{2}|$,
then we have that the absolute value of Jacobian determinant equals
\begin{eqnarray}
&&\hspace{-1cm}\left|\frac{\partial(\Delta,\nu)}{\partial(\xi_{1},\mu_{1})}\right|
=2\left|\frac{\mu_{1}}{\xi_{1}}-\frac{\mu_{2}}{\xi_{2}}\right|
\left|(\alpha+1)(\xi_{1}^{\alpha}-\xi_{2}^{\alpha})\right|\nonumber\\
&&\hspace{-1cm}=2(\alpha+1)\left|\sigma-\nu-\Delta\right|^{\frac{1}{2}}\left|\frac{\xi}{\xi_{1}\xi_{2}}\right|^{\frac{1}{2}}
\left|(\xi_{1}^{\alpha}-\xi_{2}^{\alpha})\right|\sim
\left|\sigma-\nu+\Delta\right|^{\frac{1}{2}}\left|\frac{\xi}{\xi_{1}\xi_{2}}\right|^{\frac{1}{2}}|\xi_{1}|^{\alpha}
.\label{3.010}
\end{eqnarray}
Inserting (\ref{3.010}) into (\ref{3.09}), by using (\ref{2.03}), we have that
\begin{eqnarray}
&&\frac{|\xi|}{\langle \sigma \rangle^{-b^{\prime}}}\left(\int_{\SR^{3}}
\frac{d\xi_{1}d\mu_{1}d\tau_{1}}{\prod\limits_{j=1}^{2}
\langle\sigma_{j}\rangle^{2b}}\right)^{\frac{1}{2}}
\leq C\frac{|\xi|}{\langle \sigma \rangle^{-b^{\prime}}}\left(\int_{\SR^{2}}\frac{d\xi_{1}d\mu_{1}}{
\langle\tau+\phi(\xi_{1},\mu_{1})+\phi(\xi_{2},\mu_{2})\rangle^{2b}}\right)^{\frac{1}{2}}\nonumber\\
&&\leq \frac{C}{|\xi|^{\frac{\alpha}{2}-1}\langle \sigma \rangle^{-b^{\prime}}}
\left(\int_{\SR^{2}}\frac{d\nu d\Delta}{\left|\sigma-\nu-\Delta\right|^{\frac{1}{2}}
\langle\nu\rangle^{2b}}\right)^{\frac{1}{2}}\nonumber\\&&\leq  \frac{C}{|\xi|^{\frac{\alpha}{2}-1}\langle \sigma \rangle^{-b^{\prime}}}
\left(\int_{|\Delta|<20\alpha|\xi|^{\alpha}}\frac{d\Delta}{
\langle\Delta-\sigma\rangle^{\frac{1}{2}}}\right)^{\frac{1}{2}}.
\label{3.011}
\end{eqnarray}
When $|\sigma|<20\alpha|\xi|^{\alpha},$ combining (\ref{3.011}) with (\ref{2.01}), since $\alpha\geq4,$ we have that
\begin{eqnarray}
 \frac{C}{|\xi|^{\frac{\alpha}{2}-1}\langle \sigma \rangle^{-b^{\prime}}}\left(\int_{|\Delta|<20\alpha|\xi|^{\alpha}}\frac{d\Delta}{
\langle\Delta-\sigma\rangle^{\frac{1}{2}}}\right)^{\frac{1}{2}}
\leq \frac{C}{|\xi|^{\frac{\alpha}{4}-1}\langle \sigma \rangle^{-b^{\prime}}}\leq C\label{3.012}.
\end{eqnarray}
When $|\sigma|\geq20\alpha|\xi|^{\alpha},$ from   (\ref{3.011}), since $\alpha\geq4,$ we have that
\begin{eqnarray}
 \frac{C}{|\xi|^{\frac{\alpha}{2}-1}\langle \sigma \rangle^{-b^{\prime}}}\left(\int_{|\Delta|<20\alpha|\xi|^{\alpha}}\frac{d\Delta}{
\langle\Delta-\sigma\rangle^{\frac{1}{2}}}\right)^{\frac{1}{2}}\leq C\frac{|\xi|}
{\langle \sigma \rangle^{-b^{\prime}}}\leq C\label{3.013}.
\end{eqnarray}
Combining (\ref{3.09}) with (\ref{3.010})-(\ref{3.013}), we have that
\begin{eqnarray}
&&\frac{|\xi|}{\langle \sigma \rangle^{-b^{\prime}}}\left(\int_{\SR^{3}}
\frac{d\xi_{1}d\mu_{1}d\tau_{1}}{\prod\limits_{j=1}^{2}
\langle\sigma_{j}\rangle^{2b}}\right)^{\frac{1}{2}}\leq C\label{3.014}.
\end{eqnarray}
Inserting (\ref{3.014}) into (\ref{3.08}), by using the Cauchy-Schwartz
 inequality with respect to $\xi,\mu,\tau$, we have that
\begin{eqnarray}
&&{\rm Int_{2}}\leq C\int_{\SR^{3}}\frac{|\xi|}{\langle \sigma \rangle^{-b^{\prime}}}
\left(\int_{\SR^{3}}\frac{d\xi_{1}d\mu_{1}d\tau_{1}}{\prod\limits_{j=1}^{2}
\langle\sigma_{j}\rangle^{2b}}\right)^{\frac{1}{2}}\nonumber\\&&\qquad\qquad\qquad \times\left(\int_{\SR^{3}}
\prod\limits_{j=1}^{2}\left|F_{j}(\xi_{j},\mu_{j},\tau_{j})\right|^{2}
d\xi_{1}d\mu_{1}d\tau_{1}\right)^{\frac{1}{2}}F(\xi,\mu,\tau)d\xi d\mu d\tau\nonumber\\
&&\leq C\int_{\SR^{3}}\left(\int_{\SR^{3}}
\prod\limits_{j=1}^{2}\left|F_{j}(\xi_{j},\mu_{j},\tau_{j})\right|^{2}
d\xi_{1}d\mu_{1}d\tau_{1}\right)^{\frac{1}{2}}F(\xi,\mu,\tau)d\xi d\mu d\tau\nonumber\\&&
\leq C\|F\|_{L_{\tau\xi\mu}^{2}}\left(\prod\limits_{j=1}^{2}\|F_{j}\|_{L_{\tau\xi\mu}^{2}}\right).\label{3.015}
\end{eqnarray}
 (3). Region $\Omega_{3}.$ In this region, we have that $|\xi|\sim |\xi_{1}|$.
In this region, we consider
\begin{eqnarray}
&&|\sigma-\sigma_{1}-\sigma_{2}|=\left|\xi|\xi|^{\alpha}-\xi_{1}|\xi_{1}|^{\alpha}-\xi_{2}|\xi_{2}|^{\alpha}-\frac{\xi_{1}\xi_{2}}{\xi}
\left|\frac{\mu_{1}}{\xi_{1}}-\frac{\mu_{2}}{\xi_{2}}\right|^{2}\right|\nonumber\\&&\geq
\frac{\left||\xi|^{\alpha}\xi-\xi_{1}|\xi_{1}|^{\alpha}-\xi_{2}|\xi_{2}|^{\alpha}\right|}{\alpha}\label{3.016}
\end{eqnarray}
and
\begin{eqnarray}
&&|\sigma-\sigma_{1}-\sigma_{2}|=\left|\xi|\xi|^{\alpha}-\xi_{1}|\xi_{1}|^{\alpha}-\xi_{2}|\xi_{2}|^{\alpha}-\frac{\xi_{1}\xi_{2}}{\xi}
\left|\frac{\mu_{1}}{\xi_{1}}-\frac{\mu_{2}}{\xi_{2}}\right|^{2}\right|\nonumber\\&&\leq
\frac{\left||\xi|^{\alpha}\xi-\xi_{1}|\xi_{1}|^{\alpha}-\xi_{2}|\xi_{2}|^{\alpha}\right|}{\alpha}\label{3.017},
\end{eqnarray}
respectively.

\noindent
When (\ref{3.016}) is valid, we have that one of the following three cases must occur:
\begin{eqnarray}
&&|\sigma|:={\rm max}\left\{|\sigma|,|\sigma_{1}|,|\sigma_{2}|\right\}
\geq C|\xi_{\rm min}||\xi_{\rm max}|^{\alpha},\label{3.018}\\
&&|\sigma_{1}|:={\rm max}\left\{|\sigma|,|\sigma_{1}|,|\sigma_{2}|\right\}
\geq C|\xi_{\rm min}||\xi_{\rm max}|^{\alpha},\label{3.019}\\
&&|\sigma_{2}|:={\rm max}\left\{|\sigma|,|\sigma_{1}|,|\sigma_{2}|\right\}
\geq C|\xi_{\rm min}||\xi_{\rm max}|^{\alpha}.\label{3.020}
\end{eqnarray}
When (\ref{3.018}) is valid, since  $s_{1}\geq -\frac{\alpha-1}{4}+4\alpha\epsilon$, we have that
\begin{eqnarray}
&&K_{1}(\xi_{1},\mu_{1},\tau_{1},\xi,\mu,\tau)\leq\frac{|\xi|\langle\xi\rangle^{s_{1}}}
{\langle\sigma\rangle^{-b^{\prime}}
\prod\limits_{j=1}^{2}\langle\xi_{j}\rangle^{s_{1}}\langle\sigma_{j}\rangle^{b}}\leq C\frac{|\xi|^{1-\frac{\alpha}{2}+2\alpha\epsilon}|\xi_{2}|^{-s_{1}+b^{\prime}}}
{\prod\limits_{j=1}^{2}\langle\sigma_{j}\rangle^{b}}\nonumber\\&&
\leq C\frac{|\xi_{1}|^{\frac{\alpha}{4}}|\xi_{2}|^{-\frac{1}{2}}}
{\prod\limits_{j=1}^{2}\langle\sigma_{j}\rangle^{b}}
.\label{3.021}
\end{eqnarray}
Thus, combining (\ref{2.017}) with (\ref{3.021}),  we have that
\begin{eqnarray*}
|{\rm Int_{3}}|\leq C\|F\|_{L_{\tau\xi\mu}^{2}}\left(\prod\limits_{j=1}^{2}\|F_{j}\|_{L_{\tau\xi\mu}^{2}}\right).
\end{eqnarray*}
When (\ref{3.019}) is valid, since  $s_{1}\geq-\frac{\alpha-1}{4}+4\alpha\epsilon$
 and $\langle \sigma\rangle^{b^{\prime}}
 \langle \sigma_{1}\rangle^{-b}\leq \langle \sigma\rangle^{-b}\langle
  \sigma_{1}\rangle^{b^{\prime}},$
we have that
\begin{eqnarray}
\hspace{-0.5cm}K_{1}(\xi_{1},\mu_{1},\tau_{1},\xi,\mu,\tau)\leq\frac{|\xi|\langle\xi\rangle^{s_{1}}}
{\langle\sigma\rangle^{-b^{\prime}}
\prod\limits_{j=1}^{2}\langle\xi_{j}\rangle^{s_{1}}\langle\sigma_{j}\rangle^{b}}\leq C\frac{|\xi|^{1-\frac{\alpha}{2}+2\alpha\epsilon}|\xi_{2}|^{-s_{1}+b^{\prime}}}
{\langle\sigma\rangle^{b}\langle\sigma_{2}\rangle^{b}}
\leq C\frac{|\xi|^{-\frac{1}{2}}|\xi_{2}|^{\frac{\alpha}{4}}}
{\langle\sigma\rangle^{b}\langle\sigma_{2}\rangle^{b}}
.\label{3.022}
\end{eqnarray}
Thus, combining (\ref{2.019}) with (\ref{3.022}),  we have that
\begin{eqnarray*}
|{\rm Int_{3}}|\leq C\|F\|_{L_{\tau\xi\mu}^{2}}\left(\prod\limits_{j=1}^{2}\|F_{j}\|_{L_{\tau\xi\mu}^{2}}\right).
\end{eqnarray*}
When (\ref{3.020}) is valid, since  $s_{1}\geq -\frac{\alpha-1}{4}+4\alpha\epsilon$
 and $\langle \sigma\rangle^{b^{\prime}}\langle \sigma_{2}
 \rangle^{-b}\leq \langle \sigma\rangle^{-b}
 \langle \sigma_{2}\rangle^{b^{\prime}},$
we have that
\begin{eqnarray}
\hspace{-0.5cm}K_{1}(\xi_{1},\mu_{1},\tau_{1},\xi,\mu,\tau)\leq\frac{|\xi|\langle\xi\rangle^{s_{1}}}
{\langle\sigma\rangle^{-b^{\prime}}
\prod\limits_{j=1}^{2}\langle\xi_{j}\rangle^{s_{1}}\langle\sigma_{j}\rangle^{b}}
\leq C\frac{|\xi|^{1-\frac{\alpha}{2}+2\alpha\epsilon}|\xi_{1}|^{-s_{1}+b^{\prime}}}
{\langle\sigma_{1}\rangle^{b}\langle\sigma\rangle^{b}}
\leq C\frac{|\xi_{1}|^{-\frac{1}{2}}|\xi|^{\frac{\alpha}{4}}}
{\langle\sigma_{1}\rangle^{b}\langle\sigma\rangle^{b}}
.\label{3.023}
\end{eqnarray}
Thus, combining (\ref{2.018}) with (\ref{3.023}),  we have that
\begin{eqnarray*}
|{\rm Int_{3}}|\leq C\|F\|_{L_{\tau\xi\mu}^{2}}\left(\prod\limits_{j=1}^{2}\|F_{j}\|_{L_{\tau\xi\mu}^{2}}\right).
\end{eqnarray*}
When (\ref{3.017}) is valid, from Lemma 2.8, we have that
\begin{eqnarray}
\left|\frac{\mu_{1}}{\xi_{1}}-\frac{\mu_{2}}{\xi_{2}}\right|\sim |\xi|^{\frac{\alpha}{2}}.\label{3.024}
\end{eqnarray}
We consider $|\sigma|\geq |\xi|^{\frac{2\alpha-1}{2}}$ and $|\sigma|<|\xi|^{\frac{2\alpha-1}{2}}$, respectively.

\noindent
When $|\sigma|\geq |\xi|^{\frac{2\alpha-1}{2}}$, since $-\frac{\alpha-1}{4}+4\alpha\epsilon\leq s_{1}<0$ and $\alpha\geq4$, we have that
\begin{eqnarray}
&&K_{1}(\xi_{1},\mu_{1},\tau_{1},\xi,\mu,\tau)\leq C\frac{|\xi|\langle\xi\rangle^{s_{1}}}
{\langle\sigma\rangle^{-b^{\prime}}
\prod\limits_{j=1}^{2}\langle\xi_{j}\rangle^{s_{1}}\langle\sigma_{j}\rangle^{b}}\leq C\frac{|\xi_{1}|^{\frac{5-2\alpha}{4}+(2\alpha-1)\epsilon}|\xi_{2}|^{\frac{\alpha-1}{4}-16\alpha\epsilon}}
{\prod\limits_{j=1}^{2}\langle\sigma_{j}\rangle^{b}}\nonumber\\&&
\leq C\frac{|\xi_{1}|^{-\frac{1}{2}}|\xi_{2}|^{\frac{\alpha}{4}}}
{\prod\limits_{j=1}^{2}\langle\sigma_{j}\rangle^{b}}
.\label{3.025}
\end{eqnarray}
Thus, combining (\ref{2.017}) with (\ref{3.025}),  we have that
\begin{eqnarray*}
|{\rm Int_{3}}|\leq C\|F\|_{L_{\tau\xi\mu}^{2}}\left(\prod\limits_{j=1}^{2}\|F_{j}\|_{L_{\tau\xi\mu}^{2}}\right).
\end{eqnarray*}
Now we consider case $|\sigma|<|\xi|^{\frac{2\alpha-1}{2}}$. We dyadically decompose with respect to
\begin{eqnarray*}
\langle\sigma\rangle\sim 2^{j},\langle\sigma_{1}\rangle\sim 2^{j_{1}},
\langle\sigma_{2}\rangle\sim 2^{j_{2}},|\xi|\sim 2^{m},|\xi_{1}|\sim 2^{m_{1}},|\xi_{2}|\sim 2^{m_{2}}.
\end{eqnarray*}
Let  $D_{j,j_{1},j_{2},m,m_{1},m_{2}}$ be the image of set of all points
$(\xi_{1},\mu_{1},\tau_{1},\xi,\mu,\tau)\in D^{*}$ satisfying
\begin{eqnarray}
&&|\xi_1|\geq 2A, |\xi_{1}|\gg|\xi_{2}|,|\xi_{2}|> 2A,|\xi|\sim 2^{m},|\xi_{1}|\sim 2^{m_{1}},|\xi_{2}|\sim 2^{m_{2}}, \nonumber\\
&&\langle\sigma\rangle\sim 2^{j}\leq 2^{\frac{2\alpha-1}{2}m},\langle\sigma_{1}\rangle\sim 2^{j_{1}},
\langle\sigma_{2}\rangle\sim 2^{j_{2}}\label{3.026}
\end{eqnarray}
under the transformation $(\xi_{1},\mu_{1},\tau_{1},\xi_{2},\mu_{2},\tau_{2})\longrightarrow
(\xi_{1},\mu_{1},\sigma_{1},\xi_{2},\mu_{2},\sigma_{2})$. We define
\begin{eqnarray}
&&f_{m_{k},j_k}:=\left|\eta_{m_{k}}(\xi_{k})\eta_{j_{k}}(\sigma_{k})F_{k}(\xi_{k},\mu_{k},\tau_{k})\right|(k=1,2),\label{3.027}\\
&&f_{m,j}:=\left|\eta_{m}(\xi)\eta_{j}(\sigma)F(\xi,\mu,\sigma_{1}-\phi(\xi_{1},\mu_{1})+\sigma_{2}-\phi(\xi_{2},\mu_{2}))\right|,\label{3.028}\\
&&P:=\left|f_{m,j}\prod\limits_{k=1}^{2}f_{m_{k},j_k}\right|,dV=d\xi_{1}d\mu_{1}d\sigma_{1}d\xi_{2}d\mu_2d\sigma_{2}.\label{3.029}
\end{eqnarray}
Thus, we have that
\begin{eqnarray}
\hspace{-0.6cm}{\rm Int_{3}}\leq  C\sum\limits_{\>m,m_{1},m_{2}>0}\sum\limits_{j_{1},j_{2}\geq0,\>0<j\leq \frac{(2\alpha-1)m}{2}}
\int_{D_{j,j_{1},j_{2},m,m_{1},m_{2}}}
2^{jb^{\prime}-(j_{1}+j_{2})b-m_{2}s_{1}+m}PdV
.\label{3.030}
\end{eqnarray}
In this case, we consider
\begin{eqnarray}
&&\left|(\alpha+1)(|\xi_{1}|^{\alpha}-|\xi_{2}|^{\alpha})
-\left[\left(\frac{\mu_{1}}{\xi_{1}}\right)^{2}
-\left(\frac{\mu_{2}}{\xi_{2}}\right)^{2}\right]\right|>2^{j+\frac{(\alpha-1)m_{1}}{2}},\label{3.031}\\
&&\left|(\alpha+1)(|\xi_{1}|^{\alpha}-|\xi_{2}|^{\alpha})
-\left[\left(\frac{\mu_{1}}{\xi_{1}}\right)^{2}
-\left(\frac{\mu_{2}}{\xi_{2}}\right)^{2}\right]\right|\leq 2^{j+\frac{(\alpha-1)m_{1}}{2}},\label{3.032}
\end{eqnarray}
respectively.

\noindent Now  we consider case (\ref{3.031}).   We make the change of variables
\begin{eqnarray}
u=\xi_{1}+\xi_{2},v=\mu_{1}+\mu_{2},w=\sigma_{1}
-\phi(\xi_{1},\mu_{1})+\sigma_{2}-\phi(\xi_{2},\mu_{2}),\mu_{2}=\mu_{2},\label{3.033}
\end{eqnarray}
thus the Jacobian determinant equals
\begin{eqnarray}
\frac{\partial(u,v,w,\mu_{2})}{\partial(\xi_{1},\xi_{2},\mu_{1},\mu_{2})}
=(\alpha+1)(|\xi_{1}|^{\alpha}-|\xi_{2}|^{\alpha})
-\left[\left(\frac{\mu_{1}}{\xi_{1}}\right)^{2}-\left(\frac{\mu_{2}}{\xi_{2}}\right)^{2}\right].\label{3.034}
\end{eqnarray}
We assume that $D_{j,j_{1},j_{2},m,m_{1},m_{2}}^{(1)}$
 is the image of the subset of all points
$$(\xi_{1},\mu_{1},\sigma_{1},\xi_{2},\mu_{2},\sigma_{2})\in D_{j,j_{1},j_{2},m,m_{1},m_{2}},$$
which satisfies (\ref{3.031}) under the transformation
(\ref{3.033}). Combining (\ref{3.031}) with (\ref{3.034}), we have that
\begin{eqnarray}
\left|\frac{\partial(u,v,w,\mu_{2})}{\partial(\xi_{1},\xi_{2},\mu_{1},\mu_{2})}\right|>2^{j+\frac{(\alpha-1)m_{1}}{2}}.\label{3.035}
\end{eqnarray}
Let $G_{1}(u,v,w,\mu_{2},\sigma_{1},\sigma_{2})$ be $\eta_{m}(\xi)\eta_{j}(\sigma)\prod\limits_{k=1}^{2}f_{m_{k},j_{k}}$
 under the change of the variables (\ref{3.033}) and
\begin{eqnarray}
M_{1}=\left|F(u,v,w)G_{1}(u,v,w,\mu_{2},\sigma_{1},\sigma_{2})\right|,dV^{(1)}=dudvdwd\mu_{2}d\sigma_{1}d\sigma_{2}
\label{3.036}.
\end{eqnarray}
Thus, (\ref{3.030}) can be controlled by
\begin{eqnarray}
\hspace{-0.5cm}C\sum\limits_{m,m_{1},m_{2}>0}\sum\limits_{j_{1},j_{2}\geq0,\>0<j\leq \frac{7m}{2}}\int_{D_{j,j_{1},j_{2},m,m_{1},m_{2}}^{(1)}}
2^{jb^{\prime}-(j_{1}+j_{2})b-m_{2}s_{1}+m}\frac{M_{1}dV^{(1)}}
{\left|\frac{\partial(u,v,w,\mu_{2})}{\partial(\xi_{1},\xi_{2},\mu_{1},\mu_{2})}\right|}\label{3.037}.
\end{eqnarray}
Inspired by \cite{CKS-GAFA, LX},
we define
\begin{eqnarray}
&&f(\mu):=\sigma_{1}+\sigma_{2}-
(\xi|\xi|^{\alpha}-\xi_{1}|\xi_{1}|^{\alpha}-\xi_{2}|\xi_{2}|^{\alpha})+\frac{\xi_{1}\xi_{2}}{\xi}
\left[\frac{\mu_{1}}{\xi_{1}}-\frac{\mu}{\xi_{2}}\right]^{2}\label{3.038}.
\end{eqnarray}
From (\ref{3.038}) and (\ref{3.024}), for fixed $\sigma_{1},\sigma_{2},\xi_{1},\xi_{2},\mu_{1}$,
 we have that
\begin{eqnarray}
&&\hspace{-0.5cm}|f(\mu_{2})|:=|\sigma_{1}+\sigma_{2}+\phi(\xi_{1},\mu_{1})+\phi(\xi_{2},\mu_{2})
-\phi(\xi,\mu)|=|\tau-\phi(\xi,\mu)|\sim 2^{j},\label{3.039}\\&&|f^{\prime}(\mu_{2})|\sim
 \left|\frac{\mu_{1}}{\xi_{1}}-\frac{\mu_{2}}{\xi_{2}}\right|\sim 2^{\frac{\alpha m}{2}} .\label{3.040}
\end{eqnarray}
For fixed $\sigma_{1},\sigma_{2},\xi_{1},\xi_{2},\mu_{1}$, combining (\ref{3.039}), (\ref{3.040})
 with Lemma 2.6,  we have that
the Lebesgue measure of $\mu_{2}$ can be controlled by $C2^{j-\frac{\alpha m}{2}}$.
By using the Cauchy-Schwartz inequality with respect to $\mu_{2}$
  and the Cauchy-Schwartz inequality with respect to $u,v,w$ and the inverse change of variables related to (\ref{3.033})
  and the Cauchy-Schwartz inequality with respect to $\sigma_{1}$
   and $\sigma_{2}$, since $s\geq -\frac{\alpha-1}{4}+4\alpha\epsilon,$ we have that (\ref{3.037})
     can be  bounded by
\begin{eqnarray*}
&&\hspace{-0.5cm}C\sum
\int_{D_{j,j_{1},j_{2},m,m_{1},m_{2}}^{(1)}}
2^{jb^{\prime}-(j_{1}+j_{2})b-m_{2}s_{1}+m}\frac{M_{1}dV^{(1)}}
{\left|\frac{\partial(u,v,w,\mu_{2})}{\partial(\xi_{1},\xi_{2},\mu_{1},\mu_{2})}\right|}\nonumber\\
&&\leq C\sum2^{2j
\epsilon-(j_{1}+j_{2})b-m(s_{1}+\frac{\alpha-4}{4})}\int F(u,v,w)
\left(\int\frac{G^{2}(u,v,w,\mu_{2},\sigma_{1},\sigma_{2})}{
\left|\frac{\partial(u,v,w,\mu_{2})}{\partial(\xi_{1},\xi_{2},\mu_{1},\mu_{2})}\right|^{2}}
d\mu_{2}\right)^{\frac{1}{2}}dV^{(2)}\nonumber\\
&&\leq C\sum2^{jb^{\prime}-(j_{1}+j_{2})b-m_{1}(s_{1}+\frac{2\alpha-5}{4})}\|F\|_{L^{2}}\int
\left(\int\frac{G^{2}(u,v,w,\mu_{2},\sigma_{1},\sigma_{2})}{
\left|\frac{\partial(u,v,w,\mu_{2})}{\partial(\xi_{1},\xi_{2},\mu_{1},\mu_{2})}\right|}
dV^{(3)}\right)^{\frac{1}{2}}
d\sigma_{1}d\sigma_{2}\nonumber\\
&&\leq C\sum
2^{jb^{\prime}-(j_{1}+j_{2})b-16\alpha m_{1}\epsilon }
\|F\|_{L^{2}}\int
\left(\int \prod\limits_{k=1}^{2}f_{m_{k},j_{k}}^{2}d\xi_{1}d\mu_{1}d\xi_{2}d\mu_{2}\right)^{\frac{1}{2}}
d\sigma_{1}d\sigma_{2}\nonumber\\&&
\leq C\sum
2^{jb^{\prime}-(j_{1}+j_{2})\epsilon-16\alpha m_{1}\epsilon}
\|F\|_{L^{2}}
\left(\int \prod\limits_{k=1}^{2}f_{m_{k},j_{k}}^{2}dV\right)^{\frac{1}{2}}
\leq C\|F\|_{L^{2}}\left(\prod\limits_{j=1}^{2}\|F_{j}\|_{L^{2}}\right),
\end{eqnarray*}
where
\begin{eqnarray}
dV^{(2)}=dudvdwd\sigma_{1}d\sigma_{2}, dV^{(3)}=dudvdwd\mu_{2},\sum=\sum\limits_{m, m_{1},m_{2}>0}\sum\limits_{j_{1},j_{2}\geq0,\>0<j\leq \frac{(2\alpha-1)m}{2}}.\label{3.041}
\end{eqnarray}
Now we consider (\ref{3.032}).
We make the change of variables
\begin{eqnarray}
u=\xi_{1}+\xi_{2},v=\mu_{1}+\mu_{2},w=\sigma_{1}
+\phi(\xi_{1},\mu_{1})+\sigma_{2}+\phi(\xi_{2},\mu_{2}),\xi_{1}=\xi_{1}.\label{3.042}
\end{eqnarray}
 From (\ref{3.042}) and (\ref{3.024}), we have that the absolute value of  the Jacobian determinant equals
\begin{eqnarray}
\left|\frac{\partial(u,v,w,\xi_{1})}{\partial(\xi_{1},\xi_{2},\mu_{1},\mu_{2})}\right|
=2\left|\frac{\mu_{1}}{\xi_{1}}-\frac{\mu_{2}}{\xi_{2}}\right|\sim 2^{\frac{\alpha m}{2}}.\label{3.043}
\end{eqnarray}
We assume that $D_{j,j_{1},j_{2},m,m_{1},m_{2}}^{(2)}$ is the image of the subset of all points
$$(\xi_{1},\mu_{1},\sigma_{1},\xi_{2},\mu_{2},\sigma_{2})\in D_{j,j_{1},j_{2},m,m_{1},m_{2}},$$
 which satisfies (\ref{3.032}) under the transformation
(\ref{3.042}).
Let
$$H_{1}(u,v,w,\xi_{1},\sigma_{1},\sigma_{2})$$
be $\eta_{m}(\xi)\eta_{j}(\sigma)\prod\limits_{k=1}^{2}f_{m_{k},j_{k}}$
under the change of the variables as in (\ref{3.042}) and
\begin{eqnarray}
M_{2}=\left|F(u,v,w)H_{1}(u,v,w,\xi_{1},\sigma_{1},\sigma_{2})\right|, dV^{(4)}
=dudvdwd\xi_{1}d\sigma_{1}d\sigma_{2}.\label{3.044}
\end{eqnarray}
Thus, (\ref{3.030}) can be controlled by
\begin{eqnarray}
\hspace{-0.5cm}C\sum\limits_{m, m_{1},m_{2}>0}\sum\limits_{j_{1},j_{2}\geq0,\>0<j\leq \frac{7m}{2}}
\int_{D_{j,j_{1},j_{2},m,m_{1},m_{2}}^{(2)}}
2^{jb^{\prime}-(j_{1}+j_{2})b-m_{2}s_{1}+m}\frac{M_{2}dV^{(4)}}
{\left|\frac{\partial(u,v,w,\xi_{1})}{\partial(\xi_{1},\xi_{2},\mu_{1},\mu_{2})}\right|}\label{3.045}.
\end{eqnarray}
Inspired by \cite{CKS-GAFA,LX}, we define
\begin{eqnarray}
h(\xi):=(\alpha+1)(|\xi|^{\alpha}-|\xi_{2}|^{\alpha})
-\left[\left(\frac{\mu_{1}}{\xi}\right)^{2}-\left(\frac{\mu_{2}}{\xi_{2}}\right)^{2}\right]\label{3.046},
\end{eqnarray}
for fixed $\xi_{2},\mu_{1},\mu_{2}$, from (\ref{3.046}), we have that
\begin{eqnarray}
&&|h^{\prime}(\xi_{1})|=\left|\alpha(\alpha+1)\xi_{1}|\xi_{1}|^{\alpha-2}+2
\left(\frac{\mu_{1}}{\xi_{1}}\right)^{2}\xi_{1}\right|\geq \alpha(\alpha+1)|\xi_{1}|^{\alpha-1}\geq C2^{(\alpha-1)m_{1}},\nonumber\\&&
|h(\xi_{1})|\leq C2^{j+\frac{(\alpha-1)m_{1}}{2}}\label{3.047}.
\end{eqnarray}
For fixed $\xi_{2},\mu_{1},\mu_{2},$ combining  (\ref{3.047}) with Lemma 2.6,  we have that
the Lebesgue measure of $\xi_{1}$ can be controlled by
 $C2^{j-\frac{(\alpha-1)m_{1}}{2}}$.
By using the Cauchy-Schwartz inequality with respect to $\xi_{1}$
  and the Cauchy-Schwartz inequality with respect to $u,v,w$ and the inverse change of variables related to (\ref{3.042})
  and the Cauchy-Schwartz inequality with respect to $\sigma_{1}$ and $\sigma_{2}$, since $s\geq-\frac{\alpha-1}{4}+4\alpha\epsilon,$
   we have that (\ref{3.045}) can be bounded by
\begin{eqnarray}
&&\hspace{-1cm}C\sum
\int_{D_{j,j_{1},j_{2},m,m_{1},m_{2}}^{(2)}}
2^{jb^{\prime}-(j_{1}+j_{2})b-m_{2}s_{1}+m}\frac{M_{2}dV^{(4)}}
{\left|\frac{\partial(u,v,w,\xi_{1})}{\partial(\xi_{1},\xi_{2},\mu_{1},\mu_{2})}\right|}\nonumber\\
&&\hspace{-1cm}\leq C\sum2^{2j\epsilon-(j_{1}+j_{2})b+m(-s_{1}+\frac{5-\alpha}{4})}
\int |F|
\left(\int\frac{H^{2}(u,v,w,\xi_{1},\sigma_{1},\sigma_{2})}{
\left|\frac{\partial(u,v,w,\xi_{1})}{\partial(\xi_{1},\xi_{2},\mu_{1},\mu_{2})}
\right|^{2}}d\xi_{1}\right)^{\frac{1}{2}}dV^{(2)}\nonumber\\
&&\hspace{-1cm}\leq C\sum2^{2j\epsilon-(j_{1}+j_{2})b
+m(-s_{1}+\frac{5-2\alpha}{4})\epsilon}\|F\|_{L^{2}}\int
\left(\int\frac{H^{2}(u,v,w,\xi_{1},\sigma_{1},\sigma_{2})}{
\left|\frac{\partial(u,v,w,\xi_{1})}{\partial(\xi_{1},\xi_{2},\mu_{1},\mu_{2})}\right|}
dudvdwd\xi_{1}\right)^{\frac{1}{2}}
d\sigma_{1}d\sigma_{2}\nonumber\\
&&\hspace{-1cm}\leq C\sum2^{2j\epsilon-(j_{1}+j_{2})b- 4\alpha m\epsilon}
\|F\|_{L^{2}}\int
\left(\int \prod\limits_{k=1}^{2}f_{m_{k},j_{k}}^{2}d\xi_{1}d\mu_{1}d\xi_{2}d\mu_{2}\right)^{\frac{1}{2}}
d\sigma_{1}d\sigma_{2}\nonumber\\
&&\hspace{-1cm}\leq C\sum2^{2j\epsilon-(j_{1}+j_{2})
\epsilon-4\alpha m\epsilon}\|F\|_{L^{2}}
\left(\int \prod\limits_{k=1}^{2}f_{m_{k},j_{k}}^{2}dV\right)^{\frac{1}{2}}
\leq C\|F\|_{L^{2}}\left(\prod\limits_{j=1}^{2}\|F_{j}\|_{L^{2}}\right),\label{3.048}
\end{eqnarray}
where $\sum$ is defined as in (\ref{3.041}).

\noindent (4). Region $\Omega_{4}.$
In this case, we consider (\ref{3.016}), (\ref{3.017}), respectively.

\noindent When (\ref{3.016}) is valid,  one of (\ref{3.018})-(\ref{3.020}) must occur.

\noindent When (\ref{3.018}) is valid, since  $s_{1}\geq -\frac{\alpha-1}{4}+4\alpha\epsilon$, we have that
\begin{eqnarray}
\hspace{-0.5cm}K_{1}(\xi_{1},\mu_{1},\tau_{1},\xi,\mu,\tau)\leq\frac{|\xi|\langle\xi\rangle^{s_{1}}}
{\langle\sigma\rangle^{-b^{\prime}}
\prod\limits_{j=1}^{2}\langle\xi_{j}\rangle^{s_{1}}\langle\sigma_{j}\rangle^{b}}
\leq C\frac{|\xi|^{\frac{1}{2}+2\epsilon}|\xi_{2}|^{-2s_{1}-2+8\epsilon}}
{\prod\limits_{j=1}^{2}\langle\sigma_{j}\rangle^{b}}
\leq C\frac{|\xi_{1}|^{-\frac{1}{2}}|\xi_{2}|}
{\prod\limits_{j=1}^{2}\langle\sigma_{j}\rangle^{b}}
.\label{3.049}
\end{eqnarray}
Thus, combining (\ref{2.017}) with (\ref{3.049}),  we have that
\begin{eqnarray*}
|{\rm Int_{4}}|\leq C\|F\|_{L_{\tau\xi\mu}^{2}}\left(\prod\limits_{j=1}^{2}\|F_{j}\|_{L_{\tau\xi\mu}^{2}}\right).
\end{eqnarray*}
When (\ref{3.019}) is valid, since  $s_{1}\geq -\frac{\alpha-1}{4}+4\alpha\epsilon$ and
 $\langle \sigma\rangle^{b^{\prime}}\langle \sigma_{1}\rangle^{-b}
 \leq \langle \sigma\rangle^{-b}\langle \sigma_{1}\rangle^{b^{\prime}},$
we have that
\begin{eqnarray}
\hspace{-0.5cm}K_{1}(\xi_{1},\mu_{1},\tau_{1},\xi,\mu,\tau)\leq\frac{|\xi|\langle\xi\rangle^{s_{1}}}
{\langle\sigma\rangle^{-b^{\prime}}
\prod\limits_{j=1}^{2}\langle\xi_{j}\rangle^{s_{1}}\langle\sigma_{j}\rangle^{b}}
\leq C\frac{|\xi|^{\frac{1}{2}+2\epsilon}|\xi_{2}|^{-2s_{1}-2+8\epsilon}}
{\langle\sigma\rangle^{b}\langle\sigma_{2}\rangle^{b}}
\leq C\frac{|\xi|^{-\frac{1}{2}}|\xi_{2}|}
{\langle\sigma\rangle^{b}\langle\sigma_{2}\rangle^{b}}
.\label{3.050}
\end{eqnarray}
Thus, combining (\ref{2.019}) with (\ref{3.050}),  we have that
\begin{eqnarray*}
|{\rm Int_{4}}|\leq C\|F\|_{L_{\tau\xi\mu}^{2}}\left(\prod\limits_{j=1}^{2}\|F_{j}\|_{L_{\tau\xi\mu}^{2}}\right).
\end{eqnarray*}
When (\ref{3.020}) is valid, this case can be proved similarly to (\ref{3.019}).

\noindent  When (\ref{3.017}) is valid, from Lemma 2.8, we have that
\begin{eqnarray}
\left|\frac{\mu_{1}}{\xi_{1}}-\frac{\mu_{2}}{\xi_{2}}\right|\sim |\xi||\xi_{1}|^{\frac{\alpha-2}{2}}.\label{3.051}
\end{eqnarray}
We consider $|\sigma|\geq |\xi_{1}|^{\frac{2\alpha-1}{2}}$ and $|\sigma|<|\xi_{1}|^{\frac{2\alpha-1}{2}}$, respectively.

\noindent
When $|\sigma|\geq |\xi_{1}|^{\frac{2\alpha-1}{2}}$, since $s_{1}\geq-\frac{\alpha-1}{4}+4\alpha\epsilon$, we have that
\begin{eqnarray*}
&&K_{1}(\xi_{1},\mu_{1},\tau_{1},\xi,\mu,\tau)\leq\frac{|\xi|\langle\xi\rangle^{s_{1}}}
{\langle\sigma\rangle^{-b^{\prime}}
\prod\limits_{j=1}^{2}\langle\xi_{j}\rangle^{s_{1}}\langle\sigma_{j}
\rangle^{b}}\leq C\frac{|\xi||\xi_{2}|^{-2s_{1}-\frac{2\alpha-1}{4}+(2\alpha-1)\epsilon}}
{\prod\limits_{j=1}^{2}\langle\sigma_{j}\rangle^{b}}\nonumber\\&&
\leq C\frac{|\xi_{1}|^{-\frac{1}{2}}|\xi_{2}|^{\frac{\alpha}{4}}}
{\prod\limits_{j=1}^{2}\langle\sigma_{j}\rangle^{b}}
.
\end{eqnarray*}
This case can be proved similarly to (\ref{3.049}).

\noindent We consider case $|\sigma|<|\xi_{1}|^{\frac{2\alpha-1}{2}}$.
In this case, we consider
\begin{eqnarray}
&&\left|(\alpha+1)(|\xi_{1}|^{\alpha}-|\xi_{2}|^{\alpha})
-\left[\left(\frac{\mu_{1}}{\xi_{1}}\right)^{2}
-\left(\frac{\mu_{2}}{\xi_{2}}\right)^{2}\right]\right|>2^{j+\frac{(\alpha-2)m_{1}}{2}},\label{3.052}\\
&&\left|(\alpha+1)(|\xi_{1}|^{\alpha}-|\xi_{2}|^{\alpha})
-\left[\left(\frac{\mu_{1}}{\xi_{1}}\right)^{2}
-\left(\frac{\mu_{2}}{\xi_{2}}\right)^{2}\right]\right|\leq 2^{j+\frac{(\alpha-2)m_{1}}{2}},\label{3.053}
\end{eqnarray}
respectively.
We dyadically decompose with respect to
\begin{eqnarray*}
\langle\sigma\rangle\sim 2^{j},\langle\sigma_{1}\rangle\sim 2^{j_{1}},
\langle\sigma_{2}\rangle\sim 2^{j_{2}},|\xi|\sim 2^{m},|\xi_{1}|\sim 2^{m_{1}},|\xi_{2}|\sim 2^{m_{2}}.
\end{eqnarray*}
Let  $D_{j,j_{1},j_{2},m,m_{1},m_{2}}^{(3)}$ be the image of set of all points
$(\xi_{1},\mu_{1},\tau_{1},\xi,\mu,\tau)\in D^{*}$ satisfying
\begin{eqnarray}
&&|\xi_{1}|\geq 2A,4|\xi|\leq |\xi_{2}|\sim|\xi_{1}|,|\xi|\leq 2A,\xi_{1}\xi_{2}<0,
|\xi|\sim 2^{m},|\xi_{1}|\sim 2^{m_{1}},|\xi_{2}|\sim 2^{m_{2}},
\nonumber\\
&&\langle\sigma\rangle\sim 2^{j}\leq 2^{\frac{(2\alpha-1)m_{1}}{2}},\langle\sigma_{1}\rangle\sim 2^{j_{1}},
\langle\sigma_{2}\rangle\sim 2^{j_{2}}\label{3.054}
\end{eqnarray}
under the transformation $(\xi_{1},\mu_{1},\tau_{1},\xi_{2},\mu_{2},\tau_{2})\longrightarrow
(\xi_{1},\mu_{1},\sigma_{1},\xi_{2},\mu_{2},\sigma_{2})$.

\noindent Thus, we have
\begin{eqnarray}
\hspace{-0.8cm}{\rm Int_{4}}\leq C\sum\limits_{m_{1}, m_{2}>0,\>m}
\sum\limits_{j_{1},j_{2}\geq0,\>0<j\leq \frac{(2\alpha-1)m_{1}}{2}}\int_{D_{j,j_{1},j_{2},m,m_{1},m_{2}}^{(3)}}
2^{jb^{\prime}-(j_{1}+j_{2})b-2m_{2}s_{1}+m}PdV
,\label{3.055}
\end{eqnarray}
where $P$ and $dV$ are defined as in (\ref{3.029}).

\noindent In this case, we consider (\ref{3.052}), (\ref{3.053}), respectively.

\noindent When (\ref{3.052}) is valid, we make the change of variables as in (\ref{3.033}).

\noindent Thus the Jacobian determinant equals
\begin{eqnarray}
\frac{\partial(u,v,w,\mu_{2})}{\partial(\xi_{1},\xi_{2},\mu_{1},\mu_{2})}=(\alpha+1)(|\xi_{1}|^{\alpha}-|\xi_{2}|^{\alpha})
-\left[\left(\frac{\mu_{1}}{\xi_{1}}\right)^{2}-\left(\frac{\mu_{2}}{\xi_{2}}\right)^{2}\right]\label{3.056}.
\end{eqnarray}
We assume that $D_{j,j_{1},j_{2},m,m_{1},m_{2}}^{(4)}$ is the image of the subset of all points
$$(\xi_{1},\mu_{1},\sigma_{1},\xi_{2},\mu_{2},\sigma_{2})\in D_{j,j_{1},j_{2},m,m_{1},m_{2}}^{(3)},$$
which satisfies (\ref{3.052}) under the transformation
(\ref{3.033}). Combining (\ref{3.056}) with (\ref{3.052}), we have that
\begin{eqnarray}
\left|\frac{\partial(u,v,w,\mu_{2})}{\partial(\xi_{1},\xi_{2},\mu_{1},\mu_{2})}\right|>2^{j+\frac{(\alpha-2)m_{1}}{2}}\label{3.057}.
\end{eqnarray}
Let $G_{2}(u,v,w,\mu_{2},\sigma_{1},\sigma_{2})$ be $\eta_{m}(\xi)\eta_{j}(\sigma)\prod\limits_{k=1}^{2}f_{m_{k},j_{k}}$
 under the change of the variables (\ref{3.033}) and
\begin{eqnarray}
M_{3}=F(u,v,w)G_{2}(u,v,w,\mu_{2},\sigma_{1},\sigma_{2}),
dV^{(1)}=dudvdwd\mu_{2}d\sigma_{1}d\sigma_{2}\label{3.058}.
\end{eqnarray}
Thus, (\ref{3.055}) can be controlled by
\begin{eqnarray}
\hspace{-0.5cm}C\sum\limits_{m_{1}, m_{2}>0,\>m}\sum\limits_{j_{1},j_{2}\geq0,\>0<j\leq
\frac{(2\alpha-1)m_{1}}{2}}\int_{D_{j,j_{1},j_{2},m,m_{1},m_{2}}^{(4)}}
2^{jb^{\prime}-(j_{1}+j_{2})b-2m_{1}s_{1}+m}\frac{M_{3}dV^{(1)}}
{\left|\frac{\partial(u,v,w,\mu_{2})}{\partial(\xi_{1},\xi_{2},\mu_{1},\mu_{2})}\right|}\label{3.059}.
\end{eqnarray}
Inspired by \cite{CKS-GAFA, LX},
we define
\begin{eqnarray}
&&f(\mu):=\sigma_{1}+\sigma_{2}-
(\xi|\xi|^{\alpha}-\xi_{1}|\xi_{1}|^{\alpha}-\xi_{2}|\xi_{2}|^{\alpha})+\frac{\xi_{1}\xi_{2}}{\xi}
\left[\frac{\mu_{1}}{\xi_{1}}-\frac{\mu}{\xi_{2}}\right]^{2}\label{3.060}.
\end{eqnarray}
From (\ref{3.060}) and (\ref{3.051}), for fixed $\xi_{1},\xi_{2},\mu_{1},\sigma_{1},\sigma_{2},$ we have that
\begin{eqnarray}
&&\hspace{-0.5cm}|f(\mu_{2})|=|\sigma_{1}+\sigma_{2}+\phi(\xi_{1},\mu_{1})+\phi(\xi_{2},\mu_{2})-\phi(\xi,\mu)|
=|\tau-\phi(\xi,\mu)|\sim 2^{j},\label{3.061}\\&&|f^{\prime}(\mu_{2})|\sim \left|\frac{\xi_{1}}{\xi}\right|
\left|\frac{\mu_{1}}{\xi_{1}}-\frac{\mu_{2}}{\xi_{2}}\right|\sim 2^{\frac{\alpha m_{1}}{2}} .\label{3.062}
\end{eqnarray}
Combining (\ref{3.061}), (\ref{3.062}) with Lemma 2.6, for fixed $\xi_{1},\xi_{2},\mu_{1},\sigma_{1},\sigma_{2},$
 we have that
the Lebesgue measure of $\mu_{2}$ can be controlled by $C2^{j-\frac{\alpha m_{1}}{2}}$.
By using the Cauchy-Schwartz inequality with respect to $\mu_{2}$
 and the inverse change of variables related to (\ref{3.029}) and the Cauchy-Schwartz inequality with respect to $u,v,w$
  and the Cauchy-Schwartz inequality with respect to $\sigma_{1}$
   and $\sigma_{2}$, we have that (\ref{3.059}) can be controlled by
\begin{eqnarray}
&&\hspace{-0.5cm}C\sum
\int_{D_{j,j_{1},j_{2},m,m_{1},m_{2}}^{(4)}}
2^{jb^{\prime}-(j_{1}+j_{2})b-2m_{1}s_{1}+m}\frac{M_{3}dV^{(1)}}
{\left|\frac{\partial(u,v,w,\mu_{2})}{\partial(\xi_{1},\xi_{2},\mu_{1},\mu_{2})}\right|}\nonumber\\
&&\leq C\sum2^{2j\epsilon-(j_{1}+j_{2})b-(\frac{\alpha }{4}+2s_{1})m_{1}+m}
\int F
\left(\int\frac{G_{2}^{2}(u,v,w,\mu_{2},\sigma_{1},\sigma_{2})}{
\left|\frac{\partial(u,v,w,\mu_{2})}{\partial(\xi_{1},\xi_{2},\mu_{1},\mu_{2})}\right|^{2}}
d\mu_{2}\right)^{\frac{1}{2}}dV^{(2)}\nonumber\\
&&\leq C\sum2^{jb^{\prime}-(j_{1}+j_{2})b
-(\frac{\alpha-1}{2}+2s_{1}) m_{1}+m}\|F\|_{L^{2}}\int
\left(\int\frac{G_{2}^{2}(u,v,w,\mu_{2},\sigma_{1},\sigma_{2})}{
\left|\frac{\partial(u,v,w,\mu_{2})}{\partial(\xi_{1},\xi_{2},\mu_{1},\mu_{2})}\right|}
dV^{(3)}\right)^{\frac{1}{2}}
d\sigma_{1}d\sigma_{2}\nonumber\\
&&\leq C\sum
2^{jb^{\prime}-(j_{1}+j_{2})b-32\alpha m_{1}\epsilon +m}
\|F\|_{L^{2}}\int
\left(\int \prod\limits_{k=1}^{2}f_{m_{k},j_{k}}^{2}d\xi_{1}d\mu_{1}d\xi_{2}d\mu_{2}\right)^{\frac{1}{2}}
d\sigma_{1}d\sigma_{2}\nonumber\\&&
\leq C\sum
2^{jb^{\prime}-(j_{1}+j_{2})\epsilon-32\alpha m_{1}\epsilon+m}
\|F\|_{L^{2}}
\left(\int \prod\limits_{k=1}^{2}f_{m_{k},j_{k}}^{2}dV\right)^{\frac{1}{2}}
\leq C\|F\|_{L^{2}}\prod\limits_{j=1}^{2}\|F_{j}\|_{L^{2}}.\label{3.063}
\end{eqnarray}
Here $\sum=\sum\limits_{m_{1}, m_{2}>0,\>m}\sum\limits_{j_{1},j_{2}\geq0,\>j\leq \frac{(2\alpha-1)m_{1}}{2}}$.

\noindent Now we consider (\ref{3.053}). We make the change of variables (\ref{3.042}).

\noindent Thus the Jacobian determinant equals
\begin{eqnarray}
\frac{\partial(u,v,w,\mu_{2})}{\partial(\xi_{1},\xi_{2},\mu_{1},\mu_{2})}
=2\left[\frac{\mu_{1}}{\xi_{1}}-\frac{\mu_{2}}{\xi_{2}}\right].\label{3.064}
\end{eqnarray}
We assume that $D_{j,j_{1},j_{2},m,m_{1},m_{2}}^{(5)}$ is the image of the subset of all points
$$(\xi_{1},\mu_{1},\sigma_{1},\xi_{2},\mu_{2},\sigma_{2})\in D_{j,j_{1},j_{2},m,m_{1},m_{2}}^{(3)},$$
 which satisfies (\ref{3.052}) under the transformation
(\ref{3.042}). Combining (\ref{3.051}) with (\ref{3.064}), we have that
\begin{eqnarray}
\left|\frac{\partial(u,v,w,\xi_{1})}{\partial(\xi_{1},\xi_{2},\mu_{1},\mu_{2})}\right|\sim 2^{m+\frac{(\alpha-2)m_{1}}{2}}.\label{3.065}
\end{eqnarray}
Let
$
f_{m_{k},j_{k}}=\eta_{m_{k}}(\xi_{k})\eta_{j_{k}}(\sigma_{k})f_{k}(\xi_{k},\mu_{k},\tau_{k})(k=1,2)
 $
 and $$H_{2}(u,v,w,\xi_{1},\sigma_{1},\sigma_{2})$$
be $\eta_{m}(\xi)\eta_{j}(\sigma)\prod\limits_{k=1}^{2}f_{m_{k},j_{k}}$
under the change of the variables (\ref{3.042}) and
\begin{eqnarray}
M_{4}=F(u,v,w)H_{2}(u,v,w,\xi_{1},\sigma_{1},\sigma_{2}), dV^{(4)}=dudvdwd\xi_{1}d\sigma_{1}d\sigma_{2}.\label{3.066}
\end{eqnarray}
Thus, (\ref{3.055}) can be controlled by
\begin{eqnarray}
\hspace{-0.5cm}C\sum\limits_{{\rm min}\{ j,j_{1},j_{2},m_{1},m_{2}\}\geq0,m}\int_{D_{j,j_{1},j_{2},m,m_{1},m_{2}}^{(2)}}
2^{jb^{\prime}-(j_{1}+j_{2})b-2m_{1}s_{1}+m}\frac{M_{4}dV^{(4)}}
{\left|\frac{\partial(u,v,w,\xi_{1})}{\partial(\xi_{1},\xi_{2},\mu_{1},\mu_{2})}\right|}\label{3.067},
\end{eqnarray}
We assume that  $h(\xi)$ is defined as in (\ref{3.046}), for fixed $\xi_{2},\mu_{1},\mu_{2},$ from (\ref{3.053}), we have that
\begin{eqnarray}
&&|h^{\prime}(\xi_{1})|=\left|\alpha(\alpha+1)\xi_{1}|\xi_{1}|^{\alpha-2}+2
\left(\frac{\mu_{1}}{\xi_{1}}\right)^{2}\xi_{1}\right|\geq \alpha(\alpha+1)|\xi_{1}|^{\alpha-1}
\geq C2^{(\alpha-1)m_{1}},\nonumber\\&&|h(\xi_{1})|\leq C2^{j+\frac{(\alpha-2)m_{1}}{2}}\label{3.068}.
\end{eqnarray}
Combining  (\ref{3.068}) with Lemma 2.6, for fixed $\xi_{2},\mu_{1},\mu_{2},$ we have that
the measure of $\xi_{1}$ can be controlled by
 $C2^{j-\frac{\alpha m_{1}}{2}}$.
By using the Cauchy-Schwartz inequality with respect to $\xi_{1}$
 and the inverse change of variables and the Cauchy-Schwartz inequality with respect to $u,v,w$
  and the Cauchy-Schwartz inequality with respect to $\sigma_{1}$
   and $\sigma_{2}$, we have that (\ref{3.067}) can be bounded by
\begin{eqnarray}
&&\hspace{-0.9cm}C\sum
\int_{D_{j,j_{1},j_{2},m,m_{1},m_{2}}^{(2)}}
2^{jb^{\prime}-(j_{1}+j_{2})b-2m_{1}s_{1}+m}\frac{M_{4}dV^{(4)}}
{\left|\frac{\partial(u,v,w,\xi_{1})}{\partial(\xi_{1},\xi_{2},\mu_{1},\mu_{2})}\right|}\nonumber\\
&&\hspace{-0.9cm}\leq C\sum2^{2j\epsilon-(j_{1}+j_{2})b-(2s_{1}+\frac{\alpha}{4})m_{1}+m}
\int F(u,v,w)\left(\int\frac{H_{2}^{2}(u,v,w,\xi_{1},\sigma_{1},\sigma_{2})}{
\left|\frac{\partial(u,v,w,\xi_{1})}{\partial(\xi_{1},\xi_{2},\mu_{1},\mu_{2})}\right|^{2}}
d\xi_{1}\right)^{\frac{1}{2}}dV^{(2)}\nonumber\\
&&\hspace{-0.9cm}\leq C\sum2^{2j\epsilon-(j_{1}+j_{2})b- m_{1}(2s_{1}+\frac{\alpha-1}{2})+\frac{m}{2}}
\|F\|_{L^{2}}\int
\left(\int\frac{H_{2}^{2}(u,v,w,\xi_{1},\sigma_{1},\sigma_{2})}{
\left|\frac{\partial(u,v,w,\xi_{1})}{\partial(\xi_{1},\xi_{2},\mu_{1},\mu_{2})}\right|}dV^{(5)}\right)^{\frac{1}{2}}
d\sigma_{1}d\sigma_{2}\nonumber\\
&&\hspace{-0.9cm}\leq C\sum2^{2j\epsilon-(j_{1}+j_{2})b-(2s_{1}+\frac{\alpha-1}{2})m_{1}+\frac{m}{2}}
\|F\|_{L^{2}}\int
\left(\int \prod\limits_{k=1}^{2}f_{m_{k},j_{k}}^{2}d\xi_{1}d\mu_{1}d\xi_{2}d\mu_{2}\right)^{\frac{1}{2}}
d\sigma_{1}d\sigma_{2}\nonumber\\
&&\hspace{-0.9cm}\leq C\sum2^{2j\epsilon-(j_{1}+j_{2})
\epsilon-32\alpha m_{1}\epsilon+\frac{m}{2}}\|F\|_{L^{2}}
\left(\int \prod\limits_{k=1}^{2}f_{m_{k},j_{k}}^{2}dV\right)^{\frac{1}{2}}\nonumber\\&&
\hspace{-0.9cm}\leq C\|F\|_{L^{2}}\left(\prod\limits_{j=1}^{2}\|F_{j}\|_{L^{2}}\right),\label{3.069}
\end{eqnarray}
where $\sum=\sum\limits_{m_{1}, m_{2}>0,\>m}\sum\limits_{j_{1},j_{2}\geq0,\>j\leq
\frac{(2\alpha-1)m_{1}}{2}}$ and $dV^{(5)}=dudvdwd\xi_{1}$.

\noindent (5). Region $\Omega_{5}.$

\noindent
In this region, we consider (\ref{3.016}), (\ref{3.017}), respectively.

\noindent When (\ref{3.016}) is valid,  one of (\ref{3.018})-(\ref{3.020}) must occur.

\noindent When (\ref{3.018}) is valid, since  $s_{1}\geq -\frac{\alpha-1}{4}+4\alpha\epsilon$, we have that
\begin{eqnarray}
&&K_{1}(\xi_{1},\mu_{1},\tau_{1},\xi,\mu,\tau)\leq\frac{|\xi|\langle\xi\rangle^{s_{1}}}
{\langle\sigma\rangle^{-b^{\prime}}
\prod\limits_{j=1}^{2}\langle\xi_{j}\rangle^{s_{1}}\langle\sigma_{j}\rangle^{b}}
\leq C\frac{|\xi|^{\frac{5-\alpha}{4}}|\xi_{1}|^{\frac{\alpha-1}{4}+s_{1}}}
{\langle\sigma\rangle^{-b^{\prime}}
\prod\limits_{j=1}^{2}\langle\xi_{j}\rangle^{s_{1}}\langle\sigma_{j}\rangle^{b}}\nonumber\\
&&\leq C\frac{|\xi|^{-\frac{\alpha-3}{4}+2\epsilon}|\xi_{2}|^{\frac{\alpha-1}{4}-s_{1}-\frac{\alpha}{2}+2\alpha\epsilon}}
{\prod\limits_{j=1}^{2}\langle\sigma_{j}\rangle^{b}}
\leq C\frac{|\xi|^{-\frac{\alpha-3}{4}+2\epsilon}|\xi_{2}|^{-\frac{1}{2}-2\alpha\epsilon}}
{\prod\limits_{j=1}^{2}\langle\sigma_{j}\rangle^{b}}
\leq C\frac{|\xi_{1}|^{-\frac{1}{2}}|\xi_{2}|^{\frac{\alpha}{4}}}
{\prod\limits_{j=1}^{2}\langle\sigma_{j}\rangle^{b}}
.\label{3.070}
\end{eqnarray}
Thus, combining (\ref{2.017}) with (\ref{3.070}),  we have that
\begin{eqnarray*}
|{\rm Int_{5}}|\leq C\|F\|_{L_{\tau\xi\mu}^{2}}\left(\prod\limits_{j=1}^{2}\|F_{j}\|_{L_{\tau\xi\mu}^{2}}\right).
\end{eqnarray*}
When (\ref{3.019}) is valid, since  $s_{1}\geq -\frac{\alpha-1}{4}+4\alpha\epsilon$ and
 $\langle \sigma\rangle^{b^{\prime}}\langle \sigma_{1}\rangle^{b}
 \leq \langle \sigma\rangle^{-b}\langle \sigma_{1}\rangle^{b^{\prime}},$
we have that
\begin{eqnarray}
&&K_{1}(\xi_{1},\mu_{1},\tau_{1},\xi,\mu,\tau)\leq\frac{|\xi|\langle\xi\rangle^{s_{1}}}
{\langle\sigma\rangle^{-b^{\prime}}
\prod\limits_{j=1}^{2}\langle\xi_{j}\rangle^{s_{1}}\langle\sigma_{j}\rangle^{b}}
\leq C\frac{|\xi|^{\frac{5-\alpha}{4}}|\xi_{1}|^{\frac{\alpha-1}{4}+s_{1}}}
{\langle\sigma\rangle^{-b^{\prime}}
\prod\limits_{j=1}^{2}\langle\xi_{j}\rangle^{s_{1}}\langle\sigma_{j}\rangle^{b}}\nonumber\\
&&\leq C\frac{|\xi|^{\frac{\alpha-3}{4}+2\epsilon}|\xi_{2}|^{\frac{\alpha-1}{4}-s_{1}-\frac{\alpha}{2}+2\alpha\epsilon}}
{\langle\sigma\rangle^{b}\langle\sigma_{2}\rangle^{b}}
\leq C\frac{|\xi|^{-\frac{\alpha-3}{4}+2\epsilon}|\xi_{2}|^{-\frac{1}{2}-2\alpha\epsilon}}
{\langle\sigma\rangle^{b}\langle\sigma_{2}\rangle^{b}}\leq C\frac{|\xi|^{-\frac{1}{2}}|\xi_{2}|^{\frac{\alpha}{4}}}
{\langle\sigma\rangle^{b}\langle\sigma_{2}\rangle^{b}}
.\label{3.071}
\end{eqnarray}
Thus, combining (\ref{2.019}) with (\ref{3.071}),  we have that
\begin{eqnarray*}
|{\rm Int_{5}}|\leq C\|F\|_{L_{\tau\xi\mu}^{2}}\left(\prod\limits_{j=1}^{2}\|F_{j}\|_{L_{\tau\xi\mu}^{2}}\right).
\end{eqnarray*}
When (\ref{3.020}) is valid, this case can be proved similarly to (\ref{3.019}) with the aid of (\ref{2.018}).

\noindent  When (\ref{3.017}) is valid, from Lemma 2.8, we have that
\begin{eqnarray}
\left|\frac{\mu_{1}}{\xi_{1}}-\frac{\mu_{2}}{\xi_{2}}\right|\sim |\xi||\xi_{1}|^{\frac{\alpha-2}{2}}.\label{3.072}
\end{eqnarray}
We consider $|\sigma|\geq |\xi_{1}|^{\frac{2\alpha-1}{2}}$ and $|\sigma|<|\xi_{1}|^{\frac{2\alpha-1}{2}}$, respectively.

\noindent
When $|\sigma|\geq |\xi_{1}|^{\frac{2\alpha-1}{2}}$, since $s_{1}\geq-\frac{\alpha-1}{4}+4\alpha\epsilon$, we have that
\begin{eqnarray*}
&&K_{1}(\xi_{1},\mu_{1},\tau_{1},\xi,\mu,\tau)\leq\frac{|\xi|\langle\xi\rangle^{s_{1}}}
{\langle\sigma\rangle^{-b^{\prime}}
\prod\limits_{j=1}^{2}\langle\xi_{j}\rangle^{s_{1}}\langle\sigma_{j}\rangle^{b}}
\leq C\frac{|\xi|^{\frac{5-\alpha}{4}}|\xi_{1}|^{\frac{\alpha-1}{4}+s_{1}}}
{\langle\sigma\rangle^{-b^{\prime}}
\prod\limits_{j=1}^{2}\langle\xi_{j}\rangle^{s_{1}}\langle\sigma_{j}\rangle^{b}}\nonumber\\
&&\leq C\frac{|\xi|^{\frac{5-\alpha}{4}}|\xi_{2}|^{-\frac{\alpha}{4}-s_{1}+(2\alpha-1)\epsilon}}
{\prod\limits_{j=1}^{2}\langle\sigma_{j}\rangle^{b}}
\leq C\frac{|\xi|^{\frac{5-\alpha}{4}}|\xi_{2}|^{-\frac{1}{4}-(14\alpha-1)\epsilon}}
{\prod\limits_{j=1}^{2}\langle\sigma_{j}\rangle^{b}}
\leq C\frac{|\xi_{1}|^{-\frac{1}{2}}|\xi_{2}|^{\frac{\alpha}{4}}}
{\prod\limits_{j=1}^{2}\langle\sigma_{j}\rangle^{b}}.
\end{eqnarray*}
When $|\sigma|<|\xi_{1}|^{\frac{2\alpha-1}{2}}$ is valid,
 we consider
\begin{eqnarray}
&&\left|(\alpha+1)(|\xi_{1}|^{\alpha}-|\xi_{2}|^{\alpha})
-\left[\left(\frac{\mu_{1}}{\xi_{1}}\right)^{2}
-\left(\frac{\mu_{2}}{\xi_{2}}\right)^{2}\right]\right|>2^{j+\frac{m}{2}+\frac{(\alpha-2)m_{1}}{2}},\label{3.073}\\
&&\left|(\alpha+1)(|\xi_{1}|^{\alpha}-|\xi_{2}|^{\alpha})
-\left[\left(\frac{\mu_{1}}{\xi_{1}}\right)^{2}
-\left(\frac{\mu_{2}}{\xi_{2}}\right)^{2}\right]\right|\leq 2^{j+\frac{m}{2}+\frac{(\alpha-2)m_{1}}{2}},\label{3.074}
\end{eqnarray}
respectively.
We dyadically decompose with respect to
\begin{eqnarray*}
\langle\sigma\rangle\sim 2^{j},\langle\sigma_{1}\rangle\sim 2^{j_{1}},
\langle\sigma_{2}\rangle\sim 2^{j_{2}},|\xi|\sim 2^{m},|\xi_{1}|\sim 2^{m_{1}},|\xi_{2}|\sim 2^{m_{2}}.
\end{eqnarray*}
Let $D_{j,j_{1},j_{2},m,m_{1},m_{2}}^{(6)}$ be the image of set of all points
$(\xi_{1},\mu_{1},\tau_{1},\xi,\mu,\tau)\in D^{*}$ satisfying
\begin{eqnarray}
&&|\xi_{1}|\geq 2A,4|\xi|\leq |\xi_{2}|\sim|\xi_{1}|,|\xi|> 2A,\xi_{1}\xi_{2}<0,|\xi|\sim 2^{m},
|\xi_{1}|\sim 2^{m_{1}},|\xi_{2}|\sim 2^{m_{2}}\nonumber\\
&&\langle\sigma\rangle\sim 2^{j_{1}},\langle\sigma_{1}\rangle\sim 2^{j_{1}},
\langle\sigma_{2}\rangle\sim 2^{j_{2}}\label{3.075}
\end{eqnarray}
under the transformation $(\xi_{1},\mu_{1},\tau_{1},\xi_{2},\mu_{2},\tau_{2})\longrightarrow
(\xi_{1},\mu_{1},\sigma_{1},\xi_{2},\mu_{2},\sigma_{2})$.

\noindent Thus, we have that
\begin{eqnarray}
&&\hspace{-2cm}{\rm Int_{5}}\leq C\sum\limits_{m_{1}, m_{2},\>m>0}\sum\limits_{j_{1},j_{2}\geq0,
\>0<j\leq \frac{(2\alpha-1)m_{1}}{2}}
\left(\int_{D_{j,j_{1},j_{2},m,m_{1},m_{2}}^{(6)}}
2^{jb^{\prime}-(j_{1}+j_{2})b-2m_{2}s_{1}+m}PdV
\right),\label{3.076}
\end{eqnarray}
where $P$ and $dV$ are defined as in (\ref{3.029}).

\noindent In this case, we consider (\ref{3.073}), (\ref{3.074}), respectively.

\noindent When (\ref{3.073}) is valid, we make the change of variables as in (\ref{3.033}).

\noindent Thus the Jacobian determinant equals
\begin{eqnarray}
\frac{\partial(u,v,w,\mu_{2})}{\partial(\xi_{1},\xi_{2},\mu_{1},\mu_{2})}=(\alpha+1)(|\xi_{1}|^{\alpha}-|\xi_{2}|^{\alpha})
-\left[\left(\frac{\mu_{1}}{\xi_{1}}\right)^{2}-\left(\frac{\mu_{2}}{\xi_{2}}\right)^{2}\right]\label{3.077}.
\end{eqnarray}
We assume that $D_{j,j_{1},j_{2},m,m_{1},m_{2}}^{(7)}$ is the image of the subset of all points
$$(\xi_{1},\mu_{1},\sigma_{1},\xi_{2},\mu_{2},\sigma_{2})\in D_{j,j_{1},j_{2},m,m_{1},m_{2}}^{(6)},$$
which satisfies (\ref{3.073}) under the transformation
(\ref{3.033}). Combining (\ref{3.073}) with (\ref{3.077}), we have that
\begin{eqnarray}
\left|\frac{\partial(u,v,w,\mu_{2})}{\partial(\xi_{1},\xi_{2},\mu_{1},\mu_{2})}\right|>2^{j+\frac{m}{2}+\frac{(\alpha-2)m_{1}}{2}}\label{3.078}.
\end{eqnarray}
Let $G_{3}(u,v,w,\mu_{2},\sigma_{1},\sigma_{2})$ be $\eta_{m}(\xi)\eta_{j}(\sigma)\prod\limits_{k=1}^{2}f_{m_{k},j_{k}}$
 under the change of the variables (\ref{3.033})  and
\begin{eqnarray}
M_{5}=F(u,v,w)G_{3}(u,v,w,\mu_{2},\sigma_{1},\sigma_{2}),
dV^{(1)}=dudvdwd\mu_{2}d\sigma_{1}d\sigma_{2}\label{3.079}.
\end{eqnarray}
Thus, (\ref{3.076}) can be controlled by
\begin{eqnarray}
&&\hspace{-1cm}C\sum\limits_{m_{1}, m_{2},m>0}\sum\limits_{j_{1},j_{2}\geq0,
\>0<j\leq \frac{(2\alpha-1)m_{1}}{2}}\int_{D_{j,j_{1},j_{2},m,m_{1},m_{2}}^{(7)}}
2^{jb^{\prime}-(j_{1}+j_{2})b-s_{1}m_{1}+m}\frac{M_{5}dV^{(1)}}
{\left|\frac{\partial(u,v,w,\mu_{2})}{\partial(\xi_{1},\xi_{2},\mu_{1},\mu_{2})}\right|}\label{3.080}.
\end{eqnarray}
Combining (\ref{3.060})-(\ref{3.062}) with Lemma 2.6, for fixed $\xi_{1},\xi_{2},\mu_{1},\sigma_{1},\sigma_{2},$
 we have that the Lebesgue measure of $\mu_{2}$ can be controlled by $C2^{j-\frac{\alpha m_{1}}{2}}$.
By using the Cauchy-Schwartz inequality with respect to $\mu_{2}$
 and the inverse change of variables related to  (\ref{3.033}) and the Cauchy-Schwartz inequality with respect to $u,v,w$
  and the Cauchy-Schwartz inequality with respect to $\sigma_{1}$
   and $\sigma_{2}$, we have that (\ref{3.080}) can be controlled by
\begin{eqnarray}
&&\hspace{-0.5cm}C\sum
\int_{D_{j,j_{1},j_{2},m,m_{1},m_{2}}^{(7)}}
2^{jb^{\prime}-(j_{1}+j_{2})b-s_{1}m_{1}+m}\frac{M_{5}dV^{(1)}}
{\left|\frac{\partial(u,v,w,\mu_{2})}{\partial(\xi_{1},\xi_{2},\mu_{1},\mu_{2})}\right|}\nonumber\\
&&\leq C\sum2^{2j\epsilon-(j_{1}+j_{2})b+(-\frac{\alpha}{4}-s_{1})m_{1}+m}
\int F(u,v,w)
\left(\int\frac{G_{3}^{2}(u,v,w,\mu_{2},\sigma_{1},\sigma_{2})}{
\left|\frac{\partial(u,v,w,\mu_{2})}{\partial(\xi_{1},\xi_{2},\mu_{1},\mu_{2})}\right|^{2}}
d\mu_{2}\right)^{\frac{1}{2}}dV^{(2)}\nonumber\\
&&\leq C\sum2^{jb^{\prime}-(j_{1}+j_{2})b+(-s_{1}-\frac{\alpha-1}{2}) m_{1}+\frac{3 m}{4}}\|F\|_{L^{2}}\int
\left(\int\frac{G_{3}^{2}(u,v,w,\mu_{2},\sigma_{1},\sigma_{2})}{
\left|\frac{\partial(u,v,w,\mu_{2})}{\partial(\xi_{1},\xi_{2},\mu_{1},\mu_{2})}\right|}
dV^{(3)}\right)^{\frac{1}{2}}
d\sigma_{1}d\sigma_{2}\nonumber\\
&&\leq C\sum
2^{jb^{\prime}-(j_{1}+j_{2})b+ m_{1}(-s_{1}+\frac{5-2\alpha}{4})}
\|F\|_{L^{2}}\int
\left(\int \prod\limits_{k=1}^{2}f_{m_{k},j_{k}}^{2}d\xi_{1}d\mu_{1}d\xi_{2}d\mu_{2}\right)^{\frac{1}{2}}
d\sigma_{1}d\sigma_{2}\nonumber\\&&
\leq C\sum
2^{jb^{\prime}-(j_{1}+j_{2})\epsilon-16m_{1}\epsilon}
\|F\|_{L^{2}}
\left(\int \prod\limits_{k=1}^{2}f_{m_{k},j_{k}}^{2}dV\right)^{\frac{1}{2}}
\leq C\|F\|_{L^{2}}\left(\prod\limits_{j=1}^{2}\|F_{j}\|_{L^{2}}\right),\label{3.081}
\end{eqnarray}
where $\sum=\sum\limits_{m_{1}, m_{2},m>0}\sum\limits_{j_{1},j_{2}\geq0,\>0<j\leq \frac{(2\alpha-1)m_{1}}{2}}$.

\noindent Now we consider (\ref{3.074}). We make the change of  variables as in (\ref{3.042}).

\noindent Thus the Jacobian determinant equals
\begin{eqnarray}
\frac{\partial(u,v,w,\mu_{2})}{\partial(\xi_{1},\xi_{2},\mu_{1},\mu_{2})}
=2\left[\frac{\mu_{1}}{\xi_{1}}-\frac{\mu_{2}}{\xi_{2}}\right].\label{3.082}
\end{eqnarray}
We assume that $D_{j,j_{1},j_{2},m,m_{1},m_{2}}^{(8)}$ is the image of the subset of all points
$$(\xi_{1},\mu_{1},\sigma_{1},\xi_{2},\mu_{2},\sigma_{2})\in D_{j,j_{1},j_{2},m,m_{1},m_{2}}^{(6)},$$
 which satisfies (\ref{3.074}) under the transformation
(\ref{3.042}). Combining (\ref{3.072}) with (\ref{3.082}), we have that
\begin{eqnarray}
\left|\frac{\partial(u,v,w,\xi_{1})}{\partial(\xi_{1},\xi_{2},\mu_{1},\mu_{2})}\right|\sim 2^{m+\frac{(\alpha-2)m_{1}}{2}}.\label{3.083}
\end{eqnarray}
Let
$
f_{m_{k},j_{k}}=\eta_{m_{k}}(\xi_{k})\eta_{j_{k}}(\sigma_{k})f_{k}(\xi_{k},\mu_{k},\tau_{k})(k=1,2)
$
 and  $$H_{3}(u,v,w,\xi_{1},\sigma_{1},\sigma_{2})$$
be $\eta_{m}(\xi)\eta_{j}(\sigma)\prod\limits_{k=1}^{2}f_{m_{k},j_{k}}$
under the change of the variables (\ref{3.042}) and
\begin{eqnarray}
M_{6}=F(u,v,w)H_{3}(u,v,w,\xi_{1},\sigma_{1},\sigma_{2}), dV^{(4)}=dudvdwd\xi_{1}d\sigma_{1}d\sigma_{2}.\label{3.084}
\end{eqnarray}
Thus, (\ref{3.076}) can be controlled by
\begin{eqnarray}
&&\hspace{-1cm}C\sum\limits_{m_{1}, m_{2},m>0}\sum\limits_{j_{1},j_{2}\geq0,\>0<j\leq \frac{(2\alpha-1)m_{1}}{2}}
\int_{D_{j,j_{1},j_{2},m,m_{1},m_{2}}^{(8)}}
2^{jb^{\prime}-(j_{1}+j_{2})b-s_{1}m_{1}+m}\frac{M_{6}dV^{(4)}}
{\left|\frac{\partial(u,v,w,\xi_{1})}{\partial(\xi_{1},\xi_{2},\mu_{1},\mu_{2})}\right|}\label{3.085}.
\end{eqnarray}
We assume  that $h(\xi)$ is defined as in (\ref{3.046}), for fixed $\xi_{2},\mu_{1},\mu_{2},$
from (\ref{3.074}), we have that
\begin{eqnarray}
&&|h^{\prime}(\xi_{1})|=\left|\alpha(\alpha+1)\xi_{1}|\xi_{1}|^{\alpha-2}+2
\left(\frac{\mu_{1}}{\xi_{1}}\right)^{2}\xi_{1}\right|\geq \alpha(\alpha+1)|\xi_{1}|^{\alpha-1}
\geq C2^{(\alpha-1)m_{1}},\nonumber\\&&|h(\xi_{1})|\leq C2^{j+\frac{(\alpha-2)m_{1}}{2}+\frac{m}{2}}.\label{3.086}
\end{eqnarray}
For fixed $\xi_{2},\mu_{1},\mu_{2},$
combining (\ref{3.086}) with Lemma 2.6,
we have that
the Lebesgue measure of $\xi_{1}$ can be controlled by
 $C2^{j-\frac{\alpha m_{1}}{2}+\frac{m}{2}}$.
By using the Cauchy-Schwartz inequality with respect to $\xi_{1}$
 and the inverse change of variables related to  (\ref{3.042}) and the Cauchy-Schwartz inequality with respect to $u,v,w$
  and the Cauchy-Schwartz inequality with respect to $\sigma_{1}$
   and  $\sigma_{2}$, we have that (\ref{3.085}) can be bounded by
\begin{eqnarray}
&&\hspace{-0.5cm}C\sum\int_{D_{j,j_{1},j_{2},m,m_{1},m_{2}}^{(2)}}
2^{jb^{\prime}-(j_{1}+j_{2})b-s_{1}m_{1}+m}\frac{M_{6}dV^{(4)}}
{\left|\frac{\partial(u,v,w,\xi_{1})}{\partial(\xi_{1},\xi_{2},\mu_{1},\mu_{2})}\right|}\nonumber\\
&&\leq C\sum2^{2j\epsilon-(j_{1}+j_{2})b+(-s_{1}-\frac{\alpha-1}{2})m_{1}+\frac{3m}{4}}
\int F(u,v,w)\left(\int\frac{H_{3}^{2}(u,v,w,\xi_{1},\sigma_{1},\sigma_{2})}{
\left|\frac{\partial(u,v,w,\xi_{1})}{\partial(\xi_{1},\xi_{2},\mu_{1},\mu_{2})}\right|^{2}}d\xi_{1}\right)^{\frac{1}{2}}dV^{(2)}\nonumber\\
&&\leq C\sum2^{2j\epsilon-(j_{1}+j_{2})b
+ m_{1}(-s_{1}+\frac{5-2\alpha}{4})}\|F\|_{L^{2}}\int
\left(\int\frac{H^{2}(u,v,w,\xi_{1},\sigma_{1},\sigma_{2})}{
\left|\frac{\partial(u,v,w,\xi_{1})}{\partial(\xi_{1},\xi_{2},\mu_{1},\mu_{2})}\right|}dV^{(5)}\right)^{\frac{1}{2}}
d\sigma_{1}d\sigma_{2}\nonumber\\
&&\leq C\sum2^{2j\epsilon-(j_{1}+j_{2})b-16\alpha m_{1}\epsilon}
\|F\|_{L^{2}}\int
\left(\int \prod\limits_{k=1}^{2}f_{m_{k},j_{k}}^{2}d\xi_{1}d\mu_{1}d\xi_{2}d\mu_{2}\right)^{\frac{1}{2}}
d\sigma_{1}d\sigma_{2}\nonumber\\
&&\leq C\sum2^{2j\epsilon-(j_{1}+j_{2})
\epsilon-16\alpha m_{1}\epsilon}\|F\|_{L^{2}}
\left(\int \prod\limits_{k=1}^{2}f_{m_{k},j_{k}}^{2}dV
\right)^{\frac{1}{2}}\leq
C\|F\|_{L^{2}}\left(\prod\limits_{j=1}^{2}\|F_{j}\|_{L^{2}}\right),\label{3.087}
\end{eqnarray}
where $\sum=\sum\limits_{m_{1}, m_{2},m>0}\sum\limits_{j_{1},j_{2}\geq0,
\>0<j\leq \frac{(2\alpha-1)m_{1}}{2}}$ and $dV^{(5)}:=dudvdw d\xi_{1}$.

\noindent (6). Region $\Omega_{6}.$

\noindent
In this region, we consider (\ref{3.016}),  (\ref{3.017}),  respectively.

\noindent When (\ref{3.016}) is valid,  one of (\ref{3.018})-(\ref{3.020}) must occur.

\noindent When (\ref{3.018}) is valid, since  $s_{1}\geq -\frac{\alpha-1}{4}+4\alpha\epsilon$, we have that
\begin{eqnarray}
K_{1}(\xi_{1},\mu_{1},\tau_{1},\xi,\mu,\tau)\leq C\frac{|\xi|^{1-s_{1}}}
{\langle\sigma\rangle^{-b^{\prime}}
\prod\limits_{j=1}^{2}\langle\sigma_{j}\rangle^{b}}\leq C\frac{|\xi_{2}|^{-s_{1}+\frac{1-\alpha}{2}+2(\alpha+1)\epsilon}}
{\prod\limits_{j=1}^{2}\langle\sigma_{j}\rangle^{b}}
\leq C\frac{|\xi_{1}|^{-\frac{1}{2}}|\xi_{2}|^{\frac{\alpha}{4}}}
{\prod\limits_{j=1}^{2}\langle\sigma_{j}\rangle^{b}}
.\label{3.088}
\end{eqnarray}
Thus, combining (\ref{2.017}) with (\ref{3.088}),  we have that
\begin{eqnarray*}
|{\rm Int_{6}}|\leq C\|F\|_{L_{\tau\xi\mu}^{2}}\left(\prod\limits_{j=1}^{2}\|F_{j}\|_{L_{\tau\xi\mu}^{2}}\right).
\end{eqnarray*}
When (\ref{3.019}) is valid, since  $s_{1}\geq -\frac{\alpha-1}{4}+4\alpha\epsilon$ and
$\langle \sigma\rangle^{b^{\prime}}\langle \sigma_{1}\rangle^{-b}
\leq \langle \sigma\rangle^{-b}\langle \sigma_{1}\rangle^{b^{\prime}},$
we have that
\begin{eqnarray}
K_{1}(\xi_{1},\mu_{1},\tau_{1},\xi,\mu,\tau)\leq C\frac{|\xi|^{1-s_{1}}}
{\langle\sigma\rangle^{-b^{\prime}}
\prod\limits_{j=1}^{2}\langle\sigma_{j}\rangle^{b}}\leq C\frac{|\xi_{2}|^{-s_{1}+\frac{1-\alpha}{2}+2(\alpha+1)\epsilon}}
{\langle\sigma\rangle^{b}\langle\sigma_{2}\rangle^{b}}\leq
 C\frac{|\xi|^{-\frac{1}{2}}|\xi_{2}|^{\frac{\alpha}{4}}}
{\langle\sigma\rangle^{b}\langle\sigma_{2}\rangle^{b}}
.\label{3.089}
\end{eqnarray}
Thus, combining (\ref{2.020}) with (\ref{3.089}),  we have that
\begin{eqnarray*}
{\rm |Int_{6}|}\leq C\|F\|_{L_{\tau\xi\mu}^{2}}\left(\prod\limits_{j=1}^{2}\|F_{j}\|_{L_{\tau\xi\mu}^{2}}\right).
\end{eqnarray*}
When (\ref{3.020}) is valid, this case can be proved similarly to (\ref{3.019}).

\noindent  When (\ref{3.017}) is valid, from Lemma 2.8,
we have that
\begin{eqnarray}
\left|\frac{\mu_{1}}{\xi_{1}}-\frac{\mu_{2}}{\xi_{2}}\right|\sim |\xi_{1}|^{\frac{\alpha}{2}}.\label{3.090}
\end{eqnarray}
We consider $|\sigma|\geq |\xi_{1}|^{\frac{2\alpha-1}{2}}$ and $|\sigma|<|\xi_{1}|^{\frac{2\alpha-1}{2}}$, respectively.

\noindent
When $|\sigma|\geq |\xi_{1}|^{\frac{2\alpha-1}{2}}$, since $s_{1}\geq-\frac{\alpha-1}{4}+4\alpha\epsilon$, we have that
\begin{eqnarray}
K_{1}(\xi_{1},\mu_{1},\tau_{1},\xi,\mu,\tau)\leq C\frac{|\xi|^{1-s_{1}}}
{\langle\sigma\rangle^{-b^{\prime}}
\prod\limits_{j=1}^{2}\langle\sigma_{j}\rangle^{b}}\leq C\frac{|\xi_{1}|^{-\frac{\alpha-1}{4}-(14\alpha-2)\epsilon}}
{\prod\limits_{j=1}^{2}\langle\sigma_{j}\rangle^{b}}
\leq C\frac{|\xi_{1}|^{-\frac{1}{2}}|\xi_{2}|^{\frac{\alpha}{4}}}
{\prod\limits_{j=1}^{2}\langle\sigma_{j}\rangle^{b}}
.\label{3.091}
\end{eqnarray}
Thus, combining (\ref{2.017}) with (\ref{3.091}),  we have that
\begin{eqnarray*}
|{\rm Int_{6}}|\leq C\|F\|_{L_{\tau\xi\mu}^{2}}\left(\prod\limits_{j=1}^{2}\|F_{j}\|_{L_{\tau\xi\mu}^{2}}\right).
\end{eqnarray*}
When $|\sigma|<|\xi_{1}|^{\frac{2\alpha-1}{2}}$ is valid, we dyadically decompose with respect to
\begin{eqnarray*}
\langle\sigma\rangle\sim 2^{j},\langle\sigma_{1}\rangle\sim 2^{j_{1}},
\langle\sigma_{2}\rangle\sim 2^{j_{2}},|\xi|\sim 2^{m},|\xi_{1}|\sim 2^{m_{1}},|\xi_{2}|\sim 2^{m_{2}}.
\end{eqnarray*}
Let $D_{j,j_{1},j_{2},m,m_{1},m_{2}}^{(8)}$ be the image of set of all points
$(\xi_{1},\mu_{1},\tau_{1},\xi,\mu,\tau)\in D^{*}$ satisfying
\begin{eqnarray}
&&|\xi_{1}|\geq 2A, |\xi_{1}|\sim |\xi_{2}|,\xi_{1}\xi_{2}<0,
|\xi|\geq \frac{|\xi_{2}|}{4},|\xi|\sim 2^{m},|\xi_{1}|\sim 2^{m_{1}},|\xi_{2}|\sim 2^{m_{2}},
\nonumber\\
&&\langle\sigma\rangle\sim 2^{j}\leq C2^{3m_{1}},
\langle\sigma_{1}\rangle\sim 2^{j_{1}},\langle\sigma_{2}\rangle\sim 2^{j_{2}}\label{3.092}
\end{eqnarray}
under the transformation $(\xi_{1},\mu_{1},\tau_{1},\xi_{2},\mu_{2},\tau_{2})\longrightarrow
(\xi_{1},\mu_{1},\sigma_{1},\xi_{2},\mu_{2},\sigma_{2})$.

\noindent Thus, we have that
\begin{eqnarray}
{\rm Int_{6}}\leq C\sum\limits_{m_{1}, m_{2}, m>0}
\sum\limits_{j_{1},j_{2}\geq0,\>j\leq \frac{(2\alpha-1)m_{1}}{2}}
\int_{D_{j,j_{1},j_{2},m,m_{1},m_{2}}^{(8)}}
2^{jb^{\prime}-(j_{1}+j_{2})b-m_{2}s_{1}+m}PdV,\label{3.093}
\end{eqnarray}
where $P$ and $dV$ are defined in (\ref{3.093}).

\noindent In this case, we consider (\ref{3.031}), (\ref{3.032}), respectively.

\noindent When (\ref{3.031}) is valid, we make the change of variables as in (\ref{3.033}).

\noindent Thus the Jacobian determinant equals
\begin{eqnarray}
\frac{\partial(u,v,w,\mu_{2})}{\partial(\xi_{1},\xi_{2},\mu_{1},\mu_{2})}=(\alpha+1)(|\xi_{1}|^{\alpha}-|\xi_{2}|^{\alpha})
-\left[\left(\frac{\mu_{1}}{\xi_{1}}\right)^{2}-\left(\frac{\mu_{2}}{\xi_{2}}\right)^{2}\right]\label{3.094}.
\end{eqnarray}
We assume that $D_{j,j_{1},j_{2},m,m_{1},m_{2}}^{(9)}$ is the image of the subset of all points
$$(\xi_{1},\mu_{1},\sigma_{1},\xi_{2},\mu_{2},\sigma_{2})\in D_{j,j_{1},j_{2},m,m_{1},m_{2}}^{(8)},$$
which satisfies (\ref{3.031}) under the transformation
(\ref{3.033}). Combining (\ref{3.094}) with (\ref{3.031}), we have that
\begin{eqnarray}
\left|\frac{\partial(u,v,w,\mu_{2})}{\partial(\xi_{1},\xi_{2},\mu_{1},\mu_{2})}
\right|>2^{j+\frac{(\alpha-1)m_{1}}{2}}\label{3.095}.
\end{eqnarray}
Let $G_{4}(u,v,w,\mu_{2},\sigma_{1},\sigma_{2})$ be $\eta_{m}(\xi)\eta_{j}(\sigma)
\prod\limits_{k=1}^{2}f_{m_{k},j_{k}}$ under the change of the variables (\ref{3.033})  and
\begin{eqnarray}
M_{7}=F(u,v,w)G_{4}(u,v,w,\mu_{2},\sigma_{1},\sigma_{2}),
dV^{(1)}=dudvdwd\mu_{2}d\sigma_{1}d\sigma_{2}\label{3.096}.
\end{eqnarray}
Thus, (\ref{3.093}) can be controlled by
\begin{eqnarray}
\hspace{-0.5cm}C\sum\limits_{m, m_{1}, m_{2}>0}\sum\limits_{j_{1},j_{2}\geq0,
\>j\leq \frac{(2\alpha-1)m_{1}}{2},}\int_{D_{j,j_{1},j_{2},m,m_{1},m_{2}}^{(9)}}
2^{jb^{\prime}-(j_{1}+j_{2})b+(1-s_{1})m_{1}}\frac{M_{7}dV^{(1)}}
{\left|\frac{\partial(u,v,w,\mu_{2})}{\partial(\xi_{1},\xi_{2},\mu_{1},\mu_{2})}\right|}\label{3.097}.
\end{eqnarray}
Inspired by \cite{CKS-GAFA, LX},
we define
\begin{eqnarray}
&&f(\mu):=\sigma_{1}+\sigma_{2}-
(\xi|\xi|^{\alpha}-\xi_{1}|\xi_{1}|^{\alpha}-\xi_{2}|\xi_{2}|^{\alpha})+\frac{\xi_{1}\xi_{2}}{\xi}
\left[\frac{\mu_{1}}{\xi_{1}}-\frac{\mu}{\xi_{2}}\right]^{2}\label{3.098}.
\end{eqnarray}
From (\ref{3.098}) and (\ref{3.090}), we have that
\begin{eqnarray}
&&\hspace{-0.5cm}|f(\mu_{2})|=|\sigma_{1}+\sigma_{2}+\phi(\xi_{1},\mu_{1})+\phi(\xi_{2},\mu_{2})
-\phi(\xi,\mu)|=|\tau-\phi(\xi,\mu)|\sim 2^{j},\label{3.099}\\&&|f^{\prime}(\mu_{2})|\sim \left|\frac{\xi_{1}}{\xi}\right|
\left|\frac{\mu_{1}}{\xi_{1}}-\frac{\mu_{2}}{\xi_{2}}\right|\sim 2^{\frac{\alpha m_{1}}{2}} .\label{3.0100}
\end{eqnarray}
Combining (\ref{3.099}), (\ref{3.0100}) with Lemma 2.6, for fixed $\xi_{1},\xi_{2},\mu_{1},\sigma_{1},\sigma_{2},$
the Lebesgue measure of $\mu_{2}$ can be controlled by $C2^{j-\frac{\alpha m_{1}}{2}}$.
By using the Cauchy-Schwartz inequality with respect to $\mu_{2}$
 and the inverse change of variables related to (\ref{3.033}) and the Cauchy-Schwartz inequality with respect to $u,v,w$
  and the Cauchy-Schwartz inequality with respect to $\sigma_{1}$
   and $\sigma_{2}$, we have that (\ref{3.098}) can be bounded by
\begin{eqnarray}
&&\hspace{-0.5cm}C\sum\int_{D_{j,j_{1},j_{2},m,m_{1},m_{2}}^{(9)}}
2^{jb^{\prime}-(j_{1}+j_{2})b+(1-s_{1})m_{1}}\frac{M_{7}dV^{(1)}}
{\left|\frac{\partial(u,v,w,\mu_{2})}{\partial(\xi_{1},\xi_{2},\mu_{1},\mu_{2})}\right|}\nonumber\\
&&\leq C\sum2^{2j\epsilon-(j_{1}+j_{2})b+m_{1}(1-s_{1}-\frac{\alpha}{4})}\int F(u,v,w)
\left(\int\frac{G_{4}^{2}(u,v,w,\mu_{2},\sigma_{1},\sigma_{2})}{
\left|\frac{\partial(u,v,w,\mu_{2})}{\partial(\xi_{1},\xi_{2},\mu_{1},\mu_{2})}\right|^{2}}
d\mu_{2}\right)^{\frac{1}{2}}dV^{(2)}\nonumber\\
&&\leq C\sum2^{jb^{\prime}-(j_{1}+j_{2})b+m_{1}(-s_{1}+\frac{5-2\alpha}{4})}\|F\|_{L^{2}}\int
\left(\int\frac{G_{4}^{2}(u,v,w,\mu_{2},\sigma_{1},\sigma_{2})}{
\left|\frac{\partial(u,v,w,\mu_{2})}{\partial(\xi_{1},\xi_{2},\mu_{1},\mu_{2})}\right|}
dV^{(3)}\right)^{\frac{1}{2}}
d\sigma_{1}d\sigma_{2}\nonumber\\
&&\leq C\sum2^{jb^{\prime}-(j_{1}+j_{2})b-16\alpha m_{1}\epsilon}
\|F\|_{L^{2}}\int
\left(\int \prod\limits_{k=1}^{2}f_{m_{k},j_{k}}^{2}d\xi_{1}d\mu_{1}d\xi_{2}d\mu_{2}\right)^{\frac{1}{2}}
d\sigma_{1}d\sigma_{2}\nonumber\\&&
\leq C\sum
2^{jb^{\prime}-(j_{1}+j_{2})b-16 \alpha m_{1}\epsilon }
\|F\|_{L^{2}}
\left(\int \prod\limits_{k=1}^{2}f_{m_{k},j_{k}}^{2}dV\right)^{\frac{1}{2}}\nonumber\\&&
\leq C\|F\|_{L^{2}}\left(\prod\limits_{j=1}^{2}\|F_{j}\|_{L^{2}}\right),\label{3.0101}
\end{eqnarray}
where $\sum=\sum\limits_{m, m_{1}, m_{2}>0}\sum\limits_{j_{1},j_{2}\geq0,\>0<j\leq \frac{(2\alpha-1)m_{1}}{2}}.$

\noindent Now we consider (\ref{3.032}). We make the change of variables  as in (\ref{3.042}).

\noindent Thus the Jacobian determinant equals
\begin{eqnarray}
\frac{\partial(u,v,w,\mu_{2})}{\partial(\xi_{1},\xi_{2},\mu_{1},\mu_{2})}
=2\left[\frac{\mu_{1}}{\xi_{1}}-\frac{\mu_{2}}{\xi_{2}}\right].\label{3.0102}
\end{eqnarray}
We assume that $D_{j,j_{1},j_{2},m,m_{1},m_{2}}^{(10)}$ is the image of the subset of all points
$$(\xi_{1},\mu_{1},\sigma_{1},\xi_{2},\mu_{2},\sigma_{2})\in D_{j,j_{1},j_{2},m,m_{1},m_{2}}^{(8)},$$
 which satisfies (\ref{3.032}) under the transformation
(\ref{3.042}). Combining (\ref{3.032}) with (\ref{3.0102}), we have that
\begin{eqnarray}
\left|\frac{\partial(u,v,w,\xi_{1})}{\partial(\xi_{1},\xi_{2},\mu_{1},\mu_{2})}\right|\sim 2^{\frac{\alpha m_{1}}{2}}.\label{3.0103}
\end{eqnarray}
Let
$$
f_{m_{k},j_{k}}=\eta_{m_{k}}(\xi_{k})\eta_{j_{k}}(\sigma_{k})f_{k}(\xi_{k},\mu_{k},\tau_{k})(k=1,2)
$$
and
$H_{4}(u,v,w,\xi_{1},\sigma_{1},\sigma_{2})$
be $\eta_{m}(\xi)\eta_{j}(\sigma)\prod\limits_{k=1}^{2}f_{m_{k},j_{k}}$
under the change of the variables (\ref{3.042}) and
\begin{eqnarray}
M_{8}=F(u,v,w)H_{4}(u,v,w,\xi_{1},\sigma_{1},\sigma_{2}), dV^{(4)}=dudvdwd\xi_{1}d\sigma_{1}d\sigma_{2}.\label{3.0104}
\end{eqnarray}
Thus, we have that
\begin{eqnarray}
&&\hspace{-2cm}{\rm Int_{6}}\leq C\sum\limits_{m, m_{1}, m_{2}>0}
\sum\limits_{j_{1},j_{2}\geq0,\>0<j\leq \frac{7m_{1}}{2}}\int_{D_{j,j_{1},j_{2},m,m_{1},m_{2}}^{(10)}}
2^{jb^{\prime}-(j_{1}+j_{2})b+(1-s_{1})m_{1}}\frac{M_{8}dV^{(4)}}
{\left|\frac{\partial(u,v,w,\xi_{1})}{\partial(\xi_{1},\xi_{2},\mu_{1},\mu_{2})}\right|}\label{3.0105}.
\end{eqnarray}
We assume that $h(\xi)$ is defined as in (\ref{3.046}), from (\ref{3.032}), for fixed $\xi_{2},\mu_{1},\mu_{2},$ we have that
\begin{eqnarray}
&&|h^{\prime}(\xi_{1})|=\left|\alpha(\alpha+1)\xi_{1}|\xi_{1}|^{\alpha-2}+2
\left(\frac{\mu_{1}}{\xi_{1}}\right)^{2}\xi_{1}\right|\geq \alpha(\alpha+1)|\xi_{1}|^{\alpha-1}
\geq C2^{(\alpha-1)m_{1}},\nonumber\\&&|h(\xi_{1})|\leq C2^{j+\frac{(\alpha-1)m_{1}}{2}}\label{3.0106}.
\end{eqnarray}
 Combining (\ref{3.0106}) with Lemma 2.6, for fixed $\xi_{2},\mu_{1},\mu_{2},$ we have that
the Lebesgue measure of $\xi_{1}$ can be controlled by
 $C2^{j-\frac{(\alpha-1)m_{1}}{2}}$.
By using the Cauchy-Schwartz inequality with respect to $\xi_{1}$
 and the inverse change of variables related to  (\ref{3.042}) and the Cauchy-Schwartz inequality with respect to $u,v,w$
  and the Cauchy-Schwartz inequality with respect to $\sigma_{1}$
   and  $\sigma_{2}$, we have that (\ref{3.0105}) can be bounded by
\begin{eqnarray*}
&&\hspace{-0.5cm}C\sum\int_{D_{j,j_{1},j_{2},m,m_{1},m_{2}}^{(10)}}
2^{jb^{\prime}-(j_{1}+j_{2})b+(1-s_{1})m_{1}}\frac{M_{8}dV^{(4)}}
{\left|\frac{\partial(u,v,w,\xi_{1})}{\partial(\xi_{1},\xi_{2},\mu_{1},\mu_{2})}\right|}\nonumber\\
&&\leq C\sum2^{2j\epsilon-(j_{1}+j_{2})b+(\frac{5-\alpha}{4}-s_{1})m_{1}}
\int F(u,v,w)
\left(\int\frac{H_{4}^{2}(u,v,w,\xi_{1},\sigma_{1},\sigma_{2})}{
\left|\frac{\partial(u,v,w,\xi_{1})}{\partial(\xi_{1},\xi_{2},\mu_{1},\mu_{2})}\right|^{2}}
d\xi_{1}\right)^{\frac{1}{2}}dV^{(2)}\nonumber\\
&&\leq C\sum2^{2j\epsilon-(j_{1}+j_{2})b+
 (-s_{1}+\frac{5-2\alpha}{4})m_{1}}\|F\|_{L^{2}}\int
\left(\int\frac{H_{4}^{2}(u,v,w,\xi_{1},\sigma_{1},\sigma_{2})}{
\left|\frac{\partial(u,v,w,\xi_{1})}{\partial(\xi_{1},\xi_{2},\mu_{1},\mu_{2})}\right|}dV^{(5)}\right)^{\frac{1}{2}}
d\sigma_{1}d\sigma_{2}\nonumber\\
&&\leq C\sum2^{2j\epsilon-(j_{1}+j_{2})b-4\alpha m_{1}\epsilon}
\|F\|_{L^{2}}\int
\left(\int \prod\limits_{k=1}^{2}f_{m_{k},j_{k}}^{2}d\xi_{1}d\mu_{1}d\xi_{2}d\mu_{2}\right)^{\frac{1}{2}}
d\sigma_{1}d\sigma_{2}\nonumber\\
&&\leq  C\sum2^{2j\epsilon-(j_{1}+j_{2})
\epsilon-4\alpha m_{1}\epsilon}\|F\|_{L^{2}}
\left(\int \prod\limits_{k=1}^{2}f_{m_{k},j_{k}}^{2}dV\right)^{\frac{1}{2}}
\leq C\|F\|_{L^{2}}\left(\prod\limits_{j=1}^{2}\|F_{j}\|_{L^{2}}\right).
\end{eqnarray*}
where $\sum=\sum\limits_{m, m_{1}, m_{2}>0}\sum\limits_{j_{1},j_{2}\geq0,\>0<j\leq
\frac{(2\alpha-1)m_{1}}{2}}$ and $dV^{(5)}=dudvdwd\xi_{1}$.

\noindent (7). Region $\Omega_{7}.$ This case can be proved similarly to Region $\Omega_{6}.$

This ends the proof of Lemma 3.1.

 \begin{Lemma}\label{Lemma3.2}
Let $-\frac{2\alpha-5}{8}+2\alpha \epsilon\leq s<0$ and $b=\frac{1}{2}+\epsilon$ and $b^{\prime}=-\frac{1}{2}+2\epsilon$.
Then, we have that
\begin{eqnarray}
\hspace{-0.6cm}\|\partial_{x}\left[I_{N}(u_{1}u_{2})-I_{N}u_{1}I_{N}u_{2}\right]
\|_{X_{b^{\prime}}^{0,0}}\leq CN^{3\alpha \epsilon}{\rm max}\left\{N^{-\frac{\alpha}{4}}, N^{-\frac{2\alpha-5}{4}}\right\}
\prod_{j=1}^{2}\|I_{N}u_{j}\|_{X_{b}^{0,0}}.\label{3.0107}
\end{eqnarray}
\end{Lemma}
\noindent{\bf Proof.}
To prove (\ref{3.0107}),  by duality, it suffices to  prove that
\begin{eqnarray}
&&\left|\int_{\SR^{3}}\bar{h}\partial_{x}\left[I_{N}(u_{1}u_{2})-I_{N}u_{1}I_{N}u_{2}\right]dxdydt\right|\nonumber\\&&\leq
CN^{3\alpha \epsilon}{\rm max}\left\{N^{-\frac{\alpha}{4}}, N^{-\frac{2\alpha-5}{4}}\right\}
\|h\|_{X_{-b^{\prime}}^{0,0}}\left(\prod_{j=1}^{2}
\|I_{N}u_{j}\|_{X_{b}^{0,0}}\right).\label{3.0108}
\end{eqnarray}
for $h\in X_{-b^{\prime}}^{0,0}.$
Let
\begin{eqnarray}
&&F(\xi,\mu,\tau)=
\langle \sigma\rangle^{-b^{\prime}}M(\xi)\mathscr{F}h(\xi,\mu,\tau),\nonumber\\&&
F_{j}(\xi_{j},\mu_{j},\tau_{j})=M(\xi_{j})
\langle \sigma_{j}\rangle^{b}
\mathscr{F}u_{j}(\xi_{j},\mu,\tau_{j})(j=1,2).\label{3.0109}
\end{eqnarray}
To obtain (\ref{3.0108}), from (\ref{3.0109}), it suffices to prove that
\begin{eqnarray}
&&\int_{D}\frac{|\xi|G(\xi_{1},\xi_{2})
F(\xi,\mu,\tau)\prod\limits_{j=1}^{2}F_{j}(\xi_{j},\mu_{j},\tau_{j})}{\langle\sigma_{j}\rangle^{-b^{\prime}}
\prod\limits_{j=1}^{2}\langle\sigma_{j}\rangle^{b}}
d\xi_{1}d\mu_{1}d\tau_{1}d\xi d\mu d\tau\nonumber\\&&\leq CN^{2\alpha \epsilon}{\rm max}\left\{N^{-\frac{\alpha}{4}}, N^{-\frac{2\alpha-5}{4}}\right\}
\|F\|_{L_{\xi\mu\tau}^{2}}\left(\prod_{j=1}^{2}\|F_{j}\|_{L_{\xi\mu\tau}^{2}}\right),\label{3.0110}
\end{eqnarray}
where
$
G(\xi_{1},\xi_{2})=\frac{M(\xi_{1})M(\xi_{2})-M(\xi)}{M(\xi_{1})M(\xi_{2})}
$
and
$D$ is defined as in Lemma 3.1.

\noindent Without loss of generality,  we assume that
 $F(\xi,\mu,\tau)\geq 0,F_j(\xi_{j},\mu_{j},\tau_{j})\geq 0(j=1,2)$.
By symmetry, we can assume that $|\xi_{1}|\geq |\xi_{2}|.$

 \noindent
We define
\begin{eqnarray*}
&&\hspace{-0.8cm}A_1=\left\{(\xi_1,\mu_{1},\tau_1,\xi,\mu,\tau)\in D^{*},
  |\xi_{2}|\leq |\xi_{1}|\leq \frac{N}{2}\right\},\\
&&\hspace{-0.8cm} A_2=\{ (\xi_1,\mu_{1},\tau_1,\xi,\mu,\tau)\in D^{*},
|\xi_1|>\frac{N}{2},|\xi_{1}|\geq |\xi_{2}|, |\xi_{2}|\leq 2\},\\
&&\hspace{-0.8cm} A_3=\{ (\xi_1,\mu_{1},\tau_1,\xi,\mu,\tau)\in D^{*},
|\xi_1|>\frac{N}{2},|\xi_{1}|\geq |\xi_{2}|, 2<|\xi_{2}|\leq N\},\\
&&\hspace{-0.8cm}A_4=\{(\xi_1,\mu_{1},\tau_1,\xi,\mu,\tau)\in D^{*},
|\xi_{1}|>\frac{N}{2},|\xi_{1}|\geq |\xi_{2}|,|\xi_{2}|>N\}.
\end{eqnarray*}
Here $D^{*}$ is defined as in Lemma 3.1.
Obviously, $D^{*}\subset\bigcup\limits_{j=1}^{4}A_{j}.$
We define
\begin{equation}
    K_{2}(\xi_{1},\mu_{1},\tau_{1},\xi,\mu,\tau):=\frac{|\xi|G(\xi_{1},\xi_{2})}
    {\langle\sigma_{j}\rangle^{-b^{\prime}}
\prod\limits_{j=1}^{2}\langle\sigma_{j}\rangle^{b}}\label{3.0111}
\end{equation}
and
\begin{eqnarray*}
J_{k}:=\int_{A_{j}} K_{2}(\xi_{1},\mu_{1},\tau_{1},\xi,\mu,\tau)F(\xi,\mu,\tau)
\prod_{j=1}^{2}F_{j}(\xi_{j},\mu_{j},\tau_{j})
d\xi_{1}d\mu_{1}d\tau_{1}d\xi d\mu d\tau
\end{eqnarray*}
with $1\leq k\leq 4, k\in N.$

\noindent We consider (\ref{3.016}) and (\ref{3.017}), respectively.

\noindent When (\ref{3.016}) is valid,
one of (\ref{3.018})-(\ref{3.020}) must occur,
from Lemma 3.2 of \cite{YLZ}, we have that
\begin{eqnarray*}
\sum\limits_{k=1}^{4}{\rm J}_{k}\leq CN^{-\frac{3\alpha}{4}+1+(2\alpha+2)\epsilon}
\|F\|_{L_{\tau\xi\mu}^{2}}\left(\prod_{j=1}^{2}\|F_{j}\|_{L_{\tau\xi\mu}^{2}}\right).
\end{eqnarray*}
Thus, we only consider the case  (\ref{3.017}).

\noindent (1) Region $A_{1}.$ In this case, since $M(\xi_{1},\xi_{2})=0$, thus we have that ${\rm J}_{1}=0$.

\noindent (2) Region $A_{2}$. From page 902 of \cite{ILM-CPAA}, we have that
\begin{eqnarray}
G(\xi_{1},\xi_{2})\leq C\frac{|\xi_{2}|}{|\xi_{1}|}\label{3.0112}.
\end{eqnarray}
Inserting (\ref{3.0112}) into (\ref{3.0111}) yields
\begin{eqnarray}
    K_{2}(\xi_{1},\mu_{1},\tau_{1},\xi,\mu,\tau)\leq C\frac{|\xi|G(\xi_{1},\xi_{2})}
    {\langle\sigma\rangle^{-b^{\prime}}
\prod\limits_{j=1}^{2}\langle\sigma_{j}\rangle^{b}}\leq
\frac{C}{\langle\sigma\rangle^{-b^{\prime}}
\prod\limits_{j=1}^{2}\langle\sigma_{j}\rangle^{b}}\label{3.0113}.
\end{eqnarray}
 By using the Cauchy-Schwartz inequality with respect to $\xi_{1},\mu_{1},\tau_{1}$,
 from (\ref{3.0113}) and (\ref{3.015}),  we have that
\begin{eqnarray}
&&{\rm J}_{2}\leq C\int_{\SR^{3}}\frac{|\xi||\xi|^{-1}}{\langle \sigma \rangle^{-b^{\prime}}}
\left(\int_{\SR^{3}}\frac{d\xi_{1}d\mu_{1}d\tau_{1}}{\prod\limits_{j=1}^{2}
\langle\sigma_{j}\rangle^{2b}}\right)^{\frac{1}{2}}\nonumber\\&&\qquad\qquad \times\left(\int_{\SR^{3}}
\prod\limits_{j=1}^{2}\left|F_{j}(\xi_{j},\mu_{j},\tau_{j})\right|^{2}
d\xi_{1}d\mu_{1}d\tau_{1}\right)^{\frac{1}{2}}F(\xi,\mu,\tau)d\xi d\mu d\tau\nonumber\\
&&\leq  CN^{-\frac{\alpha}{4}}\|F\|_{L_{\tau\xi\mu}^{2}}\left(\prod\limits_{j=1}^{2}
\|F_{j}\|_{L_{\tau\xi\mu}^{2}}\right).\label{3.0114}
\end{eqnarray}
(3) Region $A_{3}$.
From page 902 of \cite{ILM-CPAA}, we have that (\ref{3.0112}) is valid.
Combining (\ref{3.0111}) with (\ref{3.0112}), we have that
\begin{equation}
    K_{2}(\xi_{1},\mu_{1},\tau_{1},\xi,\mu,\tau)\leq C\frac{{\rm min}
    \left\{|\xi|,|\xi_{1}|,|\xi_{2}|\right\}}{\langle\sigma\rangle^{-b^{\prime}}
\prod\limits_{j=1}^{2}\langle\sigma_{j}\rangle^{b}}\label{3.0115}
\end{equation}
When (\ref{3.017}) is valid, from Lemma 2.8, we have that
\begin{eqnarray}
\left|\frac{\mu_{1}}{\xi_{1}}-\frac{\mu_{2}}{\xi_{2}}\right|\sim |\xi||\xi_{1}|^{\frac{\alpha}{2}-1}.\label{3.0116}
\end{eqnarray}
We consider $|\sigma|\geq |\xi_{1}|^{\frac{2\alpha-1}{2}}$ and $|\sigma|<|\xi_{1}|^{\frac{2\alpha-1}{2}}$, respectively.

\noindent
When $|\sigma|\geq |\xi_{1}|^{\frac{2\alpha-1}{2}}$, since $s_{1}\geq-\frac{\alpha-1}{4}+4\alpha\epsilon$, from (\ref{3.0115}),  we have that
\begin{eqnarray}
&&K_{2}(\xi_{1},\mu_{1},\tau_{1},\xi,\mu,\tau)\leq C
\frac{|\xi_{1}|^{-\frac{2\alpha-1}{4}+(2\alpha-1)\epsilon}{\rm min}\left\{|\xi|,|\xi_{1}|,|\xi_{2}|\right\}}
{\prod\limits_{j=1}^{2}\langle\sigma_{j}\rangle^{b}}\nonumber\\&&
\leq CN^{-\frac{2\alpha-3}{4}+(2\alpha-1)\epsilon}\frac{|\xi_{1}|^{-\frac{1}{2}}|\xi_{2}|}
{\prod\limits_{j=1}^{2}\langle\sigma_{j}\rangle^{b}}
.\label{3.0117}
\end{eqnarray}
Thus, combining (\ref{2.017}) with (\ref{3.0117}),  we have that
\begin{eqnarray*}
|{\rm J}_{3}|\leq CN^{-\frac{2\alpha-3}{4}+(2\alpha-1)\epsilon}\|F\|_{L_{\tau\xi\mu}^{2}}
\left(\prod\limits_{j=1}^{2}\|F_{j}\|_{L_{\tau\xi\mu}^{2}}\right).
\end{eqnarray*}
Now we consider case $|\sigma|<|\xi_{1}|^{\frac{2\alpha-1}{2}}$. We dyadically decompose with respect to
\begin{eqnarray*}
\langle\sigma\rangle\sim 2^{j},\langle\sigma_{1}\rangle\sim 2^{j_{1}},
\langle\sigma_{2}\rangle\sim 2^{j_{2}},|\xi|\sim 2^{m},|\xi_{1}|\sim 2^{m_{1}},|\xi_{2}|\sim 2^{m_{2}}.
\end{eqnarray*}
Let  $D_{j,j_{1},j_{2},m,m_{1},m_{2}}^{(11)}$ be the image of set of all points
$(\xi_{1},\mu_{1},\tau_{1},\xi,\mu,\tau)\in D^{*}$ satisfy
\begin{eqnarray}
&&|\xi_1|>\frac{N}{2},|\xi_{1}|\geq |\xi_{2}|, 2A<|\xi_{2}|\leq N,
|\xi|\sim 2^{m},|\xi_{1}|\sim 2^{m_{1}},|\xi_{2}|\sim 2^{m_{2}},\nonumber\\
&&\langle\sigma\rangle\sim 2^{j}\leq 2^{\frac{2\alpha-1}{2}m},\langle\sigma_{1}\rangle\sim 2^{j_{1}},
\langle\sigma_{2}\rangle\sim 2^{j_{2}}\label{3.0118}
\end{eqnarray}
under the transformation $(\xi_{1},\mu_{1},\tau_{1},\xi_{2},\mu_{2},\tau_{2})\longrightarrow
(\xi_{1},\mu_{1},\sigma_{1},\xi_{2},\mu_{2},\sigma_{2})$. We define
\begin{eqnarray}
&&g_{m_{k},j_k}:=\eta_{m_{k}}(\xi_{k})\eta_{j_{k}}(\sigma_{k})F_{k}(\xi_{k},\mu_{k},\tau_{k})(k=1,2),\label{3.0119}\\
&&g_{m,j}:=\eta_{m}(\xi)\eta_{j}(\sigma)\left|F(\xi,\mu,\sigma_{1}+\sigma_{2}-\phi(\xi_{1},\mu_{1})-\phi(\xi_{2},\mu_{2}))\right|,\label{3.0120}\\
&&Q:=g_{m,j}\prod\limits_{k=1}^{2}g_{m_{k},j_k},dV=d\xi_{1}d\mu_{1}d\sigma_{1}d\xi_{2}d\mu_2d\sigma_{2}.\label{3.0121}
\end{eqnarray}
Thus, we have that ${\rm J}_{3}$ can be bounded by
\begin{eqnarray}
 C\sum\limits_{m_{1},m_{2}>0,\>m}\sum\limits_{j_{1},j_{2}\geq0,\>0<j\leq \frac{(2\alpha-1)m_{1}}{2}}
\left(\int_{D_{j,j_{1},j_{2},m,m_{1},m_{2}}^{(11)}}
2^{jb^{\prime}-(j_{1}+j_{2})b+{\rm min}\left\{m,m_{1},m_{2}\right\}}QdV
\right).\label{3.0122}
\end{eqnarray}
In this case, we consider (\ref{3.031}), (\ref{3.032}),
respectively.

\noindent When (\ref{3.031}) is valid,  we make the change of variables (\ref{3.032}),
thus the Jacobian determinant equals
\begin{eqnarray}
\frac{\partial(u,v,w,\mu_{2})}{\partial(\xi_{1},\xi_{2},\mu_{1},\mu_{2})}
=(\alpha+1)(|\xi_{1}|^{\alpha}-|\xi_{2}|^{\alpha})
-\left[\left(\frac{\mu_{1}}{\xi_{1}}\right)^{2}-\left(\frac{\mu_{2}}{\xi_{2}}\right)^{2}\right].\label{3.0123}
\end{eqnarray}
We assume that $D_{j,j_{1},j_{2},m,m_{1},m_{2}}^{(12)}$
 is the image of the subset of all points
$$(\xi_{1},\mu_{1},\sigma_{1},\xi_{2},\mu_{2},\sigma_{2})\in D_{j,j_{1},j_{2},m,m_{1},m_{2}}^{(11)},$$
which satisfies (\ref{3.031}) under the transformation as in
(\ref{3.033}). Combining (\ref{3.0123}) with (\ref{3.031}), we have that
\begin{eqnarray}
\left|\frac{\partial(u,v,w,\mu_{2})}{\partial(\xi_{1},\xi_{2},\mu_{1},\mu_{2})}\right|>2^{j+\frac{\alpha-1}{2}m_{1}}\label{3.0124}
\end{eqnarray}
Let $G_{5}(u,v,w,\mu_{2},\sigma_{1},\sigma_{2})$ be $\eta_{m}(\xi)\eta_{j}(\sigma)\prod\limits_{k=1}^{2}g_{m_{k},j_{k}}$
 under the change of the variables (\ref{3.033}) and
\begin{eqnarray}
M_{9}=F(u,v,w)G_{5}(u,v,w,\mu_{2},\sigma_{1},\sigma_{2}),
dV^{(1)}=dudvdwd\mu_{2}d\sigma_{1}d\sigma_{2}\label{3.0125}.
\end{eqnarray}
Thus, (\ref{3.0122}) can be controlled by
\begin{eqnarray}
&&\hspace{-2cm}C\sum\limits_{m_{1},m_{2}\geq0,\>m}\sum\limits_{j_{1},j_{2}\geq0,\>0<j\leq
\frac{(2\alpha-1)m_{1}}{2}}\int_{D_{j,j_{1},j_{2},m,m_{1},m_{2}}^{(12)}}
2^{jb^{\prime}-(j_{1}+j_{2})b+{\rm min}\left\{m,m_{1},m_{2}\right\}}\frac{M_{5}dV^{(1)}}
{\left|\frac{\partial(u,v,w,\mu_{2})}{\partial(\xi_{1},\xi_{2},\mu_{1},\mu_{2})}\right|}\label{3.0126}.
\end{eqnarray}
Inspired by \cite{CKS-GAFA, LX},
we define
\begin{eqnarray}
&&f(\mu):=\sigma_{1}+\sigma_{2}-
(\xi|\xi|^{\alpha}-\xi_{1}|\xi_{1}|^{\alpha}-\xi_{2}|\xi_{2}|^{\alpha})+\frac{\xi_{1}\xi_{2}}{\xi}
\left[\frac{\mu_{1}}{\xi_{1}}-\frac{\mu}{\xi_{2}}\right]^{2}\label{3.0127}.
\end{eqnarray}
For fixed $\sigma_{1},\sigma_{2},\xi_{1},\xi_{2},\mu_{1}$, from (\ref{3.0127}),  we have that
\begin{eqnarray}
&&\hspace{-0.5cm}|f(\mu_{2})|=|\sigma_{1}+\sigma_{2}+\phi(\xi_{1},\mu_{1})+\phi(\xi_{2},\mu_{2})
-\phi(\xi,\mu)|=|\tau-\phi(\xi,\mu)|\sim 2^{j},\label{3.0128}\\&&|f^{\prime}(\mu_{2})|\sim \left|\frac{\xi_{1}}{\xi}\right|
\left|\frac{\mu_{1}}{\xi_{1}}-\frac{\mu_{2}}{\xi_{2}}\right|\sim 2^{\frac{\alpha m_{1}}{2}} .\label{3.0129}
\end{eqnarray}
Combining (\ref{3.0128}), (\ref{3.0129}) with Lemma 2.6, for fixed $\sigma_{1},\sigma_{2},\xi_{1},\xi_{2},\mu_{1}$, we have that
the Lebesgue measure of $\mu_{2}$ can be controlled by $C2^{j-\frac{\alpha m_{1}}{2}}$.
By using the Cauchy-Schwartz inequality with respect to $\mu_{2}$
 and the inverse change of variables related to  (\ref{3.033}) and the Cauchy-Schwartz inequality with respect to $u,v,w$
  and the Cauchy-Schwartz inequality with respect to $\sigma_{1}$
   and  $\sigma_{2}$, we have that (\ref{3.0126})  can be  bounded by
\begin{eqnarray}
&&\hspace{-0.5cm}C\sum\int_{D_{j,j_{1},j_{2},m,m_{1},m_{2}}^{(12)}}
2^{jb^{\prime}-(j_{1}+j_{2})b+m}\frac{M_{9}dV^{(1)}}
{\left|\frac{\partial(u,v,w,\mu_{2})}{\partial(\xi_{1},\xi_{2},\mu_{1},\mu_{2})}\right|}\nonumber\\
&&\leq C\sum
2^{2j\epsilon-(j_{1}+j_{2})b+m-\frac{\alpha m_{1}}{4}}\int F(u,v,w)
\left(\int\frac{G_{5}^{2}(u,v,w,\mu_{2},\sigma_{1},\sigma_{2})}{
\left|\frac{\partial(u,v,w,\mu_{2})}{\partial(\xi_{1},\xi_{2},\mu_{1},\mu_{2})}\right|^{2}}
d\mu_{2}\right)^{\frac{1}{2}}dV^{(2)}\nonumber\\
&&\leq C\sum2^{jb^{\prime}-(j_{1}+j_{2})b
-\frac{(2\alpha-1)m_{1}}{4}+m}\|F\|_{L^{2}}\int
\left(\int\frac{G_{5}^{2}(u,v,w,\mu_{2},\sigma_{1},\sigma_{2})}{
\left|\frac{\partial(u,v,w,\mu_{2})}{\partial(\xi_{1},\xi_{2},\mu_{1},\mu_{2})}\right|}
dV^{(3)}\right)^{\frac{1}{2}}
d\sigma_{1}d\sigma_{2}\nonumber\\
&&\leq C\sum
2^{jb^{\prime}-(j_{1}+j_{2})b
-\frac{(2\alpha-1)m_{1}}{4}+m}\|F\|_{L^{2}}\int
\left(\int \prod\limits_{k=1}^{2}g_{m_{k},j_{k}}^{2}d\xi_{1}d\mu_{1}d\xi_{2}d\mu_{2}\right)^{\frac{1}{2}}
d\sigma_{1}d\sigma_{2}\nonumber\\&&
\leq C\sum
2^{jb^{\prime}-(j_{1}+j_{2})\epsilon-\frac{(2\alpha-1)m_{1}}{4}+m}
\|F\|_{L^{2}}
\left(\int \prod\limits_{k=1}^{2}g_{m_{k},j_{k}}^{2}dV\right)^{\frac{1}{2}}\nonumber\\
&&\leq C\sum
2^{jb^{\prime}-(j_{1}+j_{2})\epsilon-\frac{(2\alpha-1)m_{1}}{4}+m}
\|F\|_{L^{2}}
\left(\prod\limits_{j=1}^{2}\|F_{j}\|_{L^{2}}\right)\nonumber\\
&&
\leq C\sum\limits_{m_{1},m_{2}>0}2^{-\frac{(2\alpha-5)m_{1}}{4}}
\|F\|_{L^{2}}
\left(\prod\limits_{j=1}^{2}\|F_{j}\|_{L^{2}}\right)\nonumber\\
&&\leq C\sum\limits_{m_{1},m_{2}>0}2^{-\frac{(2\alpha-5)m_{1}}{4}+m_{1}\epsilon}2^{-m_{2}\epsilon}
\|F\|_{L^{2}}
\left(\prod\limits_{j=1}^{2}\|F_{j}\|_{L^{2}}\right)\nonumber\\
&&\leq C
N^{\frac{5-2\alpha}{4}+\epsilon}\|F\|_{L^{2}}
\left(\prod\limits_{j=1}^{2}\|F_{j}\|_{L^{2}}\right),\label{3.0130}
\end{eqnarray}
where $\sum=\sum\limits_{j_{1},j_{2}\geq0,\>0<j\leq \frac{(2\alpha-1)m_{1}}{2}}.$

\noindent When  (\ref{3.032}) is  valid, we make the change of variables (\ref{3.042}),
thus the Jacobian determinant equals
\begin{eqnarray}
\frac{\partial(u,v,w,\xi_{1})}{\partial(\xi_{1},\xi_{2},\mu_{1},\mu_{2})}
=2\left[\frac{\mu_{1}}{\xi_{1}}-\frac{\mu_{2}}{\xi_{2}}\right].\label{3.0131}
\end{eqnarray}
We assume that $D_{j,j_{1},j_{2},m,m_{1},m_{2}}^{(13)}$ is the image of the subset of all points
$$(\xi_{1},\mu_{1},\sigma_{1},\xi_{2},\mu_{2},\sigma_{2})\in D_{j,j_{1},j_{2},m,m_{1},m_{2}}^{(11)},$$
 which satisfies (\ref{3.032}) under the transformation
(\ref{3.042}). Combining (\ref{3.0116}) with (\ref{3.0131}), we have that
\begin{eqnarray}
\left|\frac{\partial(u,v,w,\xi_{1})}{\partial(\xi_{1},\xi_{2},\mu_{1},\mu_{2})}\right|\sim 2^{m+\frac{(\alpha-2)m_{1}}{2}}.\label{3.0132}
\end{eqnarray}
Let $H_{5}(u,v,w,\xi_{1},\sigma_{1},\sigma_{2})$
be $\eta_{m}(\xi)\eta_{j}(\sigma)\prod\limits_{k=1}^{2}g_{m_{k},j_{k}}$
under the change of the variables  (\ref{3.042}) and
\begin{eqnarray}
M_{10}=F(u,v,w)H_{5}(u,v,w,\xi_{1},\sigma_{1},\sigma_{2}), dV^{(4)}=dudvdwd\xi_{1}d\sigma_{1}d\sigma_{2}.\label{3.0133}
\end{eqnarray}
Thus, (\ref{3.0122}) can be controlled by
\begin{eqnarray}
&&\hspace{-2cm}C\sum\limits_{m_{1},m_{2}>0,\>m}\sum\limits_{j_{1},j_{2}\geq0,
\>0<j\leq \frac{(2\alpha-1)m_{1}}{2}}\int_{D_{j,j_{1},j_{2},m,m_{1},m_{2}}^{(13)}}
2^{jb^{\prime}-(j_{1}+j_{2})b+{\rm min}\left\{m,m_{1},m_{2}\right\}}\frac{M_{10}dV^{(4)}}
{\left|\frac{\partial(u,v,w,\xi_{1})}{\partial(\xi_{1},\xi_{2},\mu_{1},\mu_{2})}\right|}\label{3.0134}.
\end{eqnarray}
We assume that  $h(\xi)$ is defined as in (\ref{3.046}), from (\ref{3.032}), for fixed $\xi_{2},\mu_{1},\mu_{2},$ we have that
\begin{eqnarray}
&&|h^{\prime}(\xi_{1})|=\left|\alpha(\alpha+1)\xi_{1}|\xi_{1}|^{\alpha-2}+2
\left(\frac{\mu_{1}}{\xi_{1}}\right)^{2}\xi_{1}\right|\geq \alpha(\alpha+1)|\xi_{1}|^{\alpha-1}\geq C2^{(\alpha-1)m_{1}},\nonumber\\&&
|h(\xi_{1})|\leq C2^{j+\frac{(\alpha-1)m_{1}}{2}}\label{3.0135},
\end{eqnarray}
for fixed $\xi_{2},\mu_{1},\mu_{2}$, combining (\ref{3.0136}) with Lemma 2.6,
we have that
the Lebesgue measure of $\xi_{1}$ can be controlled by
 $C2^{j-\frac{(\alpha-1)m_{1}}{2}}$.
By using the Cauchy-Schwartz inequality with respect to $\xi_{1}$
 and the inverse change of variables related to  (\ref{3.042}) and the Cauchy-Schwartz inequality with respect to $u,v,w$
  and the Cauchy-Schwartz inequality with respect to $\sigma_{1}$
   and $\sigma_{2}$, we have that (\ref{3.0134}) can be bounded by
\begin{eqnarray}
&&\hspace{-0.8cm}C\sum\int_{D_{j,j_{1},j_{2},m,m_{1},m_{2}}^{(13)}}
2^{jb^{\prime}-(j_{1}+j_{2})b+m}\frac{M_{10}dV^{(4)}}
{\left|\frac{\partial(u,v,w,\xi_{1})}{\partial(\xi_{1},\xi_{2},\mu_{1},\mu_{2})}\right|}\nonumber\\
&&\hspace{-0.8cm}\leq C\sum2^{2j\epsilon-(j_{1}+j_{2})b
+m-\frac{(\alpha-1)m_{1}}{4}}
\int F(u,v,w)
\left(\int\frac{H_{5}^{2}(u,v,w,\xi_{1},\sigma_{1},\sigma_{2})}{
\left|\frac{\partial(u,v,w,\xi_{1})}{\partial(\xi_{1},\xi_{2},\mu_{1},\mu_{2})}\right|^{2}}d\xi_{1}\right)^{\frac{1}{2}}dV^{(2)}\nonumber\\
&&\hspace{-0.8cm}\leq C\sum2^{2j\epsilon-(j_{1}+j_{2})b +m-\frac{(2\alpha-3)m_{1}}{4}-\frac{m}{2}}\|F\|_{L^{2}}\int
\left(\int\frac{H_{5}^{2}(u,v,w,\xi_{1},\sigma_{1},\sigma_{2})}{
\left|\frac{\partial(u,v,w,\xi_{1})}{\partial(\xi_{1},\xi_{2},\mu_{1},\mu_{2})}\right|}
dV^{(5)}\right)^{\frac{1}{2}}
d\sigma_{1}d\sigma_{2}\nonumber\\
&&\hspace{-0.8cm}\leq C\sum2^{2j\epsilon-(j_{1}+j_{2})b+\frac{(2\alpha-3)m_{1}}{4}+\frac{m}{2}}\|F\|_{L^{2}}\int
\left(\int \prod\limits_{k=1}^{2}g_{m_{k},j_{k}}^{2}d\xi_{1}d\mu_{1}d\xi_{2}d\mu_{2}\right)^{\frac{1}{2}}
d\sigma_{1}d\sigma_{2}\nonumber\\
&&\hspace{-0.8cm}\leq C\sum2^{2j\epsilon-(j_{1}+j_{2})
\epsilon-\frac{(2\alpha-3)m_{1}}{4}+\frac{m}{2}}\|F\|_{L^{2}}
\left(\int \prod\limits_{k=1}^{2}g_{m_{k},j_{k}}^{2}dV\right)^{\frac{1}{2}}\nonumber\\
&&\hspace{-0.8cm}\leq C\sum\limits_{m_{1},m_{2}>0}N^{-\frac{2\alpha-5}{4}+(2\alpha-1)m_{1}\epsilon}
\|F\|_{L^{2}}\left(\prod\limits_{j=1}^{2}\|F_{j}\|_{L^{2}}\right)\nonumber\\
&&\hspace{-0.8cm}\leq C\sum\limits_{m_{1},m_{2}>0}N^{-\frac{2\alpha-5}{4}+2\alpha m_{1}\epsilon}2^{-m_{2}\epsilon}
\|F\|_{L^{2}}\left(\prod\limits_{j=1}^{2}\|F_{j}\|_{L^{2}}\right)\nonumber\\
&&\hspace{-0.8cm}\leq
 CN^{-\frac{2\alpha-5}{4}+2\alpha\epsilon}\|F\|_{L^{2}}\left(\prod\limits_{j=1}^{2}\|F_{j}\|_{L^{2}}\right),\label{3.0136}
\end{eqnarray}
where $\sum=\sum\limits_{m_{1},m_{2}>0,\>m}\sum\limits_{j_{1},j_{2}\geq0,
\>0<j\leq \frac{(2\alpha-1)m_{1}}{2}}$ and $dV^{(5)}:=dudvdwd\xi_{1}$.

\noindent (4) Region $A_{4}$.
\noindent  When (\ref{3.017}) is valid, from Lemma 2.8, we have that
\begin{eqnarray}
\left|\frac{\mu_{1}}{\xi_{1}}-\frac{\mu_{2}}{\xi_{2}}\right|\sim2^{m+\frac{(\alpha-2)m_{1}}{2}}.\label{3.0137}
\end{eqnarray}
In this case, we have that
\begin{eqnarray}
M(\xi_{1},\xi_{2})\leq C\left(\frac{\prod\limits_{j=1}^{2}|\xi_{j}|}{N^{2}}\right)^{-s_{1}}.\label{3.0138}
\end{eqnarray}
In this case, we consider $|\sigma|\geq |\xi_{1}|^{\frac{2\alpha+1}{2}}$ and $|\sigma|<|\xi_{1}|^{\frac{2\alpha+1}{2}}$, respectively.

\noindent
When $|\sigma|\geq |\xi_{1}|^{\frac{2\alpha+1}{2}}$, since $s\geq-\frac{\alpha-1}{8}+4\alpha\epsilon$,  we have that
\begin{eqnarray}
&&K_{2}(\xi_{1},\mu_{1},\tau_{1},\xi,\mu,\tau)\leq CN^{2s}\frac{|\xi|\prod\limits_{j=1}^{2}|\xi_{j}|^{-s}}
{\langle\sigma\rangle^{-b^{\prime}}
\prod\limits_{j=1}^{2}\langle\sigma_{j}\rangle^{b}}\leq CN^{2s}\frac{|\xi_{1}|^{-\frac{2\alpha+1}{4}+(2\alpha+1)\epsilon-s}|\xi_{2}|^{-s}|\xi|}
{\prod\limits_{j=1}^{2}\langle\sigma_{j}\rangle^{b}}
\nonumber\\
&&\leq CN^{2s}\frac{|\xi_{1}|^{-\frac{2\alpha-3}{4}+(2\alpha+1)\epsilon-2s}}
{\prod\limits_{j=1}^{2}\langle\sigma_{j}\rangle^{b}}\leq C N^{2s}\frac{|\xi_{1}|^{-\frac{3\alpha-5}{4}-2s+(2\alpha+1)\epsilon}|\xi_{1}|^{\frac{\alpha}{4}}|\xi_{2}|^{-\frac{1}{2}}}
{\prod\limits_{j=1}^{2}\langle\sigma_{j}\rangle^{b}}\nonumber\\&&\leq C N^{-\frac{3\alpha-5}{4}+(2\alpha+1)\epsilon}
\frac{|\xi_{1}|^{\frac{\alpha}{4}}|\xi_{2}|^{-\frac{1}{2}}}
{\prod\limits_{j=1}^{2}\langle\sigma_{j}\rangle^{b}}
.\label{3.0139}
\end{eqnarray}
Thus, combining (\ref{2.017}) with (\ref{3.0139}),  we have that
\begin{eqnarray*}
{\rm |J_{4}|}\leq CN^{-\frac{3\alpha-5}{4}+(2\alpha+1)\epsilon}\|F\|_{L_{\tau\xi\mu}^{2}}
\left(\prod\limits_{j=1}^{2}\|F_{j}\|_{L_{\tau\xi\mu}^{2}}\right).
\end{eqnarray*}
We consider case $|\sigma|<|\xi_{1}|^{\frac{2\alpha+1}{2}}$. We consider (\ref{3.031}), (\ref{3.032}), respectively.

\noindent
We dyadically decompose with respect to
\begin{eqnarray*}
\langle\sigma\rangle\sim 2^{j},\langle\sigma_{1}\rangle\sim 2^{j_{1}},
\langle\sigma_{2}\rangle\sim 2^{j_{2}},|\xi|\sim 2^{m},|\xi_{1}|\sim 2^{m_{1}},|\xi_{2}|\sim 2^{m_{2}}.
\end{eqnarray*}
Let  $D_{j,j_{1},j_{2},m,m_{1},m_{2}}^{(14)}$ be the image of set of all points
$(\xi_{1},\mu_{1},\tau_{1},\xi,\mu,\tau)\in D^{*}$ satisfying
\begin{eqnarray}
&&|\xi_{1}|>\frac{N}{2},|\xi_{1}|\geq |\xi_{2}|,|\xi_{2}|>N,
|\xi|\sim 2^{m},|\xi_{1}|\sim 2^{m_{1}},|\xi_{2}|\sim 2^{m_{2}},\nonumber\\
&&\langle\sigma_{1}\rangle\sim 2^{j_{1}},
\langle\sigma_{2}\rangle\sim 2^{j_{2}},\langle\sigma\rangle\sim 2^{j}\leq 2^{\frac{(2\alpha+1)m_{1}}{2}},\label{3.0141}
\end{eqnarray}
under the transformation $(\xi_{1},\mu_{1},\tau_{1},\xi_{2},\mu_{2},\tau_{2})\longrightarrow
(\xi_{1},\mu_{1},\sigma_{1},\xi_{2},\mu_{2},\sigma_{2})$.
Thus, we have that ${\rm J}_{4}$ can be bounded by
\begin{eqnarray}
\hspace{-1cm} CN^{2s}\sum\limits_{m_{1},m_{2}>0,\>m}\sum\limits_{j_{1},j_{2}\geq0,
\>0<j\leq \frac{(2\alpha+1)m_{1}}{2}}\
\int_{D_{j,j_{1},j_{2},m,m_{1},m_{2}}^{(14)}}
2^{jb^{\prime}-(j_{1}+j_{2})b+(m_{1}+m_{2})|s|+m}QdV
,\label{3.0141}
\end{eqnarray}
where $Q,dV$ are defined as in (\ref{3.0121}).

 \noindent When (\ref{3.031}) is valid, we make the change of variables  (\ref{3.033}).

\noindent Thus the Jacobian determinant equals
\begin{eqnarray}
\frac{\partial(u,v,w,\mu_{2})}{\partial(\xi_{1},\xi_{2},\mu_{1},\mu_{2})}=(\alpha+1)(|\xi_{1}|^{\alpha}-|\xi_{2}|^{\alpha})
-\left[\left(\frac{\mu_{1}}{\xi_{1}}\right)^{2}-\left(\frac{\mu_{2}}{\xi_{2}}\right)^{2}\right]\label{3.0142}.
\end{eqnarray}
We assume that $D_{j,j_{1},j_{2},m,m_{1},m_{2}}^{(14)}$ is the image of the subset of all points
$$(\xi_{1},\mu_{1},\sigma_{1},\xi_{2},\mu_{2},\sigma_{2})\in D_{j,j_{1},j_{2},m,m_{1},m_{2}}^{(13)},$$
which satisfies (\ref{3.031}) under the transformation
(\ref{3.033}). Combining (\ref{3.0142}) with (\ref{3.031}), we have that
\begin{eqnarray}
\left|\frac{\partial(u,v,w,\mu_{2})}{\partial(\xi_{1},\xi_{2},\mu_{1},\mu_{2})}\right|>2^{j+\frac{m}{2}+\frac{(\alpha-2)m_{1}}{2}}\label{3.0143}.
\end{eqnarray}
Let $G_{6}(u,v,w,\mu_{2},\sigma_{1},\sigma_{2})$ be $\eta_{m}(\xi)\eta_{j}(\sigma)\prod\limits_{k=1}^{2}g_{m_{k},j_{k}}$
 under the change of the variables (\ref{3.033}) and
\begin{eqnarray}
&&M_{11}=F(u,v,w)G_{6}(u,v,w,\mu_{2},\sigma_{1},\sigma_{2}),\sum=\sum\limits_{m_{1},m_{2}>0,\>m}\sum\limits_{j_{1},j_{2}\geq0,\>0<j\leq \frac{(2\alpha+1)m_{1}}{2}},\nonumber\\&&
dV^{(1)}=dudvdwd\mu_{2}d\sigma_{1}d\sigma_{2}\label{3.0144}.
\end{eqnarray}
Thus, (\ref{3.0141}) can be controlled by
\begin{eqnarray}
&&\hspace{-2cm}CN^{2s}\sum
\int_{D_{j,j_{1},j_{2},m,m_{1},m_{2}}^{(14)}}
2^{jb^{\prime}-(j_{1}+j_{2})b-(m_{1}+m_{2})s+m}\frac{M_{11}dV^{(1)}}
{\left|\frac{\partial(u,v,w,\mu_{2})}{\partial(\xi_{1},\xi_{2},\mu_{1},\mu_{2})}\right|}\label{3.0145}.
\end{eqnarray}
We assume that $f(\mu)$ is defined as in (\ref{3.038}), for fixed $\sigma_{1},\sigma_{2},\xi_{1},\xi_{2},\mu_{2}$, we have that
\begin{eqnarray}
&&\hspace{-0.5cm}|f(\mu_{2})|=|\sigma_{1}+\sigma_{2}+\phi(\xi_{1},\mu_{1})+\phi(\xi_{2},\mu_{2})
-\phi(\xi,\mu)|=|\tau-\phi(\xi,\mu)|\sim 2^{j},\label{3.0146}\\
&&|f^{\prime}(\mu_{2})|\sim \left|\frac{\xi_{1}}{\xi}\right|
\left|\frac{\mu_{1}}{\xi_{1}}-\frac{\mu_{2}}{\xi_{2}}\right|\sim 2^{\frac{\alpha m_{1}}{2}} ,\label{3.0147}
\end{eqnarray}
 for fixed $\sigma_{1},\sigma_{2},\xi_{1},\xi_{2},\mu_{2}$, combining  (\ref{3.0146}),  (\ref{3.0147})  with  Lemma 2.6,
 we have that
the Lebesgue measure of $\mu_{2}$ can be controlled by $C2^{j-\frac{\alpha m_{1}}{2}}$.
By using the Cauchy-Schwartz inequality with respect to $\mu_{2}$
 and the inverse change of variables related to (\ref{3.033})  and the Cauchy-Schwartz inequality with respect to $u,v,w$
  and the Cauchy-Schwartz inequality with respect to $\sigma_{1}$
   and  $\sigma_{2}$, since  $-\frac{2\alpha-5}{8}+2\alpha \epsilon\leq s<0$, we have that (\ref{3.0145}) can be bounded by
\begin{eqnarray}
&&\hspace{-0.5cm}CN^{2s}\sum
\int_{D_{j,j_{1},j_{2},m,m_{1},m_{2}}^{(14)}}
2^{jb^{\prime}-(j_{1}+j_{2})b-(m_{1}+m_{2})s+m}\frac{M_{11}dV^{(1)}}
{\left|\frac{\partial(u,v,w,\mu_{2})}{\partial(\xi_{1},\xi_{2},\mu_{1},\mu_{2})}\right|}\nonumber\\
&&\leq CN^{2s}\sum2^{2j\epsilon-(j_{1}+j_{2})b-2m_{1}s-\frac{\alpha m_{1}}{4}+m}\int F
\left(\int\frac{G_{6}^{2}(u,v,w,\mu_{2},\sigma_{1},\sigma_{2})}{
\left|\frac{\partial(u,v,w,\mu_{2})}{\partial(\xi_{1},\xi_{2},\mu_{1},\mu_{2})}\right|^{2}}
d\mu_{2}\right)^{\frac{1}{2}}dV^{(2)}\nonumber\\
&&\leq CN^{2s}\sum2^{jb^{\prime}-(j_{1}+j_{2})b-m_{1}(2s+\frac{2\alpha-1}{4})+m}\|F\|_{L^{2}}\int
\left(\int\frac{G_{6}^{2}(u,v,w,\mu_{2},\sigma_{1},\sigma_{2})}{
\left|\frac{\partial(u,v,w,\mu_{2})}{\partial(\xi_{1},\xi_{2},\mu_{1},\mu_{2})}\right|}
dV^{(3)}\right)^{\frac{1}{2}}
d\sigma_{1}d\sigma_{2}\nonumber\\
&&\leq CN^{2s}\sum
2^{jb^{\prime}-(j_{1}+j_{2})b-m_{1}(2s+\frac{2\alpha-1}{4})+m}
\|F\|_{L^{2}}\int
\left(\int \prod\limits_{k=1}^{2}g_{m_{k},j_{k}}^{2}d\xi_{1}d\mu_{1}d\xi_{2}d\mu_{2}\right)^{\frac{1}{2}}
d\sigma_{1}d\sigma_{2}\nonumber\\&&
\leq  CN^{2s}\sum
2^{jb^{\prime}-(j_{1}+j_{2})\epsilon-m_{1}(2s+\frac{2\alpha-1}{4})+m}
\|F\|_{L^{2}}
\left(\int \prod\limits_{k=1}^{2}g_{m_{k},j_{k}}^{2}dV\right)^{\frac{1}{2}}\nonumber\\
&&\leq C N^{2s}\sum2^{jb^{\prime}-(j_{1}+j_{2})\epsilon-m_{1}(2s+\frac{2\alpha-1}{4})+m}
\|F\|_{L^{2}}\left(\prod\limits_{j=1}^{2}\|F_{j}\|_{L^{2}}\right)\nonumber\\&&
\leq CN^{2s}\sum_{m_{1},m_{2}>0,m}2^{-m_{1}(2s+\frac{2\alpha-1}{4})+m}\|F\|_{L^{2}}\left(\prod\limits_{j=1}^{2}\|F_{j}\|_{L^{2}}\right)\nonumber\\&&\leq CN^{2s}\sum_{m_{1},m_{2}>0}2^{-m_{1}(2s+\frac{2\alpha-5}{4})}\|F\|_{L^{2}}\left(\prod\limits_{j=1}^{2}\|F_{j}\|_{L^{2}}\right)\nonumber\\
&&\leq C\sum_{m_{1},m_{2}>0}2^{-m_{1}(2s+\frac{2\alpha-5}{4}-\epsilon)}2^{-m_{2}\epsilon}\|F\|_{L^{2}}\left(\prod\limits_{j=1}^{2}\|F_{j}\|_{L^{2}}\right)\nonumber\\
&&\leq CN^{-\frac{2\alpha-5}{4}+ 2\epsilon}\|F\|_{L^{2}}\left(\prod\limits_{j=1}^{2}\|F_{j}\|_{L^{2}}\right)
\label{3.0148}
\end{eqnarray}
Here $\sum=\sum\limits_{m_{1},m_{2}>0,\>m}\sum\limits_{j_{1},j_{2}\geq0,\>0<j\leq \frac{(2\alpha+1)m_{1}}{2}}$.

Now we consider (\ref{3.032}). We make the change of variables  (\ref{3.042}).

\noindent Thus the Jacobian determinant equals
\begin{eqnarray}
\frac{\partial(u,v,w,\mu_{2})}{\partial(\xi_{1},\xi_{2},\mu_{1},\xi_{1})}
=2\left[\frac{\mu_{1}}{\xi_{1}}-\frac{\mu_{2}}{\xi_{2}}\right].\label{3.0149}
\end{eqnarray}
We assume that $D_{j,j_{1},j_{2},m,m_{1},m_{2}}^{(15)}$ is the image of the subset of all points
$$(\xi_{1},\mu_{1},\sigma_{1},\xi_{2},\mu_{2},\sigma_{2})\in D_{j,j_{1},j_{2},m,m_{1},m_{2}}^{(13)},$$
 which satisfies (\ref{3.032}) under the transformation
(\ref{3.042}). Combining (\ref{3.0149}) with (\ref{3.0137}), we have that
\begin{eqnarray}
\left|\frac{\partial(u,v,w,\xi_{1})}{\partial(\xi_{1},\xi_{2},\mu_{1},\mu_{2})}\right|\sim 2^{m+\frac{(\alpha-2)m_{1}}{2}}\label{3.0150}
\end{eqnarray}
 Let $H_{6}(u,v,w,\xi_{1},\sigma_{1},\sigma_{2})$
be $\eta_{m}(\xi)\eta_{j}(\sigma)\prod\limits_{k=1}^{2}f_{m_{k},j_{k}}$
under the change of the variables  (\ref{3.042}) and
\begin{eqnarray}
M_{12}=F(u,v,w)H_{6}(u,v,w,\xi_{1},\sigma_{1},\sigma_{2}), dV^{(4)}=dudvdwd\xi_{1}d\sigma_{1}d\sigma_{2}.\label{3.0151}
\end{eqnarray}
Thus, (\ref{3.0141}) can be controlled by
\begin{eqnarray}
&&\hspace{-2cm}CN^{2s}\sum\limits_{{\rm min}\{ j,j_{1},j_{2},m_{1},m_{2}\}\geq0,\>m}\int_{D_{j,j_{1},j_{2},m,m_{1},m_{2}}^{(15)}}
2^{jb^{\prime}-(j_{1}+j_{2})b-(m_{1}+m_{2})s+m}\frac{M_{12}dV^{(4)}}
{\left|\frac{\partial(u,v,w,\xi_{1})}{\partial(\xi_{1},\xi_{2},\mu_{1},\mu_{2})}\right|}\label{3.0152}.
\end{eqnarray}
We assume that  $h(\xi)$ is defined as in (\ref{3.046}), from (\ref{3.032}), for fixed $\xi_{2},\mu_{1},\mu_{2},$ we have that
\begin{eqnarray}
&&|h^{\prime}(\xi_{1})|=\left|\alpha(\alpha+1)\xi_{1}|\xi_{1}|^{\alpha-2}+2
\left(\frac{\mu_{1}}{\xi_{1}}\right)^{2}\xi_{1}\right|\geq \alpha(\alpha+1)|\xi_{1}|^{\alpha-1}\geq C2^{(\alpha-1)m_{1}},\nonumber\\&&
|h(\xi_{1})|\leq C2^{j+\frac{(\alpha-1)m_{1}}{2}}\label{3.0153},
\end{eqnarray}
for fixed $\xi_{2},\mu_{1},\mu_{2}$,
combining (\ref{3.0153})  with  Lemma 2.6,
we have that
the Lebesgue measure of $\xi_{1}$ can be controlled by
 $C2^{j-\frac{(\alpha-1)m_{1}}{2}}$.
By using the Cauchy-Schwartz inequality with respect to $\xi_{1}$
 and the inverse change of variables related to  (\ref{3.042}) and the Cauchy-Schwartz inequality with respect to $u,v,w$
  and the Cauchy-Schwartz inequality with respect to $\sigma_{1}$
   and $\sigma_{2}$, since $-\frac{2\alpha-5}{8}+2\alpha\epsilon\leq s<0$, we have that (\ref{3.0151}) can be bounded by
\begin{eqnarray}
&&\hspace{-0.5cm}CN^{2s}\sum\int_{D_{j,j_{1},j_{2},m,m_{1},m_{2}}^{(15)}}
2^{jb^{\prime}-(j_{1}+j_{2})b+m}\frac{M_{2}dV^{(4)}}
{\left|\frac{\partial(u,v,w,\xi_{1})}{\partial(\xi_{1},\xi_{2},\mu_{1},\mu_{2})}\right|}\nonumber\\
&&\leq CN^{2s}\sum2^{2j\epsilon-(j_{1}+j_{2})b-m_{1}(s+\frac{\alpha-1}{4})+m}
\int F(u,v,w)\left(\int\frac{H_{6}^{2}(u,v,w,\xi_{1},\sigma_{1},\sigma_{2})}{
\left|\frac{\partial(u,v,w,\xi_{1})}{\partial(\xi_{1},\xi_{2},\mu_{1},\mu_{2})}\right|^{2}}d\xi_{1}\right)^{\frac{1}{2}}dV^{(2)}\nonumber\\
&&\leq CN^{2s}\sum2^{2j\epsilon-(j_{1}+j_{2})b-m_{1}(2s+\frac{2\alpha-1}{4})+m}
\|F\|_{L^{2}}\int
\left(\int\frac{H_{6}^{2}(u,v,w,\xi_{1},\sigma_{1},\sigma_{2})}{
\left|\frac{\partial(u,v,w,\xi_{1})}{\partial(\xi_{1},\xi_{2},\mu_{1},\mu_{2})}\right|}dV^{(5)}\right)^{\frac{1}{2}}
d\sigma_{1}d\sigma_{2}\nonumber\\
&&\leq CN^{2s}\sum2^{2j\epsilon-(j_{1}+j_{2})b-m_{1}(2s+\frac{2\alpha-1}{4})+m}
\|F\|_{L^{2}}\int
\left(\int \prod\limits_{k=1}^{2}g_{m_{k},j_{k}}^{2}d\xi_{1}d\mu_{1}d\xi_{2}d\mu_{2}\right)^{\frac{1}{2}}
d\sigma_{1}d\sigma_{2}\nonumber\\
&&\leq CN^{2s}
\sum2^{2j\epsilon-(j_{1}+j_{2})
\epsilon-m_{1}(2s+\frac{2\alpha-1}{4})+m}\|F\|_{L^{2}}
\left(\int \prod\limits_{k=1}^{2}g_{m_{k},j_{k}}^{2}dV\right)^{\frac{1}{2}}\nonumber\\
&&\leq CN^{2s}
\sum\limits_{m_{1},m_{2}>0,m}2^{-m_{1}(2s+\frac{2\alpha-1}{4}-(2\alpha+1)\epsilon )+m}\|F\|_{L^{2}}
\left(\prod\limits_{j=1}^{2}\|F_{j}\|_{L^{2}}\right)\nonumber\\
&&\leq CN^{2s}
\sum\limits_{m_{1},m_{2}>0}2^{-m_{1}(2s+\frac{2\alpha-5}{4}-(2\alpha+1)\epsilon )}\|F\|_{L^{2}}
\left(\prod\limits_{j=1}^{2}\|F_{j}\|_{L^{2}}\right)\nonumber\\
&&\leq CN^{2s}
\sum\limits_{m_{1},m_{2}>0}2^{-m_{1}(2s+\frac{2\alpha-5}{4}-(2\alpha+2)\epsilon )}2^{-m_{2}\epsilon}\|F\|_{L^{2}}
\left(\prod\limits_{j=1}^{2}\|F_{j}\|_{L^{2}}\right)\nonumber\\
&&\leq CN^{-\frac{2\alpha-5}{4}+3\alpha \epsilon}
\|F\|_{L^{2}}\left(\prod\limits_{j=1}^{2}\|F_{j}\|_{L^{2}}\right).\label{3.0154}
\end{eqnarray}
Here $\sum=\sum\limits_{m_{1},m_{2}>0,\>m}\sum\limits_{j_{1},j_{2}\geq0,\>0<j\leq \frac{(2\alpha+1)m_{1}}{2}}.$

This completes the proof of Lemma 3.2.

\begin{Lemma}\label{Lem3.3}
Let $s\geq-\frac{\alpha-1}{4}+4\alpha\epsilon,s_{2}\geq0$ and $u_{j}\in X_{b}^{s_{1},s_{2}}(j=1,2)$ and $b=\frac{1}{2}+\epsilon$ and $b^{\prime}=-\frac{1}{2}+2\epsilon$.
Then, we have that
\begin{eqnarray}
&&\|\partial_{x}I(u_{1}u_{2})\|_{X_{b^{\prime}}^{0,0}}\leq C
\prod_{j=1}^{2}\|Iu_{j}\|_{X_{b}^{0,0}}.\label{3.0155}
\end{eqnarray}
\end{Lemma}
\noindent{\bf Proof.} To prove (\ref{3.0155}),  by duality, it suffices to  prove that
\begin{eqnarray}
&&\left|\int_{\SR^{3}}\bar{u}\partial_{x}I(u_{1}u_{2})dxdydt\right|\leq
C\|u\|_{X_{-b^{\prime}}^{0,0}}\prod_{j=1}^{2}
\|Iu_{j}\|_{X_{b}^{0,0}}.\label{3.0156}
\end{eqnarray}
for $u\in X_{-b^{\prime}}^{0,0}.$
Let
\begin{eqnarray}
\hspace{-1cm}F(\xi,\mu,\tau)=
\langle \sigma\rangle^{-b^{\prime}}\mathscr{F}u(\xi,\mu,\tau),
F_{j}(\xi_{j},\mu_{j},\tau_{j})=M(\xi_{j})
\langle \sigma_{j}\rangle^{b}
\mathscr{F}u_{j}(\xi_{j},\mu,\tau_{j})(j=1,2),\label{3.0157}
\end{eqnarray}
$D$ is defined as in Lemma 3.1. To obtain (\ref{3.0155}), from (\ref{3.0156}) and (\ref{3.0157}), it suffices to prove that
\begin{eqnarray}
\hspace{-0.5cm}\int_{D}\frac{|\xi|M(\xi)
F(\xi,\mu,\tau)\prod\limits_{j=1}^{2}F_{j}(\xi_{j},\mu_{j},\tau_{j})}{\langle\sigma_{j}\rangle^{-b^{\prime}}
\prod\limits_{j=1}^{2}M(\xi_{j})\langle\sigma_{j}\rangle^{b}}
d\xi_{1}d\mu_{1}d\tau_{1}d\xi d\mu d\tau\leq C
\|F\|_{L_{\tau\xi\mu}^{2}}\prod_{j=1}^{2}\|F_{j}\|_{L_{\tau\xi\mu}^{2}}.\label{3.0158}
\end{eqnarray}
From (2.4) of \cite{IMEJDE}, we have that
\begin{eqnarray}
\frac{M(\xi)}{\prod\limits_{j=1}^{2}M(\xi_{j})}\leq C\frac{\langle\xi\rangle^{s}}
{\prod\limits_{j=1}^{2}\langle\xi_{j}\rangle^{s}}\label{3.0159}.
\end{eqnarray}
By using (\ref{3.0159}), we have that the left hand side of (\ref{3.0158})  can be bounded by
\begin{eqnarray}
&&\int_{D}\frac{|\xi|\langle\xi\rangle^{s}
F(\xi,\mu,\tau)\prod\limits_{j=1}^{2}F_{j}(\xi_{j},\mu_{j},\tau_{j})}{\langle\sigma_{j}\rangle^{-b^{\prime}}
\prod\limits_{j=1}^{2}\langle\xi_{j}\rangle^{s}\langle\sigma_{j}\rangle^{b}}
d\xi_{1}d\mu_{1}d\tau_{1}d\xi d\mu d\tau.\label{3.0160}
\end{eqnarray}
By using (\ref{3.05}),  we have that (\ref{3.0160}) can be bounded by
$ C
\|F\|_{L_{\tau\xi\mu}^{2}}\left(\prod\limits_{j=1}^{2}\|F_{j}\|_{L_{\tau\xi\mu}^{2}}\right).$

This completes the proof of Lemma 3.3.

\bigskip
\bigskip
\noindent {\large\bf 4. Proof of Theorem  1.1}

\setcounter{equation}{0}

 \setcounter{Theorem}{0}

\setcounter{Lemma}{0}

\setcounter{section}{4}

In this section, combining Lemmas 2.2, 3.1 with the fixed point theorem, we present the proof of Theorem 1.1. Let $b,b^{\prime}$ be defined as in Lemma 3.1.

\noindent We define
\begin{eqnarray}
&&\Phi_{1}(u):=\psi(t)W(t)u_{0}+\frac{1}{2}\psi\left(\frac{t}{\tau}\right)\int_{0}^{t}W(t-\tau)\partial_{x}(u^{2})d\tau,\label{4.01}\\
&&B_{1}(0,2C\|u_{0}\|_{H^{s_{1},s_{2}}}):=\left\{u:\|u\|_{X_{b}^{s_{1},s_{2}}}\leq 2C\|u_{0}\|_{H^{s_{1},s_{2}}}\right\}.\label{4.02}
\end{eqnarray}
Combining Lemmas 2.2, 3.1 with (\ref{4.01}), (\ref{4.02}),  we have that
\begin{eqnarray}
&&\left\|\Phi_{1}(u)\right\|_{X_{b}^{s_{1},s_{2}}}\leq \left\|\psi(t)W(t)u_{0}\right\|_{X_{b}^{s_{1},s_{2}}}
+\left\|\frac{1}{2}\psi\left(\frac{t}{\tau}\right)\int_{0}^{t}W(t-\tau)\partial_{x}(u^{2})d\tau\right\|_{X_{b}^{s_{1},s_{2}}}\nonumber\\
&&\leq C\|u_{0}\|_{H^{s_{1},s_{2}}}+CT^{\epsilon}\left\|\partial_{x}(u^{2})\right\|_{X_{b^{\prime}}^{s_{1},s_{2}}}\nonumber\\
&&\leq C\|u_{0}\|_{H^{s_{1},s_{2}}}+CT^{\epsilon}\left\|u\right\|_{X_{b}^{s_{1},s_{2}}}^{2}\leq C\|u_{0}\|_{H^{s_{1},s_{2}}}+4C^{3}T^{\epsilon}\left\|u_{0}\right\|_{H^{s_{1},s_{2}}}^{2}.\label{4.03}
\end{eqnarray}
We choose $T\in (0,1)$ such that
\begin{eqnarray}
T^{\epsilon}=\left[16C^{2}(\|u_{0}\|_{H^{s_{1},s_{2}}}+1)\right]^{-1}.\label{4.04}
\end{eqnarray}
Combining (\ref{4.03}) with (\ref{4.04}),  we have that
\begin{eqnarray}
&&\left\|\Phi_{1}(u)\right\|_{X_{b}^{s_{1},s_{2}}}
\leq C\|u_{0}\|_{H^{s_{1},s_{2}}}+C\left\|u_{0}\right\|_{H^{s_{1},s_{2}}}=2C\left\|u_{0}\right\|_{H^{s_{1},s_{2}}}.\label{4.05}
\end{eqnarray}
Thus, $\Phi_{1}$ maps $B_{1}(0,2C\|u_{0}\|_{H^{s_{1},s_{2}}})$ into $B_{1}(0,2C\|u_{0}\|_{H^{s_{1},s_{2}}})$.
By using Lemmas 2.2, 3.1, (\ref{4.04})-(\ref{4.05}), we have that
\begin{eqnarray}
&&\left\|\Phi_{1}(u_{1})-\Phi_{1}(u_{2})\right\|_{X_{b}^{s_{1},s_{2}}}
\leq C\left\|\frac{1}{2}\psi\left(\frac{t}{\tau}\right)
\int_{0}^{t}W(t-\tau)\partial_{x}(u_{1}^{2}-u_{2}^{2})d\tau\right\|_{X_{b}^{s_{1},s_{2}}}\nonumber\\
&&\leq CT^{\epsilon}\left\|u_{1}-u_{2}\right\|_{X_{b}^{s_{1},s_{2}}}
\left[\left\|u_{1}\right\|_{X_{b}^{s_{1},s_{2}}}+\left\|u_{2}\right\|_{X_{b}^{s_{1},s_{2}}}\right]\nonumber\\
&&\leq 4C^{2}T^{\epsilon}\left\|u_{0}\right\|_{H^{s_{1},s_{2}}}\left\|u_{1}-u_{2}\right\|_{X_{b}^{s_{1},s_{2}}}\leq \frac{1}{2}\left\|u_{1}-u_{2}\right\|_{X_{b}^{s_{1},s_{2}}}.\label{4.06}
\end{eqnarray}
Thus, $\Phi_{1}$ is a contraction in the closed ball $B_{1}(0,2C\|u_{0}\|_{H^{s_{1},s_{2}}})$.
 Consequently, $u$ is the fixed point of $\Phi$ in the closed ball
$B_{1}(0,2C\|u_{0}\|_{H^{s_{1},s_{2}}})$. Then $v:=u|_{[0,T]}\in X_{b}^{s_{1},s_{2}}([0,T])$
is a solution in the interval $[0,T]$ of the Cauchy problem for (\ref{1.01}) with the initial data $u_{0}$.
 For the facts that  uniqueness of the solution and
the solution to the Cauchy problem for (\ref{1.01}) is continuous
with respect to the initial data, we refer the readers  to Theorems II, III of \cite{IMS}.

This ends the proof of Theorem 1.1.
\bigskip
\bigskip

\noindent {\large\bf 5. Proof of Theorem  1.2}

\setcounter{equation}{0}

 \setcounter{Theorem}{0}

\setcounter{Lemma}{0}

\setcounter{section}{5}
In this section, we give the proof of Theorem 1.2. We present the proof of  Lemma 5.1 before giving the proof of Theorem 1.2.

\begin{Lemma}\label{Lemma5.1}
Let $s_{1}>-\frac{\alpha-1}{4}$ and $R:=\frac{1}{8(C+1)^{3}}$, where  $C$  is the largest of those constants which appear in (\ref{2.05})-(\ref{2.06}),  (\ref{3.0155}).
  Then, the Cauchy problem for (\ref{1.01}) locally well-posed for data satisfying $I_{N}u_{0}\in L^{2}(\R^{2})$ with
  \begin{eqnarray}
 \left\|I_{N}u_{0}\right\|_{L^{2}}\leq R.\label{5.01}
  \end{eqnarray}
Moreover, the solution to the Cauchy problem for (\ref{1.01}) exists on a time interval $[0,1]$.
\end{Lemma}
\noindent{\bf Proof.} We define $v:=I_{N}u$. If $u$ is the solution to the Cauchy problem for (\ref{1.01}), then
$v$ satisfies the following equation
\begin{eqnarray}
v_{t}+\partial_{x}^{5}v+\partial_{x}^{-1}
\partial_{y}^{2}v+\frac{1}{2}I_{N}\partial_{x}(I_{N}^{-1}v)^{2}=0.\label{5.02}
\end{eqnarray}
Then $v$ is formally equivalent to the following integral equation
\begin{eqnarray}
v=W(t)I_{N}u_{0}+\frac{1}{2}\int_{0}^{t}W(t-\tau)I_{N}\partial_{x}(I_{N}^{-1}v)^{2}.\label{5.03}
\end{eqnarray}
We define
\begin{eqnarray}
\Phi_{2}(v)=\psi(t)W(t)I_{N}u_{0}+\frac{1}{2}\psi(t)
\int_{0}^{t}W(t-\tau)I_{N}\partial_{x}(I_{N}^{-1}v)^{2}.\label{5.04}
\end{eqnarray}
Let $b,b^{\prime}$ be defined as in Lemmas 3.1-3.3.
By using Lemmas 2.2,  3.3, we have that
\begin{eqnarray}
&&\left\|\Phi_{2}(v)\right\|_{X_{b}^{0,0}}\leq
 \left\|\psi(t)W(t)I_{N}u_{0}\right\|_{X_{b}^{0,0}}
+C\left\|\psi(t)\int_{0}^{t}W(t-\tau)I_{N}\partial_{x}
(I_{N}^{-1}v)^{2}\right\|_{X_{b}^{0,0}}\nonumber\\
&&\leq C\left\|I_{N}u_{0}\right\|_{L^{2}}+C\left\|I_{N}\partial_{x}
(I_{N}^{-1}v)^{2}\right\|_{X_{b^{\prime}}^{0,0}}\nonumber\\
&&\leq C\left\|I_{N}u_{0}\right\|_{L^{2}}+C\left\|I_{N}\partial_{x}
(I_{N}^{-1}v)^{2}\right\|_{X_{b^{\prime}}^{0,0}}\nonumber\\
&&\leq C\left\|I_{N}u_{0}\right\|_{L^{2}}+C\|v\|_{X_{b}^{0,0}}^{2}\leq CR+C\|v\|_{X_{b}^{0,0}}^{2}.\label{5.05}
\end{eqnarray}
We define
\begin{eqnarray}
B_{2}(0,2CR):=\left\{v:\|v\|_{X_{b}^{0,0}}\leq
2CR\right\}.\label{5.06}
\end{eqnarray}
Combining   (\ref{5.05})-(\ref{5.06}) with the definition of  $R$,  we have that
\begin{eqnarray}
&&\left\|\Phi_{2}(v)\right\|_{X_{b}^{0,0}}
\leq CR+4C^{3}R^{2}
=2CR.\label{5.07}
\end{eqnarray}
Thus, $\Phi_{2}$ maps $B_{2}(0,2CR)$ into
$B_{2}(0,2CR)$. We define
\begin{eqnarray}
v_{j}=I_{N}u_{j}(j=1,2),w_{1}=I_{N}^{-1}v_{1}-I_{N}^{-1}v_{2},
w_{2}=I_{N}^{-1}v_{1}+I_{N}^{-1}v_{2}.\label{5.08}
\end{eqnarray}
By using Lemmas 2.2, 3.1, 3.2,  (\ref{5.05})-(\ref{5.06}) and the definition of $R$, we have that
\begin{eqnarray}
&&\left\|\Phi_{2}(v_{1})-\Phi_{2}(v_{2})\right\|_{X_{b}^{0,0}}
\leq C\left\|\psi(t)
\int_{0}^{t}W(t-\tau)\partial_{x}I_{N}\left[(I_{N}^{-1}v_{1})^{2}-(I_{N}^{-1}v_{2})^{2}\right]
d\tau\right\|_{X_{b}^{0,0}}\nonumber\\
&&\leq C\left\|\partial_{x}I_{N}(w_{1}w_{2})\right\|_{X_{b^{\prime}}^{0,0}}\leq C\|v_{1}-v_{2}\|_{X_{b}^{0,0}}
\left[\|v_{1}\|_{X_{b}^{0,0}}+\|v_{2}\|_{X_{b}^{0,0}}\right]\nonumber\\
&&\leq 4C^{2}R^{2}
\|v_{1}-v_{2}\|_{X_{b}^{0,0}}\leq \frac{1}{2}\|v_{1}-v_{2}\|_{X_{\frac{1}{2}+\epsilon}^{0,0}}
.\label{5.09}
\end{eqnarray}
Thus, $\Phi_{2}$ is a contraction in the closed ball $B_{2}(0,2CR)$.
 Consequently, $u$ is the fixed point of $\Phi_{2}$ in the closed ball
$B_{2}(0,2CR)$. Then $v:=u|_{[0,1]}\in X_{b}^{0,0}([0,1])$
is a solution in the interval $[0,1]$ of the Cauchy problem for (\ref{5.03}) with the initial data $I_{N}u_{0}$.
For the uniqueness of the solution, we  refer the readers to Theorem II of \cite{IMS}.
 For the fact that
the solution to the Cauchy problem for (\ref{5.03}) is continuous
with respect to the initial data, we refer the readers to Theorem III of  \cite{IMS}. Since the phase function $\phi(\xi,\mu)$
is singular at $\xi=0$, to  define  the  derivative  of  $W(t)u_{0}$,  the requirement $|\xi|^{-1}\mathscr{F}_{xy}u_{0}(\xi,\mu)\in \mathscr{S}^{'}(\R^{2})$  is necessary.

This ends the proof of Lemma 5.1.

Inspired by  \cite{ILM-CPAA}, we use Lemmas 2.7, 3.2, 5.1 to prove Theorem 1.2.

For $\lambda>0$, we define
\begin{eqnarray}
u_{\lambda}(x,y,t)=\lambda^{\alpha}u
\left(\lambda x,\lambda^{\frac{\alpha}{2}+1}y,\lambda^{\alpha+1} t\right),
 u_{0\lambda}(x,y)=\lambda^{\alpha}
u\left(\lambda x,\lambda^{\frac{\alpha}{2}+1}y\right).\label{5.010}
\end{eqnarray}
Thus, $u_{\lambda}(x,y,t)\in X_{b}^{s_{1},0}([0,\frac{T}{\lambda}])$ is the solution to
\begin{eqnarray}
&&\partial_{t}u_{\lambda}+|D_{x}|^{\alpha}\partial_{x}u_{\lambda}+
\partial_{x}^{-1}\partial_{y}^{2}u_{\lambda}+u_{\lambda}\partial_{x}u_{\lambda}=0,\label{5.011}\\
&&u_{\lambda}(x,y,0)=u_{0\lambda}(x,y),\label{5.012}
\end{eqnarray}
if and only if $u(x,y,t)\in X_{b}^{s,0}([0,T])$ is the solution to the
Cauchy problem for (\ref{1.01}) in $[0,T]$ with the initial data $u_{0}.$
By using a direct computation, for $\lambda \in (0,1),$  we have that
\begin{eqnarray}
\|I_{N}u_{0\lambda}\|_{L^{2}}\leq CN^{-s}\lambda^{\frac{3\alpha-4}{4}+s}\|u_{0}\|_{H^{s,0}}.\label{5.013}
\end{eqnarray}
For $u_{0}\neq0$ and $u_{0}\in H^{s,0}(\R^{2})$, we choose $\lambda,N$ such that
\begin{eqnarray}
\|I_{N}u_{0\lambda}\|_{L^{2}}\leq CN^{-s}\lambda^{\frac{3\alpha-4}{4}+s}
\|u_{0}\|_{H^{s,0}}:=\frac{R}{4}.\label{5.014}
\end{eqnarray}
Then there exist $w_{3}$ which satisfies that  $\|w_{3}\|_{X_{b}^{s,0}}\leq  2CR$ such that $v:=w_{3}\mid_{[0,1]}$
is a solution to the Cauchy problem for (\ref{5.010}) with $u_{0\lambda}$.
Multiplying (\ref{5.010}) by $2I_{N}u_{\lambda}$ and integrating with respect to $x,y$
and integrating by parts with respect to $x$ yield
\begin{eqnarray}
\frac{d}{dt}\int_{\SR^{2}}(I_{N}u)^{2}dxdy+\int_{\SR^{2}}I_{N}u\partial_{x}I_{N}
\left[(u)^{2}\right]dxdy=0.\label{5.015}
\end{eqnarray}
Combining
$
\int_{\SR^{2}}I_{N}u\partial_{x}\left[(I_{N}u)^{2}\right]dxdy=0
$
with (\ref{5.015}),
we have that
\begin{eqnarray}
\frac{d}{dt}\int_{\SR^{2}}(I_{N}u)^{2}dxdy=-
\int_{\SR^{2}}I_{N}u\partial_{x}\left[I_{N}\left((u)^{2}\right)-(I_{N}u)^{2}\right]dxdy.\label{5.016}
\end{eqnarray}
From (\ref{5.016}) and Lemma 2.7,  we have that
\begin{eqnarray}
&&\int_{\SR^{2}}(I_{N}u(x,y,1))^{2}dxdy-\int_{\SR^{2}}(I_{N}u_{0\lambda})^{2}dxdy\nonumber\\&&
=-\int_{0}^{1}\int_{\SR^{2}}I_{N}u_{\lambda}\partial_{x}\left[I_{N}\left((u_{\lambda})^{2}
\right)-(I_{N}u_{\lambda})^{2}\right]dxdydt\nonumber\\
&&=-\int_{\SR}\int_{\SR^{2}}\left(\chi_{[0,1]}(t)I_{N}u_{\lambda}\right)
\left(\chi_{[0,1]}(t)\partial_{x}\left[I_{N}\left((u_{\lambda})^{2}\right)-(I_{N}u_{\lambda})^{2}\right]\right)dxdydt
\nonumber\\
&&\leq C\left\|\chi_{[0,1]}(t)I_{N}u_{\lambda}\right\|_{X_{\frac{1}{2}-\epsilon}^{0,0}}
\left\|\chi_{[0,1]}(t)\partial_{x}\left[I_{N}\left((u_{\lambda})^{2}\right)-(I_{N}u_{\lambda})^{2}\right]
\right\|_{X_{-\frac{1}{2}+\epsilon}^{0,0}}\nonumber\\
&&\leq C\left\|I_{N}u_{\lambda}\right\|_{X_{\frac{1-\epsilon}{2}}^{0,0}}
\left\|\partial_{x}\left[I_{N}\left((u_{\lambda})^{2}\right)-(I_{N}u_{\lambda})^{2}\right]\right\|_{X_{b^{\prime}}^{0,0}}\nonumber\\
&&\leq CN^{-\frac{2\alpha-5}{4}+3\alpha\epsilon}\|I_{N}u_{\lambda}\|_{X_{b}^{0,0}}^{3}.\label{5.017}
\end{eqnarray}
Combining  (\ref{5.014}) with  (\ref{5.017}), we have that
\begin{eqnarray}
\int_{\SR^{2}}(I_{N}u(x,y,1))^{2}dxdy\leq \frac{R^{2}}{4}+ CN^{-\frac{2\alpha-5}{4}+3\alpha\epsilon}
\|I_{N}u_{\lambda}\|_{X_{b}^{0,0}}^{3}\leq \frac{R^{2}}{4}+8C^{4}N^{-\frac{2\alpha-5}{4}+3\alpha\epsilon}R^{3}.\label{5.018}
\end{eqnarray}
Thus, if we take $N$ sufficiently large such that such that $8C^{4}N^{-\frac{2\alpha-5}{4}+3\alpha\epsilon}R^{3}\leq \frac{3}{4}R^{2},$
then
\begin{eqnarray}
\left[\int_{\SR^{2}}(I_{N}u(x,y,1))^{2}dxdy\right]^{\frac12}\leq R.\label{5.019}
\end{eqnarray}
We consider $I_{N}u(x,y,1)$ as the initial data and repeat the above argument, we obtain that (\ref{5.011})-(\ref{5.012}) possess
a solution in $\R^{2}\times [1,2]$. In this way, we can extend the solution to (\ref{5.011})-(\ref{5.012}) to the time intyerval
$[0,2].$ The above argument can be repeated $L$ steps, where $L$ is the maximal positive integer  such that
\begin{eqnarray}
8C^{4}N^{-\frac{2\alpha-5}{4}+3\alpha\epsilon}R^{3}L\leq \frac{3}{4}R^{2}.\label{5.020}
\end{eqnarray}
More precisely, the solution to (\ref{5.011})-(\ref{5.012}) can be extended to the time interval $[0,L]$. Thus, we can prove that (\ref{5.011})-(\ref{5.012})
are globally well-posed in $[0,\frac{T}{\lambda^{\alpha+1}}]$ if  we can choose a number $N$ such that
\begin{eqnarray}
L\geq \frac{T}{\lambda}.\label{5.021}
\end{eqnarray}
From (\ref{5.020}), we know that
\begin{eqnarray}
L\sim N^{\frac{2\alpha-5}{4}-3\alpha\epsilon}.\label{5.022}
\end{eqnarray}
We know  that (\ref{5.021}) is valid provided that the following inequality is valid
\begin{eqnarray}
CN^{\frac{2\alpha-5}{4}-3\alpha \epsilon}\geq \frac{T}{\lambda^{\alpha+1}}\sim CTN^{\frac{-4(\alpha+1)s}{3\alpha-4+4s}}.\label{5.023}
\end{eqnarray}
In fact,  (\ref{5.023}) is valid if
$
N^{\frac{2\alpha-5}{4}}> N^{\frac{-4(\alpha+1)s}{3\alpha-4+4s}}
$
which is equivalent to $-\frac{(\alpha-1)(3\alpha-4)}{4(5\alpha+3)}<s<0.$

This completes the proof of Theorem 1.2.

\bigskip
\bigskip
\noindent {\large\bf 6. Proof of Theorem  1.3}

\setcounter{equation}{0}

 \setcounter{Theorem}{0}

\setcounter{Lemma}{0}

\setcounter{section}{6}
In this section, we give the proof of Theorem 1.3.

From Lemmas 5.1, 3.2, 2.7 and (\ref{5.010})-(\ref{5.016}), we have that
\begin{eqnarray}
&&\int_{\SR^{2}}(I_{N}u(x,y,1))^{2}dxdy-\int_{\SR^{2}}(I_{N}u_{0\lambda})^{2}dxdy\nonumber\\&&
=-\int_{0}^{1}\int_{\SR^{2}}I_{N}u_{\lambda}\partial_{x}\left[I_{N}\left((u_{\lambda})^{2}
\right)-(I_{N}u_{\lambda})^{2}\right]dxdydt\nonumber\\
&&=-\int_{\SR}\int_{\SR^{2}}\left(\chi_{[0,1]}(t)I_{N}u_{\lambda}\right)
\left(\chi_{[0,1]}(t)\partial_{x}\left[I_{N}\left((u_{\lambda})^{2}\right)-(I_{N}u_{\lambda})^{2}\right]\right)dxdydt
\nonumber\\
&&\leq C\left\|\chi_{[0,1]}(t)I_{N}u_{\lambda}\right\|_{X_{\frac{1}{2}-\epsilon}^{0,0}}
\left\|\chi_{[0,1]}(t)\partial_{x}\left[I_{N}\left((u_{\lambda})^{2}\right)-(I_{N}u_{\lambda})^{2}\right]
\right\|_{X_{-\frac{1}{2}+\epsilon}^{0,0}}\nonumber\\
&&\leq C\left\|I_{N}u_{\lambda}\right\|_{X_{\frac{1-\epsilon}{2}}^{0,0}}
\left\|\partial_{x}\left[I_{N}\left((u_{\lambda})^{2}\right)-(I_{N}u_{\lambda})^{2}\right]\right\|_{X_{b^{\prime}}^{0,0}}\nonumber\\
&&\leq CN^{-\frac{\alpha}{4}+3\alpha\epsilon}\|I_{N}u_{\lambda}\|_{X_{b}^{0,0}}^{3}.\label{6.01}
\end{eqnarray}
Combining  (\ref{5.014}) with  (\ref{6.01}), we have that
\begin{eqnarray}
\int_{\SR^{2}}(I_{N}u(x,y,1))^{2}dxdy\leq \frac{R^{2}}{4}+ CN^{-\frac{\alpha}{4}+3\alpha\epsilon}
\|I_{N}u_{\lambda}\|_{X_{b}^{0,0}}^{3}\leq \frac{R^{2}}{4}+8C^{4}N^{-\frac{\alpha}{4}+3\alpha\epsilon}R^{3}.\label{6.02}
\end{eqnarray}
Thus, if we take $N$ sufficiently large such that such that $8C^{4}N^{-\frac{\alpha}{4}+3\alpha\epsilon}R^{3}\leq \frac{3}{4}R^{2},$
then
\begin{eqnarray}
\left[\int_{\SR^{2}}(I_{N}u(x,y,1))^{2}dxdy\right]^{\frac12}\leq R.\label{6.03}
\end{eqnarray}
We consider $I_{N}u(x,y,1)$ as the initial data and repeat the above argument, we obtain that (\ref{5.011})-(\ref{5.012}) possess
a solution in $\R^{2}\times [1,2]$. In this way, we can extend the solution to (\ref{5.011})-(\ref{5.012}) to the time interval
$[0,2].$ The above argument can be repeated $L$ steps, where $L$ is the maximal positive integer  such that
\begin{eqnarray}
8C^{4}N^{-\frac{\alpha}{4}+3\alpha\epsilon}R^{3}L\leq \frac{3}{4}R^{2}.\label{6.04}
\end{eqnarray}
More precisely, the solution to (\ref{5.011})-(\ref{5.012}) can be extended to the time interval $[0,L]$. Thus, we can prove that (\ref{5.011})-(\ref{5.012})
are globally well-posed in $[0,\frac{T}{\lambda^{\alpha+1}}]$ if  we can choose a number $N$ such that
\begin{eqnarray}
L\geq \frac{T}{\lambda}.\label{6.05}
\end{eqnarray}
From (\ref{6.05}), we know that
\begin{eqnarray}
L\sim N^{\frac{\alpha}{4}-3\alpha\epsilon}.\label{6.06}
\end{eqnarray}
We know  that (\ref{6.06}) is valid provided that the following inequality is valid
\begin{eqnarray}
CN^{\frac{\alpha}{4}-3\alpha \epsilon}\geq \frac{T}{\lambda^{\alpha+1}}\sim CTN^{\frac{-4(\alpha+1)s}{3\alpha-4+4s}}.\label{6.07}
\end{eqnarray}
In fact,  (\ref{6.07}) is valid if
$
N^{\frac{\alpha}{4}}> N^{\frac{-4(\alpha+1)s}{3\alpha-4+4s}}
$
which is equivalent to $-\frac{\alpha(3\alpha-4)}{4(5\alpha+4)}<s<0.$

This completes the proof of Theorem 1.3.

\leftline{\large \bf Acknowledgments}

\bigskip

\noindent

 This work is supported by the Natural Science Foundation of China
 under grant numbers 11171116 and 11401180. The first author is also
 supported by   the Young core Teachers program of Henan Normal University and  15A110033.
 We are deeply indebted to professor Hsi-Wei Shih for his valuable
 comments.

\end{document}